\documentclass[11pt]{amsart}

\usepackage{latexsym, amssymb, amsthm,color}
\usepackage[centertags]{amsmath}
\usepackage{pictexwd}
\usepackage[normalem]{ulem}
\definecolor{gray1}{rgb}{0.7,0,0.7}
\newcommand{\todo}[1]{{\leavevmode\unskip\footnotesize\hspace{0.75em}\textcolor{gray1}{[todo: #1]}}}
\definecolor{GREEN}{rgb}{0,1,0}
\definecolor{green4}{rgb}{.1,.5,.1}
\definecolor{blue}{rgb}{0,0,1}
\definecolor{gray}{rgb}{0.5,0.5,0.5}

 \newcommand{\ep}{\end{proof}}

 \newcommand{\sm}{\smallskip}

 \newif\ifpctex

 \newcommand{\eps}{\epsilon}

\newcommand{\ups}{h}
\newcommand{\varsig}{\ell}
  \newtheorem{theorem}{Theorem}
  \newtheorem{definition}{Definition}[section]
  \newtheorem{defi}[definition]{Definition}
  
  \newtheorem{nota}[definition]{Notation}

  \newtheorem{cond}[definition]{Condition}
  \newtheorem{proposition}[definition]{Proposition}
  \newtheorem{lemma}[definition]{Lemma}
  
  \newtheorem{ass}{Assumption}
  \newtheorem{corollary}[definition]{Corollary}
  \newcommand{\beCond}[2]{\Rand{\vspace{0,6cm}\tt #1}\begin{cond}[#2]
  \label{#1}} \theoremstyle{definition}
  \newtheorem{remark}[definition]{Remark}
  
  \numberwithin{equation}{section}
  \newtheoremstyle{step}{3pt}{0pt}{\itshape}{}{\bf}{}{.5em}{}

\theoremstyle{step} \newtheorem{step}{Step}

\newcommand{\E}{\mathbb{E}}
\newcommand{\M}{\mathbb{M}}
\newcommand{\PP}{\mathbb{P}}
\newcommand{\R}{\mathbb{R}}

\newcommand{\N}{\mathbb{N}}

\newcommand{\CM}{\mathcal{M}}

\newcommand{\CS}{\mathcal{S}}

\newcommand{\Rand}[1]{\marginpar{#1}} %\renewcommand{\Rand}[1]{}
\newcommand{\be}[1]{\begin{equation}\label{#1}}
\newcommand{\ee}{\end{equation}}
\newcommand{\bew}[1]{\Rand{\vspace{0,6cm}\tt #1}\begin{equation*}\label{#1}}
\newcommand{\eew}{\end{equation*}}
\newcommand{\bea}[1]{\Rand{\vspace{0,6cm}\tt #1}\begin{eqnarray*}\label{#1}}
\newcommand{\eea}[1]{\end{eqnarray*}}

\newcommand{\beL}[2]{\Rand{\vspace{0,6cm}\tt #1}\begin{lemma}[#2]\label{#1}}
\newcommand{\beD}[2]{\Rand{\vspace{0,6cm}\tt #1}\begin{definition}[#2]\label{#1}}
\newcommand{\beT}[2]{\Rand{\vspace{0,6cm}\tt #1}\begin{theorem}[#2]\label{#1}}
\newcommand{\beP}[2]{\Rand{\vspace{0,6cm}\tt #1}\begin{proposition}[#2]\label{#1}}
\newcommand{\beC}[1]{\Rand{\vspace{0,6cm}\tt #1}\begin{corollary}\label{#1}}
\newcommand{\beR}[1]{\Rand{\vspace{0,6cm}\tt #1}\begin{remark}[#1]\label{#1}}

\newcommand{\Tno}{{_{\displaystyle\Longrightarrow\atop n\to\infty}}}
\newcommand{\tno}{{_{\displaystyle\longrightarrow\atop n\to\infty}}}

\newcommand{\tNo}{{_{\displaystyle\longrightarrow\atop N\to\infty}}}

\newcommand{\smallk}{\mathpzc{k}}

\newcommand{\smallx}{\mathpzc{x}}

\newcommand{\vk}{\underline{\kappa}}

\newcommand{\vx}{\underline{x}}

\newcommand{\ccset}{\mathcal{K}}

\newcommand{\mr}{\underline{\underline{r}}}

\newcommand{\1}{{\bf 1}}

\DeclareMathAlphabet{\mathpzc}{OT1}{pzc}{m}{it}

\begin{document}

\title[Evolving phylogenies]
{{\large Evolving phylogenies of trait-dependent branching with mutation and competition\\[1mm] Part~I: Existence}}

\author{Sandra Kliem}
\address{Sandra Kliem \\ Fakult\"at f\"ur Mathematik\\
Universit\"at Duisburg-Essen \\ Thea-Leymann-Str. 9\\
 45127 Essen \\ Germany}
\email{sandra.kliem@uni-due.de}

%\author{Alexa Manger}
%\address{Alexa Manger \\ Fakult\"at f\"ur Mathematik\\
%Universit\"at Duisburg-Essen, Campus Essen\\  Universit\"atsstra{\ss}e~2\\
% 45132 Essen \\ Germany}
%\email{alexa.manger@uni-due.de}

\author{Anita Winter}
\address{Anita Winter \\ Fakult\"at f\"ur Mathematik\\
Universit\"at Duisburg-Essen \\  Thea-Leymann-Str. 9\\
 45127 Essen \\ Germany}
\email{anita.winter@uni-due.de}

\thanks{Research was supported by the DFG through the SPP Priority Programme 1590.}

\thispagestyle{empty}
\date{\today}



  \begin{abstract}
    \sloppy We propose a type-dependent branching model with mutation and competition for modeling phylogenies of a virus population. The competition kernel depends for any two {virus particles on the particles'} types, the total mass of the population as well as genetic information available through the number of nucleotide substitutions separating the {virus particles}. We consider the evolving phylogenies of this individual based model in the huge population, short reproduction time and {frequent} mutation regime, and show tightness in the state space of marked metric measure spaces. Due to heterogeneity in the {natural branching} rates, the phylogenies are in general not ultra-metric. We therefore develop new techniques for verifying a compact containment condition.
    Finally, we characterize the limit as a solution of a martingale problem.
    \end{abstract}

 \maketitle

{\tiny
 \tableofcontents
 }

% ======================================================================

 \section{Introduction and Motivation}
 \label{S:introduction}
 For many RNA viruses the lack of a proofreading mechanism in the virus' RNA polymerase results in frequent mutation. The high mutation and replication rates cause
viral variability
which provides one of the main obstacles for the host immune response to control an infection.
Therefore it is of
great interest to understand in detail the forces which maintain this diversity.
In the last decade there has been a growing biology literature describing and classifying the epidemiological and phylogenetic pathogen patterns for
particular viruses and bacteria
within one host and within the host population (compare, for example, the excellent survey papers \cite{GPG+04,LOH07} and the references therein).
The patterns range from just a few types at any given time
and one dominating type lasting a long time (for instance, influenza on population level or HIV on intra-host level)
to many types at any given time but none of them lasting very long (for instance, HIV and HCV on population level). A simple branching model with selection which features this dichotomy between long and short lasting dominating types
has been suggested and further studied in \cite{LS09,BK14}.

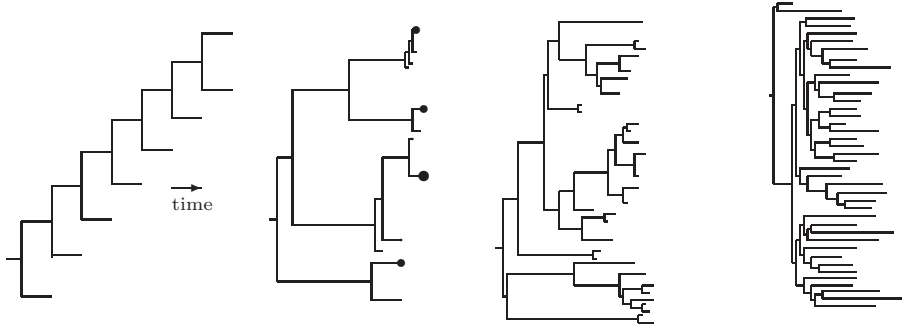
\begin{figure}
\setlength{\unitlength}{.5cm}
\hspace{-12cm}
\begin{picture}(1,7)
%\setcoordinatesystem units <.5cm,.5cm>
%\setplotarea x from -0.5 to 6, y from -0.5 to 6
%%%%%%%%%%%%5 Picture 1 (new) %%%%%%%%%%%%%%%

		\put(.0,-.30){\line(1,0){.4}}
		
		\put(.4,-.30){\line(0,1){1.}}
		\put(.4,-.30){\line(0,-1){1.}}
		\put(.4,0.7){\line(1,0){.8}}
		\put(.4,-1.3){\line(1,0){.8}}
		
		\put(1.2,0.7){\line(0,1){.95}}
		\put(1.2,.7){\line(0,-1){.95}}
		\put(1.2,1.65){\line(1,0){.8}}
		\put(1.2,-.20){\line(1,0){.8}}
		
		\put(2.0,1.65){\line(0,1){0.9}}
		\put(2.0,1.65){\line(0,-1){0.9}}
		\put(2.0,2.55){\line(1,0){0.8}}
		\put(2.0,0.75){\line(1,0){0.8}}
		
		\put(2.8,2.55){\line(0,1){0.85}}
		\put(2.8,2.55){\line(0,-1){0.85}}
		\put(2.8,3.4){\line(1,0){0.8}}
		\put(2.8,1.7){\line(1,0){0.8}}
		
		\put(3.6,3.4){\line(0,1){0.8}}
		\put(3.6,3.4){\line(0,-1){0.8}}
		\put(3.6,4.2){\line(1,0){0.8}}
		\put(3.6,2.6){\line(1,0){0.8}}
		
		\put(4.4,4.2){\line(0,1){0.75}}
		\put(4.4,4.2){\line(0,-1){0.75}}
		\put(4.4,4.95){\line(1,0){0.8}}
		\put(4.4,3.45){\line(1,0){0.8}}
		
		\put(5.2,4.95){\line(0,1){0.7}}
		\put(5.2,4.95){\line(0,-1){0.7}}
		\put(5.2,5.7){\line(1,0){0.8}}
		\put(5.2,4.2){\line(1,0){0.8}}
		
	%	\put(6.0,5.7){\line(0,1){0.65}}
	%	\put(6.0,5.7){\line(0,-1){0.65}}
	%	\put(6.0,6.35){\line(1,0){0.8}}
	%	\put(6.0,5.05){\line(1,0){0.8}}
		
%		\put(6.8,6.65){\line(0,1){0.6}}
%		\put(6.8,6.65){\line(0,-1){0.6}}
%		\put(6.8,7.25){\line(1,0){0.8}}
%		\put(6.8,6.05){\line(1,0){0.8}}

%		\put(5.3,2.35){\vector(1,0){1.4}}
%		\put(5.3,1.55){\small{time}}
				
%%%%%%%%%%%% Picture 1 (alt) %%%%%%%%%%%%%%%%%%%%%

\put(4.4,1.6){\vector(1,0){0.8}}
\put(4.4,1.0){{\tiny time}}
%	\thicklines
% %
%		\put(.0,.0){\line(1,0){.4}}
%		\put(.4,.0){\line(0,1){1.}}
%		\put(.4,.0){\line(0,-1){1.}}
%		\put(.4,1.){\line(1,0){.8}}
%		\put(.4,-1.){\line(1,0){.8}}
%		\put(1.2,1.){\line(0,1){.9}}
%		\put(1.2,1.){\line(0,-1){.9}}
%		\put(1.2,1.9){\line(1,0){.75}}
%		\put(1.2,.1){\line(1,0){.75}}
%		\put(1.95,1.9){\line(0,1){0.8}}
%		\put(1.95,1.9){\line(0,-1){0.8}}
%		\put(1.95,2.7){\line(1,0){0.7}}
%		\put(1.95,1.1){\line(1,0){0.7}}
%		\put(2.65,2.7){\line(0,1){0.7}}
%		\put(2.65,2.7){\line(0,-1){0.7}}
%		\put(2.65,3.4){\line(1,0){0.65}}
%		\put(2.65,2.0){\line(1,0){0.65}}
%		\put(3.3,3.4){\line(0,1){0.6}}
%		\put(3.3,3.4){\line(0,-1){0.6}}
%		\put(3.3,4.0){\line(1,0){0.6}}
%		\put(3.3,2.8){\line(1,0){0.6}}
%		\put(3.9,4.0){\line(0,1){0.5}}
%		\put(3.9,4.0){\line(0,-1){0.5}}
%		\put(3.9,4.5){\line(1,0){0.55}}
%		\put(3.9,3.5){\line(1,0){0.55}}
%		\put(4.45,4.5){\line(0,1){0.4}}
%		\put(4.45,4.5){\line(0,-1){0.4}}
%		\put(4.45,4.9){\line(1,0){0.5}}
%		\put(4.45,4.1){\line(1,0){0.5}}
%		\put(4.95,4.9){\line(0,1){0.3}}
%		\put(4.95,4.9){\line(0,-1){0.3}}
%		\put(4.95,5.2){\line(1,0){0.45}}
%		\put(4.95,4.6){\line(1,0){0.45}}
%		\put(5.4,5.2){\line(0,1){0.2}}
%		\put(5.4,5.2){\line(0,-1){0.2}}
%		\put(5.4,5.4){\line(1,0){0.4}}
%		\put(5.4,5.0){\line(1,0){0.4}}

%%%%%%%%%%% coming down from infinity  %%%%%%%%%%%%%%%%%
		\put(7.0,0.75){\line(1,0){.2}}

		\put(7.2,1.0){\line(0,1){1.4}}
		\put(7.2,1.0){\line(0,-1){1.9}}
		\put(7.2,2.4){\line(1,0){.4}}
		\put(7.2,-0.9){\line(1,0){2.5}}

		\put(7.6,2.4){\line(0,1){1.8}}
		\put(7.6,2.4){\line(0,-1){1.8}}
		\put(9.7,-0.9){\line(0,1){.5}}
		\put(9.7,-0.9){\line(0,-1){.5}}

		\put(7.6,4.2){\line(1,0){1.5}}
		\put(7.6,.6){\line(1,0){2.2}}
		\put(9.7,-0.4){\line(1,0){.8}}
		\put(10.5,-0.4){\circle*{.2}}
		\put(9.7,-1.4){\line(1,0){.8}}

		\put(9.1,4.2){\line(0,1){0.8}}
		\put(9.1,4.2){\line(0,-1){0.8}}
		\put(9.8,.6){\line(0,1){0.7}}
		\put(9.8,.6){\line(0,-1){0.7}}

		\put(9.1,5.0){\line(1,0){1.5}}
		\put(9.1,3.4){\line(1,0){1.7}}
		\put(9.8,1.3){\line(1,0){0.2}}
		\put(9.8,-0.1){\line(1,0){0.2}}

		\put(10.6,5.0){\line(0,1){0.2}}
		\put(10.6,5.0){\line(0,-1){0.2}}
		\put(10.8,3.4){\line(0,1){0.3}}
		\put(10.8,3.4){\line(0,-1){0.3}}
		\put(10.0,1.3){\line(0,1){1.1}}
		\put(10.0,1.3){\line(0,-1){1.1}}

		\put(10.6,5.2){\line(1,0){0.1}}
		\put(10.6,4.8){\line(1,0){0.1}}
		\put(10.8,3.7){\line(1,0){0.3}}
		\put(11.1,3.7){\circle*{.2}}
		\put(10.8,3.1){\line(1,0){0.2}}
		\put(10.0,2.4){\line(1,0){0.7}}
		\put(10.0,.2){\line(1,0){0.5}}
		\put(10.5,.2){\circle*{.1}}

		\put(10.7,5.2){\line(0,1){0.3}}
		\put(10.7,5.2){\line(0,-1){0.3}}
		\put(10.7,5.5){\line(1,0){0.1}}
		\put(10.7,4.9){\line(1,0){0.1}}
		\put(10.8,5.5){\line(0,1){0.3}}
		\put(10.8,5.5){\line(0,-1){0.3}}
		\put(10.8,5.8){\line(1,0){0.1}}
		\put(10.8,5.2){\line(1,0){0.1}}
		\put(10.9,5.8){\circle*{.2}}

		\put(10.7,2.4){\line(0,1){0.5}}
		\put(10.7,2.4){\line(0,-1){0.5}}
		\put(10.7,2.9){\line(1,0){0.1}}
		\put(10.7,1.9){\line(1,0){0.4}}
		\put(10.8,2.9){\circle*{.1}}
		\put(11.1,1.9){\circle*{.3}}

%%%%%%%%%%%%%%%%%% dust free %%%%%%%%%%%%%%%%%%%%%%%

		\put(13.0,.0){\line(1,0){.2}}

		\put(13.2,.0){\line(0,1){1.3}}
		\put(13.2,.0){\line(0,-1){1.3}}

		\put(13.2,1.3){\line(1,0){0.4}}
		\put(13.2,-1.3){\line(1,0){0.1}}

		\put(13.6,1.3){\line(0,1){1.5}}
		\put(13.6,1.3){\line(0,-1){1.5}}
		\put(13.3,-1.3){\line(0,1){.6}}
		\put(13.3,-1.3){\line(0,-1){.6}}

		\put(13.6,2.8){\line(1,0){0.7}}
		\put(13.6,-0.2){\line(1,0){2.0}}
		\put(13.3,-0.7){\line(1,0){1.8}}
		\put(13.3,-1.9){\line(1,0){3.5}}

		\put(14.3,2.8){\line(0,1){1.8}}
		\put(14.3,2.8){\line(0,-1){1.8}}
		\put(15.6,-0.2){\line(0,1){0.1}}
		\put(15.6,-0.2){\line(0,-1){0.1}}
		\put(15.1,-0.7){\line(0,1){0.3}}
		\put(15.1,-0.7){\line(0,-1){0.3}}
		\put(16.8,-1.9){\line(0,1){0.1}}
		\put(16.8,-1.9){\line(0,-1){0.1}}

		\put(14.3,4.6){\line(1,0){0.1}}
		\put(14.3,1.0){\line(1,0){0.4}}
		\put(15.6,-0.1){\line(1,0){0.2}}
		\put(15.6,-0.3){\line(1,0){0.1}}
		\put(15.1,-0.4){\line(1,0){1.6}}
		\put(15.1,-1.0){\line(1,0){1.2}}
		\put(16.8,-1.8){\line(1,0){0.1}}
		\put(16.8,-2.0){\line(1,0){0.4}}

% first branch

		\put(14.4,4.6){\line(0,1){0.9}}
		\put(14.4,4.6){\line(0,-1){0.9}}

		\put(14.4,5.5){\line(1,0){0.3}}
		\put(14.4,3.7){\line(1,0){0.8}}

		\put(14.7,5.5){\line(0,1){0.5}}
		\put(14.7,5.5){\line(0,-1){0.5}}
		\put(15.2,3.7){\line(0,1){0.1}}
		\put(15.2,3.7){\line(0,-1){0.1}}

		\put(14.7,6.0){\line(1,0){2.2}}
		\put(14.7,5.0){\line(1,0){0.7}}

		\put(15.2,3.8){\line(1,0){0.1}}
		\put(15.2,3.6){\line(1,0){0.1}}

		\put(15.4,5.0){\line(0,1){0.4}}
		\put(15.4,5.0){\line(0,-1){0.4}}

		\put(15.4,5.4){\line(1,0){1.3}}
		\put(15.4,4.6){\line(1,0){0.3}}

		\put(15.7,4.6){\line(0,1){0.3}}
		\put(15.7,4.6){\line(0,-1){0.3}}

		\put(15.7,4.9){\line(1,0){0.6}}
		\put(15.7,4.3){\line(1,0){0.1}}

		\put(16.3,4.9){\line(0,1){0.2}}
		\put(16.3,4.9){\line(0,-1){0.2}}
		\put(15.8,4.3){\line(0,1){0.2}}
		\put(15.8,4.3){\line(0,-1){0.2}}

		\put(16.3,5.1){\line(1,0){0.5}}
		\put(16.3,4.7){\line(1,0){0.3}}
		\put(15.8,4.5){\line(1,0){0.7}}
		\put(15.8,4.1){\line(1,0){0.5}}

		\put(16.7,5.4){\line(0,1){0.1}}
		\put(16.7,5.4){\line(0,-1){0.1}}

		\put(16.7,5.5){\line(1,0){0.1}}
		\put(16.7,5.3){\line(1,0){0.3}}

% second branch

		\put(14.7,1.0){\line(0,1){0.5}}
		\put(14.7,1.0){\line(0,-1){0.5}}

		\put(14.7,1.5){\line(1,0){0.4}}
		\put(14.7,0.5){\line(1,0){0.6}}

		\put(15.1,1.5){\line(0,1){0.5}}
		\put(15.1,1.5){\line(0,-1){0.5}}
		\put(15.3,0.5){\line(0,1){0.3}}
		\put(15.3,0.5){\line(0,-1){0.3}}

		\put(15.1,2.0){\line(1,0){0.9}}
		\put(15.1,1.0){\line(1,0){0.5}}
		\put(15.3,0.8){\line(1,0){0.6}}
		\put(15.3,0.2){\line(1,0){0.8}}

		\put(16.0,2.0){\line(0,1){0.6}}
		\put(16.0,2.0){\line(0,-1){0.6}}
		\put(15.9,0.8){\line(0,1){0.1}}
		\put(15.9,0.8){\line(0,-1){0.1}}

		\put(16.0,2.6){\line(1,0){0.2}}
		\put(16.0,1.4){\line(1,0){0.4}}
		\put(15.9,0.9){\line(1,0){0.3}}
		\put(15.9,0.7){\line(1,0){0.1}}

		\put(16.2,2.6){\line(0,1){0.4}}
		\put(16.2,2.6){\line(0,-1){0.4}}
		\put(16.4,1.4){\line(0,1){0.2}}
		\put(16.4,1.4){\line(0,-1){0.2}}

		\put(16.2,3.0){\line(1,0){0.2}}
		\put(16.2,2.2){\line(1,0){0.5}}
		\put(16.4,1.6){\line(1,0){0.4}}
		\put(16.4,1.2){\line(1,0){0.1}}

		\put(16.4,3.0){\line(0,1){0.2}}
		\put(16.4,3.0){\line(0,-1){0.2}}
		\put(16.7,2.2){\line(0,1){0.3}}
		\put(16.7,2.2){\line(0,-1){0.3}}

		\put(16.4,3.2){\line(1,0){0.1}}
		\put(16.4,2.8){\line(1,0){0.4}}
		\put(16.7,2.5){\line(1,0){0.3}}
		\put(16.7,1.9){\line(1,0){0.1}}

		\put(16.5,3.2){\line(0,1){0.1}}
		\put(16.5,3.2){\line(0,-1){0.1}}

		\put(16.5,3.3){\line(1,0){0.3}}
		\put(16.5,3.1){\line(1,0){0.1}}

% third branch

		\put(16.3,-1.0){\line(0,1){0.3}}
		\put(16.3,-1.0){\line(0,-1){0.3}}

		\put(16.3,-0.7){\line(1,0){0.7}}
		\put(16.3,-1.3){\line(1,0){0.1}}

		\put(16.4,-1.3){\line(0,1){0.2}}
		\put(16.4,-1.3){\line(0,-1){0.2}}

		\put(16.4,-1.1){\line(1,0){0.5}}
		\put(16.4,-1.5){\line(1,0){0.2}}

		\put(16.9,-1.1){\line(0,1){0.1}}
		\put(16.9,-1.1){\line(0,-1){0.1}}
		\put(16.6,-1.5){\line(0,1){0.1}}
		\put(16.6,-1.5){\line(0,-1){0.1}}

		\put(16.9,-1.0){\line(1,0){0.3}}
		\put(16.9,-1.2){\line(1,0){0.2}}
		\put(16.6,-1.4){\line(1,0){0.6}}
		\put(16.6,-1.6){\line(1,0){0.3}}

		\put(16.9,-1.6){\line(0,1){0.1}}
		\put(16.9,-1.6){\line(0,-1){0.1}}

		\put(16.9,-1.5){\line(1,0){0.1}}
		\put(16.9,-1.7){\line(1,0){0.2}}

%%%%%%%%%% dusty %%%%%%%%%%%%%%%%%%%%%%%%

		\put(20.3,4.05){\line(1,0){0.1}}
		\put(20.4,4.05){\line(0,1){2.35}}
		\put(20.4,4.05){\line(0,-1){2.35}}
		\put(20.4,6.4){\line(1,0){0.1}}
		\put(20.5,6.4){\line(0,1){0.1}}
		\put(20.5,6.4){\line(0,-1){0.1}}
		\put(20.5,6.5){\line(1,0){0.4}}
		\put(20.5,6.3){\line(1,0){1.7}}
		
		\put(20.4,1.7){\line(1,0){0.5}}
		\put(20.9,1.7){\line(0,1){2.1}}
		\put(20.9,1.7){\line(0,-1){2.1}}

% preparation

		\put(20.9,3.8){\line(1,0){0.1}}

		\put(21.0,3.6){\line(0,1){1.7}}
		\put(21.0,3.6){\line(0,-1){1.7}}

		\put(21.0,5.3){\line(1,0){0.1}}
		\put(21.0,1.9){\line(1,0){0.1}}
		
		\put(20.9,-0.4){\line(1,0){0.1}}

% first one

		\put(21.1,5.3){\line(0,1){0.7}}
		\put(21.1,5.3){\line(0,-1){0.7}}

		\put(21.1,6.0){\line(1,0){0.15}}
		\put(21.25,6.0){\line(0,1){0.1}}
		\put(21.25,6.0){\line(0,-1){0.1}}
		\put(21.25,6.1){\line(1,0){1.3}}
		\put(21.25,5.9){\line(1,0){1.2}}

		\put(21.1,4.6){\line(1,0){0.1}}
		\put(21.2,4.6){\line(0,1){0.9}}
		\put(21.2,4.6){\line(0,-1){0.9}}

		\put(21.2,5.5){\line(1,0){0.1}}
		\put(21.3,5.5){\line(0,1){0.2}}
		\put(21.3,5.5){\line(0,-1){0.2}}
		\put(21.3,5.7){\line(1,0){1.3}}
		\put(21.3,5.3){\line(1,0){0.1}}
		\put(21.4,5.3){\line(0,1){0.2}}
		\put(21.4,5.3){\line(0,-1){0.2}}
		\put(21.4,5.5){\line(1,0){1.1}}
		\put(21.4,5.1){\line(1,0){0.3}}
		\put(21.7,5.1){\line(0,1){0.2}}
		\put(21.7,5.1){\line(0,-1){0.2}}
		\put(21.7,5.3){\line(1,0){1.2}}
		\put(21.7,4.9){\line(1,0){0.2}}
		\put(21.9,4.9){\line(0,1){0.1}}
		\put(21.9,4.9){\line(0,-1){0.1}}
		\put(21.9,5.0){\line(1,0){0.7}}
		\put(21.9,4.8){\line(1,0){1.6}}

% second one

		\put(21.2,3.7){\line(1,0){0.1}}
		\put(21.3,3.7){\line(0,1){0.7}}
		\put(21.3,3.7){\line(0,-1){0.7}}
		\put(21.3,4.4){\line(1,0){0.2}}
		\put(21.3,3.0){\line(1,0){0.1}}
		\put(21.5,4.4){\line(0,1){0.2}}
		\put(21.5,4.4){\line(0,-1){0.2}}
		\put(21.5,4.6){\line(1,0){0.9}}
		\put(21.5,4.2){\line(1,0){0.1}}

		\put(21.6,4.2){\line(0,1){0.2}}
		\put(21.6,4.2){\line(0,-1){0.2}}
		\put(21.6,4.4){\line(1,0){1.6}}
		\put(21.6,4.0){\line(1,0){0.4}}
		\put(22.0,4.0){\line(0,1){0.1}}
		\put(22.0,4.0){\line(0,-1){0.1}}
		\put(22.0,4.1){\line(1,0){1.0}}
		\put(22.0,3.9){\line(1,0){0.7}}

% third one

		\put(21.4,3.0){\line(0,1){0.4}}
		\put(21.4,3.0){\line(0,-1){0.4}}

		\put(21.4,3.4){\line(1,0){0.2}}
		\put(21.6,3.4){\line(0,1){0.2}}
		\put(21.6,3.4){\line(0,-1){0.2}}
		\put(21.6,3.6){\line(1,0){0.3}}
		\put(21.9,3.6){\line(0,1){0.1}}
		\put(21.9,3.6){\line(0,-1){0.1}}
		\put(21.9,3.7){\line(1,0){0.7}}
		\put(21.9,3.5){\line(1,0){0.8}}				
		\put(21.6,3.2){\line(1,0){0.3}}
		\put(21.9,3.2){\line(0,1){0.1}}
		\put(21.9,3.2){\line(0,-1){0.1}}
		\put(21.9,3.3){\line(1,0){0.4}}
		\put(21.9,3.1){\line(1,0){1.3}}

		\put(21.4,2.6){\line(1,0){0.1}}
		\put(21.5,2.6){\line(0,1){0.2}}
		\put(21.5,2.6){\line(0,-1){0.2}}
		\put(21.5,2.8){\line(1,0){0.4}}
		\put(21.5,2.4){\line(1,0){0.5}}
		\put(21.9,2.8){\line(0,1){0.1}}
		\put(21.9,2.8){\line(0,-1){0.1}}
		\put(21.9,2.9){\line(1,0){0.7}}
		\put(21.9,2.7){\line(1,0){0.9}}
		\put(22.0,2.4){\line(0,1){0.1}}
		\put(22.0,2.4){\line(0,-1){0.1}}
		\put(22.0,2.5){\line(1,0){1.2}}
		\put(22.0,2.3){\line(1,0){0.6}}

% fourth one

		\put(21.1,1.9){\line(0,1){0.2}}
		\put(21.1,1.9){\line(0,-1){0.2}}
		\put(21.1,2.1){\line(1,0){1.3}}
		\put(21.1,1.7){\line(1,0){0.2}}
		\put(21.3,1.7){\line(0,1){0.2}}
		\put(21.3,1.7){\line(0,-1){0.2}}
		\put(21.3,1.9){\line(1,0){0.9}}
		\put(21.3,1.5){\line(1,0){0.5}}
		\put(21.8,1.5){\line(0,1){0.2}}
		\put(21.8,1.7){\line(1,0){1.5}}
		\put(21.8,1.3){\line(1,0){0.3}}
		\put(21.8,1.5){\line(0,-1){0.2}}
		\put(22.1,1.3){\line(0,1){0.15}}
		\put(22.1,1.3){\line(0,-1){0.15}}
		\put(22.1,1.45){\line(1,0){1.3}}
		\put(22.1,1.15){\line(1,0){0.2}}
		\put(22.3,1.15){\line(0,1){0.1}}
		\put(22.3,1.15){\line(0,-1){0.1}}
		\put(22.3,1.25){\line(1,0){0.8}}
		\put(22.3,1.05){\line(1,0){0.7}}

% fifth and sixth one

		\put(21.0,-0.4){\line(0,1){0.6}}
		\put(21.0,-0.4){\line(0,-1){0.6}}
		\put(21.0,0.2){\line(1,0){0.1}}
		\put(21.0,-1.0){\line(1,0){0.2}}
		\put(21.1,0.2){\line(0,1){0.4}}
		\put(21.1,0.2){\line(0,-1){0.4}}
		\put(21.1,0.6){\line(1,0){0.1}}
		\put(21.1,-0.2){\line(1,0){0.1}}		
		\put(21.2,-1.0){\line(0,1){0.2}}
		\put(21.2,-1.0){\line(0,-1){0.2}}
		\put(21.2,-0.8){\line(1,0){1.4}}
		\put(21.2,-1.2){\line(1,0){0.2}}

% fifth one

		\put(21.2,0.6){\line(0,1){0.25}}
		\put(21.2,0.6){\line(0,-1){0.25}}
		\put(21.2,0.85){\line(1,0){1.9}}
		\put(21.2,0.35){\line(1,0){0.1}}
		\put(21.3,0.35){\line(0,1){0.15}}
		\put(21.3,0.35){\line(0,-1){0.15}}
		\put(21.3,0.5){\line(1,0){0.1}}
		\put(21.4,0.5){\line(0,1){0.1}}
		\put(21.4,0.5){\line(0,-1){0.1}}
		\put(21.4,0.6){\line(1,0){1.2}}
		\put(21.4,0.4){\line(1,0){2.2}}
		\put(21.3,0.2){\line(1,0){1.9}}

		\put(21.2,-0.2){\line(0,1){0.2}}	
		\put(21.2,0.0){\line(1,0){1.4}}	

		\put(21.2,-0.2){\line(0,-1){0.2}}	
		\put(21.2,-0.4){\line(1,0){0.2}}	
		\put(21.4,-0.4){\line(0,1){0.15}}	
		\put(21.4,-0.4){\line(0,-1){0.15}}	
		\put(21.4,-0.25){\line(1,0){0.8}}	
		\put(21.4,-0.55){\line(1,0){0.2}}	
		\put(21.6,-0.55){\line(0,1){0.1}}	
		\put(21.6,-0.55){\line(0,-1){0.1}}	
		\put(21.6,-0.45){\line(1,0){0.9}}	
		\put(21.6,-0.65){\line(1,0){1.0}}	

% sixth one

		\put(21.4,-1.2){\line(0,1){0.2}}
		\put(21.4,-1.2){\line(0,-1){0.2}}
		\put(21.4,-1.0){\line(1,0){1.1}}
		\put(21.4,-1.4){\line(1,0){0.1}}
		\put(21.5,-1.4){\line(0,1){0.15}}
		\put(21.5,-1.4){\line(0,-1){0.15}}
		\put(21.5,-1.25){\line(1,0){0.2}}
		\put(21.7,-1.25){\line(0,1){0.1}}
		\put(21.7,-1.25){\line(0,-1){0.1}}
		\put(21.7,-1.15){\line(1,0){1.5}}
		\put(21.7,-1.35){\line(1,0){2.1}}
		\put(21.5,-1.55){\line(1,0){1.6}}
	\end{picture}
\vspace{1cm}

{\tiny \caption{shows a) one dominating strain; b) a bounded number of coexisting strains; c) an unbounded number of coexisting strains with proper frequencies;
d)  an unbounded number of coexisting strains without proper frequencies. Note that the branch lengths are according to the number of
substitutions. (The phylogeny pictures are highly inspired by the figures given in \cite{GPG+04}.)}}
\label{Fig:001}
\end{figure}

If one reconstructs the phylogenetic tree based on the number of {nucleotide}
substitutions (compare \cite{BvH99}), then the two situations described above translate into
a trivial tree consisting of just one point representing the clan of the only dominating type versus a radial spread outward from a most recent common ancestor.
Note that in the first case the tree is often drawn with temporal structure including all the fossils yielding a
very skinny tree which consists of a trunk depicting  the ancestry of the dominating type and a few short edges reflecting the high rate of extinction
(compare the phylogeny class a) in Figure~\ref{Fig:001}).  %That means on the one hand the total length of the phylogeny is negligible uniformly in the sample size while on the other hand it grows with %increasing  sample size.
In the latter case of coexisting types, a more elaborate investigation of the patterns suggests to further distinguish between the phylogeny classes b), c) and d) in Figure~\ref{Fig:001}: The main difference between them
 might be  according to statistics such
as the degree of balance of a tree, the number of coexisting strains  or the asymptotic behaviour of the total length of the subtree spanned by a huge sample.\smallskip

 %the total length of the phylogeny might stay bounded, grow slowly (say on a logarithmic scale) or grow fast (on a linear scale) with increasing %sample size. The first case shows in viruses which are known to appear world-wide with a fixed number of serotypes (for example, dengue). The second %case is common in viruses with a strong-level of cross-immunity (for example, measles).
%The last case is reflecting very persistent viruses whose phylogenies consequently must have long edge lengths (HIV, HCV).
%Compare Figure~\ref{Fig:001} for a sketch of the four classes of phylogenetic patterns we have in mind.

In this paper and future work we seek to explore
the phylogenetic patterns which are
primarily affected by natural selection that
arises from cross-immunity, that is, the differential
effect of immune responses on genetically
variable strains.
We propose a parametric model for evolving phylogenies which shows these phylodynamic patterns on the level of phylogenies.
For that purpose
we consider an individual based multi-type branching model with mutation and competition in which {virus particles} branch with a natural branch rate depending on its trait.
Additionally, for each pair of {virus particles}, a competition kernel allows for either one {virus particle} killing the other or enhancing it to give birth.
The competition kernel might depend on the current population size, on the types of the two competing {virus particles} as well as their mutual distance (counting the number of substitutions). In this manner the response of the immun system such as cross-immunity is imlicitly taken into account.
Newborn {virus particles} are clones of their parent unless mutation occurs. When mutation occurs, the parent {particle} gives birth to a mutant whose distances to
 all other {virus particles} (including its own parent) is set to be ``one unit'' longer than the corresponding distances from its parent. We are interested in the huge population, fast branching with frequent mutation but small mutation effects limit.

A similar multi-type branching model describing the evolving empirical trait distributions was introduced in \cite{DG03}.
In that paper the authors even allow for natural branching rates
which might also depend on the population size. However, their model is not suitable for modeling cross-immunity as
their competition kernel does not take into account any information on the traits' history.
Including historical information without loosing the Markov property requires to leave the measure-valued set-up and to work with more enriched state spaces.
This issue is resolved in
\cite{MT12} by
working with historical processes.
% where the mass dependence in their competition kernel models logistic branching.

Historical processes have been established earlier in
\cite{DP91,GLW05}
 for structured neutral populations. They provide an ad-hoc coding of the
  genealogical relationships which is very much intertwined with the spatial
  interaction, where in the context of multi-type models ``spatial''  refers to  mutation in type space rather than migration in geographic space.
  {That means in particular}, that genealogies {cannot} be read off this way in non-spatial
  situations. Moreover, Markov processes with values in measures on path-space
  are notationally far too involved to allow for explicit calculations.
  We therefore propose to rather encode our multi-type phylogenies as
  marked metric measure spaces, and restrict ourselves to
  a dependence of the traits's history only through the phylogenetic distances. This is also closer to applications where the raw data are samples of gene sequences.

  The space of marked metric measure spaces  equipped with the \textbf{marked} Gromov weak topology has been introduced in \cite{Pio10,DGP11}.
  It relies on the idea of encoding ``spatial trees'' as ``tree-like'' metric spaces which are  equipped with a sampling measure, while in addition each point in the tree is assigned a distribution on type-space. It extends the space of (non-spatial) metric measure spaces and the Gromov-weak topology developed in \cite{GPW09} (compare also \cite{Gl12,Loe13}).

  So far only a few examples of dynamics with values in metric measure spaces can be found in the literature. The first paper which considers dynamics
  with values in metric (probability) measure spaces is \cite{GPW13} and constructs the evolving genealogies of (neutral) resampling dynamics. These were extended into type- and state dependent resampling dynamics (including the Otha-Kimura model) in \cite{Pio10} and to resampling dynamics with selection in \cite{DGP12}.
  Note that the dynamics introduced in the latter paper agrees with our dynamics in the case of constant natural branching rates
and conditioned on the total mass being a constant.
Evolving genealogies of spatially structured resampling populations and their continuum space limit are constructed in \cite{GSW16}.
In \cite{Gl12} the evolving genealogies of a (neutral) state-dependent branching dynamics (including Feller's branching diffusion and the Anderson model) are considered. It is worthwhile to note that \cite{Gl12} and the present article do not restrict themselves to metric \textbf{probability} measure spaces but allow for arbitrary \textbf{finite measures}. As a result, the occurrences of explosion or extinction of mass now have to be taken into account as well.

Finally notice that metric measure space valued dynamics for evolving phylogenies from a non-evolutionary point of view were already studied in \cite{EW06}.
\bigskip

The paper is organized as follows: In Section~\ref{S:model} we introduce our model. We will start with recalling the notion of marked metric measure spaces which we will use
for encoding our phylogenies. We then introduce the discrete model of evolving phylogenies for a trait-dependent branching particle model with mutation and competition.
We further state that if we rescale the {virus particle} mass, the mutation steps and the length of distance growth in case of mutation appropriately, then we obtain a tight family of stochastic processes. Finally, we represent the tree-valued dynamics arising in the limit as solutions of a martingale problem. We will defer a proof of uniqueness of this martingale problem to a forthcoming paper. Section~\ref{S:proofdiscrete} is devoted to
uniform moment estimates which will imply that the discrete model is well-defined by its transition rates and which will be used in the tightness argument.
%existence of the  discrete model of evolving phylogenies and moment estimates.
In Section~\ref{S:generators} we identify potential limits {by means of} convergence of generators.
In Section~\ref{S:Tightness} we establish the compact containment condition. Due to heterogeneous exponential branching rates, the phylogenies are in general not ultra-metric. We therefore develop new techniques, including sophisticated coupling techniques.
%In Section~\ref{S:duality} we present a function-valued duality which allows for concluding uniqueness of the limiting martingale problem.
The remaining proofs of the results are collected in Section~\ref{S:proofs}.

% ----------------------------------------------------------------------

\section{The model}
\label{S:model}

In this section we introduce the discrete model of evolving phylogenies of a trait-dependent branching particle model with mutation and competition and its diffusion limit.
In Subsection~\ref{Sub:mmm-space} we start by presenting and extending results on marked metric measure spaces. The discrete model is introduced in Subsection~\ref{Sub:branchmodel}. Small individual mass, frequent mutation with small mutation steps and small substitution distance growth in case of mutation is discussed in Subsection~\ref{Sub:rescale}.
Finally in Subsection~\ref{Sub:mp} we present the diffusion limit.

 \subsection{The state space: Marked metric measure spaces}
 \label{Sub:mmm-space}
In this subsection we define the state space for our evolving phylogenies.
The goal is to capture a phylogeny of a multi-type population in such a way that it allows an explicit description of
the population size as well as of the ancestral relationships and traits in any finite sample.
We will rely on the space of metric measure spaces equipped with the Gromov-weak topology introduced in the context of mono-type populations
in \cite{GPW09}
(see also Chapter~$3\tfrac{1}{2}$ in \cite{Gro99}). Its extension to multi-type populations leads to
marked metric measure spaces equipped with the marked Gromov-weak topology introduced in \cite{DGP11}
(compare also \cite{Pio10}).

As usual, given a topological space $(X,{\mathcal O})$ we denote by
${\mathcal M}_1(X)$ and ${\mathcal M}_f(X)$ the spaces of probability measures
and finite measures on $X$, respectively, defined on
the Borel-$\sigma$-algebra of $X$, and by $\Rightarrow$ weak
convergence in $\mathcal M_1(X)$
and ${\mathcal M}_f(X)$. Recall that the support
$\mathrm{supp}(\mu)$ of $\mu\in\mathcal M_1(X)$ or $\mu\in\mathcal M_f(X)$
is the smallest closed
set $X_0\subseteq X$ such that $\mu(X\setminus X_0)=0$.  The push forward of
$\mu$ under a measurable map $\varphi$ from $X$ into another
topological space $Z$ is the probability measure
$\varphi_\ast\mu\in{\mathcal M}_1(Z)$ defined by
\begin{equation}
\label{push}
   \varphi_\ast\mu(A)
 :=
   \mu\big(\varphi^{-1}(A)\big),
\end{equation} for
all Borel subsets $A\subseteq Z$. We denote by $\mathcal B(X)$ and
${\mathcal C}_b(X)$ the bounded real-valued functions on $X$ which are
measurable and continuous, respectively.\smallskip

The states of evolving phylogenies will be marked metric measure spaces, where the marks are elements in the trait space $K$.
Throughout the paper we will assume that
%\begin{ass}[The trait space]
\begin{center}{\em $K$  is a complete and separable metric space.}\end{center}
%\label{ass:001}
%\end{ass}\sm

The following definition agrees with Definition~2.1 in \cite{DGP11} (albeit stated there for probability rather than finite measures only).
\begin{defi}[mmm-spaces] A $K$-marked metric measure space (mmm-space, for short) is a triple $(X,r,\mu)$ such that $(X,r)$ is a complete and separable metric space and $\mu \in \CM_f(X \times K)$, where $X \times K$ is equipped with the product topology.
\label{def-mmm-spaces}
\end{defi}\sm

We call two mmm-spaces $(X,r_X,\mu_X)$ and $(Y,r_Y,\mu_Y)$ equivalent if they are measure- and mark-preserving isometric, that is there is a measurable $\varphi: \mbox{supp}(\mu_X(\boldsymbol{\cdot}\times K)) \rightarrow \mbox{supp}(\mu_Y(\boldsymbol{\cdot}\times K))$ such that
\begin{equation}
  r_X(x,x') = r_Y(\varphi(x),\varphi(x')) \mbox{ for all } x, x' \in \mbox{supp}(\mu_X(\boldsymbol{\cdot}\times K))
\end{equation}
and
\begin{equation}
  \tilde{\varphi}_* \mu_X = \mu_Y \mbox{ for } \tilde{\varphi}(x,u) = (\varphi(x),u).
\end{equation}
It is clear that the property of being
measure- and mark preserving isometric is an equivalence relation. We write
$\overline{(X,r,\mu)}$ for the equivalence class of a mmm-space
$(X,r,\mu)$. Define the set of (equivalence classes of) $K$-marked metric measure spaces
\begin{equation}
\label{eq:defM}
  \mathbb M^K
  :=
  \big\{\smallx=\overline{(X,r,\mu)} \colon (X,r,\mu)\text{ $K$-marked metric
    measure space}\big\}.
\end{equation}

By Gromov's reconstruction theorem (see \cite[$3\frac12.7$]{Gro99})
a class $\smallx:=\overline{(X,r,\mu)}\in\mathbb{M}^K$ is uniquely characterized by
the total mass
\begin{equation}
\begin{aligned}
\label{e:totalmass}
   m(\smallx):=\mu\big(X\times K\big),
\end{aligned}
\end{equation}
and the {\em marked distance matrix distribution}
\begin{equation}
\begin{aligned}
\label{eq:spaces:distance-matrix-distribution}
   \nu^{\smallx}
 :=
   \big(R^{(X,r),K}\big)_\ast\bar{\mu}^{\otimes \N}
    \in
   \mathcal{M}_1\big(\R_+^{\binom{\N}{2}}\times K^{\mathbb{N}}\big),
\end{aligned}
\end{equation}
with the {\em sampling measure}
\begin{equation}
\begin{aligned}
\label{e:sample}
   \bar{\mu}:=\left\{\begin{array}{cc}\tfrac{\mu}{m(\smallx)}, &  m(\smallx)\not= 0, \\[1mm]
   \mbox{arbitrary in ${\mathcal M}_1(X\times K)$}, & m(\smallx)= 0,\end{array}\right.
\end{aligned}
\end{equation}
and where given a metric space $(X,r)$,
\begin{equation}
\label{eq:spaces:distance-matrix}
     R^{(X,r),K} :=
     \left\{\begin{array}{ll}(X\times K)^{\mathbb N} &\to \mathbb R_+^{\binom{\mathbb N}{2}}\times K^{\mathbb{N}},
\\
     \big((x_i,\kappa_i)_{i\geq 1}\big)
  &\mapsto
     \Big(\big(r(x_i, x_j)\big)_{1\leq i<j},(\kappa_i)_{i\ge 1}\Big)\end{array}\right.
\end{equation}
denotes the \emph{marked distance matrix map}. In what follows we abbreviate
\begin{equation}
\label{e:rm}
\mr:=\big( r_{i,j}\big)_{1\leq i<j}:=\big( r(x_i, x_j)\big)_{1\leq i<j},
\end{equation}
 and
 \begin{equation}
 \label{e:vk}
    \vk:=(\kappa_i)_{i\ge 1}.
\end{equation}

We base our notion of convergence
in $\mathbb{M}^K$ on weak convergence of marked distance matrix distributions.
This topology was
introduced in \cite[Definition~2.4]{DGP11} for marked metric probability measure spaces.

\begin{defi}[Marked Gromov-weak topology]
A sequence $(\smallx_n)_{n\in\N}$ in $\mathbb M^K$ is said to converge marked Gromov-weakly
to $\smallx\in\mathbb M$, as $n\to\infty$, if and only if
\begin{equation}
\label{e:convv2}
  m(\smallx_n)\cdot\nu^{\smallx_n} \Tno m(\smallx)\cdot\nu^\smallx
\end{equation}
in the weak topology on $\mathcal M_f(\mathbb R_+^{\binom{\mathbb
    N}{2}}\times K^\mathbb{N})$. Denote this by $\smallx_n \rightarrow \smallx$.
\label{Def:005}
\end{defi}\sm

Note that we cannot immediately conclude from the fact that the space of finite measures equipped with the weak topology is a Polish space
that $\mathbb{M}^K$ equipped with the marked Gromov-weak topology is Polish. The reason is that  the weak limit of a sequence of marked distance
matrix distributions is in general not the marked distance matrix distribution of a metric measure space.
The following result extents \cite[Theorem~2]{DGP11} from marked metric probability to marked metric
(finite) measure spaces.
\begin{proposition}[Polish space]
The space $\M^K$ equipped with the marked Gromov-weak topology is Polish.
\label{P:007}
\end{proposition}\sm

This follows from the extension of the characterization of compact sets given in \cite[Theorem~3]{DGP11}
(compare also \cite[Proposition~7.1, Remark~7.2]{GPW09}):
\begin{proposition}[Relative compactness in $\mathbb{M}^K$]
A family $\Gamma\subset\M^K$ is relatively compact if and only if for all $\varepsilon>0$ there exists $N_\varepsilon\in\N$ and a compact subset $K_\varepsilon\subset K$  such that for all $\smallx =\overline{(X,r,\mu)}\in\Gamma$, the following hold:
\begin{enumerate}
\item[(i)] $m(\smallx) \leq N_\varepsilon$.
\item[(ii)] $\mu(X\times K_\varepsilon^c) \leq \varepsilon$.
\item[(iii)] There is a subset $X_\varepsilon \subseteq X$ with
\begin{enumerate}
\item[(iii-a)]
$\mu(X_\varepsilon^c\times K) \leq \varepsilon$,
\item[(iii-b)]
$X_\varepsilon$ has diameter at most $N_\varepsilon$, and
\item[(iii-c)]
$X_\varepsilon$ can be covered by at most $N_\varepsilon$ balls of radius $\varepsilon$.
\end{enumerate}
\end{enumerate}
\label{pro-rel-comp}
\end{proposition}\sm

\noindent\textit{Proof of Proposition~2.3 and 2.4.} This follows from Theorem~3 in \cite{DGP11} together with the fact that $\smallx_n\to\smallx$ if and only if $m(\smallx_n)\to m(\smallx)$ and $\nu^{\smallx_n}\Tno \nu^{\smallx}$ in case $m(\smallx) \neq 0$, or equivalently, if  $m(\smallx_n)\to m(\smallx)$ and $m(\smallx)\nu^{\smallx_n}\Tno m(\smallx)\nu^{\smallx}$.
\qed\sm

%\begin{proof}[Proof of Proposition~\ref{P:001}]
%{\color{blue}\tt\tiny xxx fill in if requested}
%\end{proof}\sm

% ----------------------------------------------------------------------

 \subsection{The trait-dependent branching particle model with mutation and competition}
 \label{Sub:branchmodel}
 In this subsection we introduce our individual based model in detail.

\begin{figure}
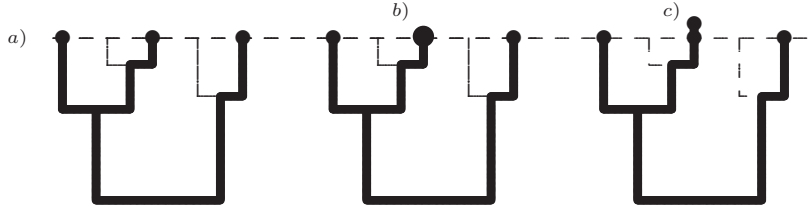

 \beginpicture
\setcoordinatesystem units <.6cm,.6cm> \setplotarea x from -6 to 8, y from 0 to 4.7
  \put{\Large $\bullet$} [cC] at 0 4
  \put{\Large $\bullet$} [cC] at 2 4
  \put{\Large $\bullet$} [cC] at 4 4
  \put{\Large $\bullet$} [cC] at 6 4
  \put{\huge $\bullet$} [cC] at 8 4
    \put{\Large $\bullet$} [cC] at 10 4
 % \put{\huge $\bullet$} [cC] at 8 4
  %\put{$\circ$} [cC] at 8.05 4.1
  \multiput {\tiny $\bullet$} at 0 4 *160 0 -.01 /
  \multiput {\tiny $\bullet$} at 0 2.4 *150 .01 0 /
  \multiput {\tiny $\bullet$} at 1.5 3.4 *100 0 -.01 /
  \multiput {\tiny $\bullet$} at 1.5 3.4 *50 .01 0 /
  \multiput {\tiny $\bullet$} at 2 3.4 *60 0 .01 /
  \multiput{\tiny $\bullet$} at 4 4 *130 0 -.01 /
  \multiput {\tiny $\bullet$} at 4 2.7 *50 -.01 0 /
  \multiput {\tiny $\bullet$} at 3.5 2.7 *230 0 -.01 /
  \multiput {\tiny $\bullet$} at .75 0.4 *275 .01 0 /
  \multiput {\tiny $\bullet$} at .75 0.4 *200 0 .01 /
  \multiput {\tiny $\bullet$} at 6 4 *160 0 -.01 /
  \multiput {\tiny $\bullet$} at 6 2.4 *150 .01 0 /
  \multiput{\tiny $\bullet$} at 7.5 3.4 *100 0 -.01 /
  \multiput {\tiny $\bullet$} at 7.5 3.4 *50 .01 0 /
  \multiput{\tiny $\bullet$} at 8 3.4 *60 0 .01 /
\multiput {\tiny $\bullet$} at 10 4 *130 0 -.01 /
\multiput {\tiny $\bullet$} at 10 2.7 *50 -.01 0 /
\multiput {\tiny $\bullet$} at 9.5 2.7 *230 0 -.01 /
\multiput {\tiny $\bullet$} at 6.75 0.4 *275 0.01 0 /
\multiput {\tiny $\bullet$} at 6.75 0.4 *200 0 .01 /
  \plot 1 4 1 3.4 2 3.4 2 4 /
  \plot 0 4 0 2.4 1.5 2.4 1.5 3.4 /
  \plot 3 4 3 2.7 4 2.7 4 4 /
  \plot 0.75 2.4 0.75 .4 3.5 .4 3.5 2.7 /
  \plot 7 4 7 3.4 8 3.4 8 4 /
  \plot 6 4 6 2.4 7.5 2.4 7.5 3.4 /
  \plot 9 4 9 2.7 10 2.7 10 4 /
  \plot 6.75 2.4 6.75 .4 9.5 .4 9.5 2.7 /
  \setdashes \plot -.2 4 10.5 4 /
%\plot -.2 3.2 7.5 3.2 /
  \put{\tiny } [cC] at 5.1 1.8
%  \multiput {$\bullet$} at 12 3 *100 .01 0 /
%\color{black}
%  \plot 12 2 13 2 /
%  \put{sample} [lC] at 13.3 3
%  \put{whole tree} [lC] at 13.3 2
%%  \put{\Large $\bullet$} [cC] at 0 4
%%  \put{\Large $\bullet$} [cC] at 2 4
%%  \put{\Large $\bullet$} [cC] at 4 4
  \put{\Large $\bullet$} [cC] at 12 4
  \put{\Large $\bullet$} [cC] at 14 4
    \put{\Large $\bullet$} [cC] at 16 4
 \put{\Large $\bullet$} [cC] at 14 4.3
  %\put{$\circ$} [cC] at 8.05 4.1
%%  \multiput {\tiny $\bullet$} at 0 4 *160 0 -.01 /
%%  \multiput {\tiny $\bullet$} at 0 2.4 *150 .01 0 /
%%  \multiput {\tiny $\bullet$} at 1.5 3.4 *100 0 -.01 /
%%  \multiput {\tiny $\bullet$} at 1.5 3.4 *50 .01 0 /
%%  \multiput {\tiny $\bullet$} at 2 3.4 *60 0 .01 /
%%  \multiput{\tiny $\bullet$} at 4 4 *130 0 -.01 /
%%  \multiput {\tiny $\bullet$} at 4 2.7 *50 -.01 0 /
%%  \multiput {\tiny $\bullet$} at 3.5 2.7 *230 0 -.01 /
%%  \multiput {\tiny $\bullet$} at .75 0.4 *275 .01 0 /
%%  \multiput {\tiny $\bullet$} at .75 0.4 *200 0 .01 /
\put{\mbox{\tiny $a)$}} at -1 4
\put{\mbox{\tiny $b)$}} at 7.5 4.6
\put{\mbox{\tiny $c)$}} at 13.5 4.6
  \multiput {\tiny $\bullet$} at 12 4 *160 0 -.01 /
  \multiput {\tiny $\bullet$} at 12 2.4 *150 .01 0 /
  \multiput{\tiny $\bullet$} at 13.5 3.4 *100 0 -.01 /
  \multiput {\tiny $\bullet$} at 13.5 3.4 *50 .01 0 /
  \multiput{\tiny $\bullet$} at 14 3.4 *60 0 .01 /
  \multiput{\tiny $\bullet$} at 14 4 *20 0 .01 /
\multiput {\tiny $\bullet$} at 16 4 *130 0 -.01 /
\multiput {\tiny $\bullet$} at 16 2.7 *50 -.01 0 /
\multiput {\tiny $\bullet$} at 15.5 2.7 *230 0 -.01 /
\multiput {\tiny $\bullet$} at 12.75 0.4 *275 0.01 0 /
\multiput {\tiny $\bullet$} at 12.75 0.4 *200 0 .01 /
%%  \plot 1 4 1 3.4 2 3.4 2 4 /
%%  \plot 0 4 0 2.4 1.5 2.4 1.5 3.4 /
%%  \plot 3 4 3 2.7 4 2.7 4 4 /
%%  \plot 0.75 2.4 0.75 .4 3.5 .4 3.5 2.7 /
  \plot 13 4 13 3.4 14 3.4 14 4 /
  \plot 14 4 14 4.4 /
  \plot 12 4 12 2.4 13.5 2.4 13.5 3.4 /
  \plot 15 4 15 2.7 16 2.7 16 4 /
  \plot 12.75 2.4 12.75 .4 15.5 .4 15.5 2.7 /
  \setdashes \plot -.2 4 16.5 4 /
\endpicture
{\caption{\small shows the phylogeny of a sample of a phylogeny of sample size $3$,  a) before the branching event, and after the second individual gave birth to b) a clone (indicated by the weight enlarged by $\zeta$ on the clan containing this clone) or c) a mutant (indicated by the mutant child being in positive genetic distance $\varsig$ to its mother).}}
\label{Fig:002}
\end{figure}

 Fix a {\em distance constant} $\varsig>0$ and an {\em individual mass constant} $\zeta>0$.
 The trait-dependent branching particle model with mutation and competition takes values in phylogenies of a finite population in which each individual has mass $\zeta>0$ and
 substitution distances are assumed to be multiples of $\varsig$ (compare Figure~\ref{Fig:002}).

 The state space is therefore
 the subspace of $(\varsig,\zeta)$-marked metric measure spaces defined as
 \begin{equation}
 \label{s:001}
 \begin{aligned}
    \mathbb{M}^{K,(\varsig,\zeta)}
  :=
    \big\{\smallx&=\overline{(X,r,\mu)}\in\mathbb{M}^K:\,
    \\
    &\zeta^{-1}\mu(\boldsymbol{\cdot}\times K)\in{\mathcal N}_f(X);\,\#\mathrm{supp}\big(\mu(\{x\}\times\boldsymbol{\cdot})\big)\in\{0,1\},\,\forall x\in X,
 %   \\
 %   &\mu(\mathrm{d}(x,\kappa))=\mu(\mathrm{d}x\times K)\otimes\delta_{\tilde{\kappa}(x)}(\mathrm{d}\kappa),\mbox{ for some $\tilde{\kappa}:X\to K$}
    \\
   &\varsig^{-1} r(x,x')\in\mathbb{N}\cup\{0\},\,\forall x,x'\in X\big\},
  \end{aligned}
 \end{equation}
 where ${\mathcal N}_f(X)$ denotes the space of non-negative integer-valued measures on $X$.

\begin{remark}[Mark functions for finite populations]
Notice that if $\smallx \in \mathbb{M}^{K,(\varsig,\zeta)}$ there exists a {measurable} {\it mark function} $\tilde{\kappa}: X \rightarrow K$ such that
\begin{equation}
  \mu(\mathrm{d}(x,\kappa)) = \mu(\mathrm{d}x \times K) \otimes \delta_{\tilde{\kappa}(x)}(\mathrm{d}\kappa).
\end{equation}
In Section~\ref{Sub:rescale} it is stated that a suitably rescaled family of tree-valued $(\varsig_N,\zeta_N,\alpha_N)$-trait-dependent branching dynamics ${\mathcal X}^N \in \mathbb{M}^{K,(\varsig_N,\zeta_N)}, N \in \mathbb{N}$ with mutation and competition is tight in $\mathbb{M}^K$. It remains as an open question, whether all subsequential limits admit a mark function. As a first step in this direction,  a criterion for the existence of a mark function has been derived in \cite{KL15} and
successfully been applied to the tree-valued Fleming-Viot dynamics with mutation and selection from \cite{DGP12}.
This criterion has been extended to marked metric boundedly finite measure spaces and then applied to the spatially interacting Fleming-Viot dynamics on $\mathbb{Z}$ in \cite{GSW16}.
Furthermore, in \cite{GSW16} a version of the continuum limit of these spatially interacting tree-valued Fleming-Viot dynamics with a mark function is constructed by means of the Brownian web.
\label{Rem:002}
\end{remark}\sm

 \begin{remark}[Link between $(\varsig,\zeta)$-marked metric measure spaces and evolving phylogenies]
 If $(X,r,\mu)$ is such that $r$ is only a pseudo-metric on $X$, (that is
$r(x,y)=0$ is possible for $x\neq y$) we can still define its measure-preserving
isometry class under the additional assumption that if $r(x,y)=0$, then $\mu(\{x\}\times\boldsymbol{\cdot})=\mu(\{y\}\times\boldsymbol{\cdot})$. Since this class contains also metric measure spaces, there
is a bijection between the set of equivalence classes of  pseudo-metric
measure spaces and the set of equivalence classes of  metric measure
spaces and we will use both notions interchangeably.

It turns out that this is particularly suited to our set-up. Given a pseudo-metric
measure space $(X,r,\mu)$,
we would like to think of two
individuals $x\not =y\in X$ with $\mu(\{x\} \times K)\cdot\mu(\{y\} \times K)>0$ but $r(x,y)=0$ as {\em clones (of each other)}.
In the corresponding equivalent metric measure
space they are all collected into one {\em clan (of clones)}. If $\smallx=\overline{(X,r,\mu)}\in\mathbb{M}^{K,(\varsig,\zeta)}$, we can recover the number of clones in clan $x\in X$ as
\begin{equation}
\label{e:003a}
   n_x:=\tfrac{1}{\zeta}\cdot\mu\big(\{x\}\times K\big),
\end{equation}
as $\zeta$ stands for an individual's mass. In particular,
\begin{equation}
\label{e:003b}
   n:=\sum\nolimits_{x\in X}n_x=\tfrac{m(\smallx)}{\zeta}
\end{equation}
relates the {\em total number of individuals} in the population with the total mass.
\hfill$\qed$
 \label{Rem:001}
 \end{remark}
 \sm

 To be in a position to introduce the discrete model, fix $p\in[0,1]$ which we refer to as the {\em mutation probability}. Moreover,
 let $\alpha$ be a stochastic kernel on $K$.
 In the following we refer to $\alpha$ as the {\em mutation kernel}.
 We will exclude the possibility of {\em non-trivial mutation}, that is, we assume
for every  $\kappa\in K$,
 \begin{equation}
 \label{e:003c}
    \alpha\big(\kappa,\{\kappa\}\big)=0.
 \end{equation}
Further let $\beta:K\to[0,\infty), \gamma^{\mathrm{birth}}:\R_+\times\R_+\times K^2\to[0,\infty)$ and $\gamma^{\mathrm{death}}:\R_+\times\R_+\times K^2\to[0,\infty)$.

The {\em trait-dependent branching particle model with mutation and competition} is a Markov process which takes values in $\mathbb{M}^{K,(\varsig,\zeta)}$. Given that it starts in
 $\smallx_0\in\mathbb{M}^{K,(\varsig,\zeta)}$, it has the following dynamics:
 \begin{itemize}
 \item {\bf Death. } For any particle from clan $x_1$ of trait $\kappa_1=\tilde{\kappa}(x_1)$ and any particle from clan $x_2$ of trait $\kappa_2=\tilde{\kappa}(x_2)$,
 at rate
\begin{equation}
\label{s:002aa}
 \tfrac{\beta(\kappa_2)}{m(\smallx)}+\zeta\cdot\tfrac{\gamma^{\mathrm{death}}(m(\smallx),r(x_1,x_2),
 \kappa_1,\kappa_2)}{m(\smallx)},
 %\end{equation}
\end{equation}
the second particle dies either due to {\em natural death} or because it gets {\em killed} by the first particle.
That is, the total death rate for a particle of clan $x_2$ and trait $\kappa_2=\tilde{\kappa}(x_2)$ is
{
 \begin{equation}
 \label{s:002a}
 \begin{aligned}
    %\zeta^{-1}\int_{X^2}\bar{\mu}\otimes\mu(\mathrm{d}((x_1,\kappa_1),(x_2,\kappa_2)))\big\{\tfrac{\beta(\kappa_2)}{\zeta}+\gamma^{\mathrm{death}}(m(\smallx),r(x_1,x_2),\kappa_1,\kappa_2)\big\},
     &\tfrac{\beta(\kappa_2)}{\zeta}+\sum\nolimits_{x_1 \in X}\tfrac{n_{x_1}}{n}\gamma^{\mathrm{death}}(m(\smallx),r(x_1,x_2),\tilde{\kappa}(x_1),\kappa_2)
     \\
     &=
     \tfrac{\beta(\kappa_2)}{\zeta}+\sum\nolimits_{x_1 \in X}\bar{\mu}(\{x_1\}\times K)\gamma^{\mathrm{death}}(m(\smallx),r(x_1,x_2),\tilde{\kappa}(x_1),\kappa_2).
     \end{aligned}
 \end{equation}

Such a death event yields the following transition:
 \begin{equation}
 \label{s:002b}
    \smallx\mapsto\big(X,r,\mu-\zeta\cdot\delta_{(x_2,\kappa_2)}\big).
 \end{equation}
Notice that if $\zeta^{-1}\cdot\mu(\{x_2\}\times K)=1$, then $(X,r,\mu-\zeta\cdot\delta_{(x_2,\kappa_2)})$ is equivalent to $(X\setminus\{x_2\},r,\mu-\zeta\cdot\delta_{(x_2,\kappa_2)})$.}

 \item  {\bf Birth. } For any particle from clan $x_1$ of trait $\kappa_1=\tilde{\kappa}(x_1)$ and any particle from clan $x_2$ of trait $\kappa_2$,
 at rate
 \begin{equation}
 \label{s:002dd}
    \tfrac{\beta(\kappa_2)}{m(\smallx)}+\zeta\cdot\tfrac{\gamma^{\mathrm{birth}}(m(\smallx),r(x_1,x_2),\kappa_1,\kappa_2)}{m(\smallx)}
 \end{equation}
 the second particle gives birth either naturally or due to {\em birth-enhancement} by the first particle.
 That is, the total birth rate of a particle from clan $x_2$ of trait $\kappa_2$ is
\begin{equation}
 \label{s:002d}
 \begin{aligned}
    %\zeta^{-1}\int_{X^2}\bar{\mu}\otimes\mu(\mathrm{d}((x_1,\kappa_1),(x_2,\kappa_2)))\big\{\tfrac{\beta(\kappa_2)}{\zeta}+\gamma^{\mathrm{death}}(m(\smallx),r(x_1,x_2),\kappa_1,\kappa_2)\big\},
     &\tfrac{\beta(\kappa_2)}{\zeta}+\sum\nolimits_{x_1 \in X}\tfrac{n_{x_1}}{n}\gamma^{\mathrm{birth}}(m(\smallx),r(x_1,x_2),\tilde{\kappa}(x_1),\kappa_2)
      \\
     &=
     \tfrac{\beta(\kappa_2)}{\zeta}+\sum\nolimits_{x_1 \in X}\bar{\mu}(\{x_1\}\times K)\gamma^{\mathrm{birth}}(m(\smallx),r(x_1,x_2),\tilde{\kappa}(x_1),\kappa_2).
     \end{aligned}
 \end{equation}

 With probability $p\in[0,1]$ the newborn $z\not\in X$ is a mutant of its parent $x_2$ whose new type $\tilde{\kappa}(z)$ is chosen with respect to $\alpha(\kappa_2,\cdot)$,
   and its distance to all other particles is given by
  \begin{equation}
  \label{s:006}
     r^{(x_2,z),\varsig}{(z,x)}:=\left\{\begin{array}{rc}r(x_2,x)+\varsig, & x\in X
     \\ 0, & x=z,\end{array}\right.
  \end{equation}
  while with probability
  $1-p$ the newborn is just a clone of its parent.
Under our assumption (\ref{e:003c}) on non-trivial mutation,
  this yields the following transition:
 \begin{equation}
 \label{s:004}
    \smallx\mapsto\big(X\uplus\{z\},r^{(x_2,z),\varsig},\mu+\zeta\cdot\delta_{(z,\tilde{\kappa}(z))} \cdot \mathbf{1}\{ \tilde{\kappa}(z) \not= \kappa_2 \} +\zeta\cdot\delta_{(x_2,\kappa_2)} \cdot \mathbf{1}\{ \tilde{\kappa}(z) = \kappa_2 \} \big),
 \end{equation}
 where $\tilde{\kappa}(z)$ is chosen with respect to
 \begin{equation}
\label{s:005}
    \widehat{\alpha}_N\big(\kappa_2,\cdot\big):=p\cdot\alpha\big(\kappa_2,\cdot\big)+(1-p)\cdot\delta_{\kappa_2}(\cdot).
 \end{equation}
% \hfill$\qed$
 \end{itemize}
We remark that if there is no mutation, $\big(X\uplus\{z\},r^{(x_2,z),\varsig},\mu+\zeta\cdot\delta_{(x_2,\kappa_2)} \big)$ is equivalent to $\big(X,r^{(x_2,z)},\mu+\zeta\cdot\delta_{(x_2,\kappa_2)} \big)$. In this case, the clan $\{z\}$ cannot effect new birth events, as there are no particles in clan $z$ due to $n_z=0$ respectively $\mu(\{z\}\times K)=0$.

\iffalse{
Notice that given $\smallx\in\mathbb{M}^{K,(\varsig,\zeta)}$, the total jump rate is
\begin{equation}
\begin{aligned}
\label{s:100}
   \lambda(\smallx)
 &:=
   \tfrac{m(\smallx)}{\zeta}\big(2\cdot\tfrac{\hat{\beta}^{\smallx}}{\zeta}+
   \widehat{\gamma}^{\mathrm{death}}\big(m(\smallx)\big)+\widehat{\gamma}^{\mathrm{birth}}
   \big(m(\smallx)\big)\big)<\infty,
\end{aligned}
\end{equation}
with
 \begin{equation}
\label{hat:beta}
   \widehat{\beta}^{\smallx}:=\int\nu^{\smallx}(\mathrm{d}(\mr,\vk))\,\beta\big(\kappa_1\big),
\end{equation}
and
\begin{equation}
\label{hat:gammad}
   \widehat{\gamma}^{\mathrm{death/birth}}\big(m(\smallx)\big)
 :=
   \int\,\nu^{\smallx}(\mathrm{d}(\mr,\vk))\,
   \gamma^{\mathrm{death/birth}}\big(m(\smallx),r_{1,2},\kappa_1,\kappa_2\big).
\end{equation}

As the total rates
are finite for each state $\smallx\in\mathbb{M}^{K,(\varsig,\zeta)}$,
the above description gives a well-defined strong Markovian pure jump process with possibly finite life time due to accumulation of jumps. }\fi

Notice that the jump process described above is well-defined if we can exclude the explosion of the total rate in finite time. To ensure this we make the following assumptions on the branching and the competition rates.
 \begin{ass}[Bounded branching rate] The function $\beta:K\to[0,\infty)$ is {measurable and} bounded, that is, there is a constant $\overline{\beta}<\infty$ such that
 \begin{equation}
 \label{e:005}
    \sup_{\kappa\in K}\beta(\kappa)\le\overline{\beta}.
 \end{equation}
\label{ass:002}
 \end{ass}\sm

 \begin{ass}[Bounded birth-enhancement rate]
The function $\gamma^{\mathrm{birth}}:\R_+\times\R_+\times K^2\to[0,\infty)$ is  {measurable and} bounded, that is,
  there is a constant $\overline{\gamma}_b<\infty$ such that
  \begin{equation}
\label{e:bds-comp-rates-b}
\begin{aligned}
  &\sup_{m,r\in\R_+,\kappa,\kappa'\in K}\gamma^{\mathrm{birth}}(m,r,\kappa,\kappa') \leq \overline{\gamma}_b.
\end{aligned}
\end{equation}
\label{ass:007b}
\end{ass}\sm

\begin{ass}[Bounded death-competition rate]
There is a continuous function $\tilde{\gamma}:\R_+\to[0,\infty)$ such that for all $m\in\mathbb{R}_+$, the function $\gamma^{\mathrm{death}}:\R_+\times\R_+\times K^2\to[0,\infty)$  {is measurable and} satisfies
 \begin{equation}
 \label{s:008}
    \sup_{r\in\R_+,\kappa,\kappa'\in K}\gamma^{\mathrm{death}}\big(m,r,\kappa,\kappa'\big)\le\tilde{\gamma}(m).
 \end{equation}
\label{ass:003}
\end{ass}\sm

Under these assumptions the total jump rate can be upper bounded by a continuous function of the total mass, which itself can be stochastically upper bounded by a Yule process with rate $\zeta^{-2}\bar{\beta}+\zeta^{-1}\bar{\gamma}_b$. We can therefore immediately conclude the following:
\begin{proposition}[Well-posed] Under Assumptions~\ref{ass:002},~\ref{ass:007b} and~\ref{ass:003}, there is a well-defined strong Markov pure jump process, ${\mathcal X}=({\mathcal X}_t)_{t\ge 0}$,
with the above described transition rates provided we start in a random ${\mathcal X}_0\in\mathbb{M}^{K,(\varsig,\zeta)}$ with $\mathbb{E}[m({\mathcal X}_0)]<\infty$.
\label{P:004}
\end{proposition}\sm

This leads to the following definition.
\begin{definition}[$(\varsig,\zeta,\alpha)$-trait-dependent branching with mutation and competition]
The process ${\mathcal X}=({\mathcal X}_t)_{t\ge 0}$ (starting in ${\mathcal X}_0$ with $\mathbb{E}[m({\mathcal X}_0)]<\infty$) is referred to as {\em tree-valued $(\varsig,\zeta,\alpha)$-trait-dependent branching dynamics with mutation and competition}.
%In the following we write $\mathbb{P}^{(\varsig,\zeta,\alpha)}:=\{\mathbb{P}^{(\varsig,\zeta,\alpha)}_{\smallx},\,\smallx\in\mathbb{M}^{K,(\varsig,\zeta)}\}$ for its law.
\label{Def:002}
\end{definition}\sm

% ----------------------------------------------------------------------

\subsection{The fast evolving small mass and large population rescaling}
\label{Sub:rescale}
In this subsection we state that the suitably rescaled family of tree-valued $(\varsig,\zeta,\alpha)$-trait-dependent branching dynamics with mutation and competition is tight.

%Recall $(\Omega^{(\varsig,\zeta)},{\mathcal D}(\Omega^{(\varsig,\zeta)}))$ from (\ref{e:OmegaN}) and (\ref{s:007}).
We are interested in a rescaling where the particles branch fast and a mutation event is a typical event while the mutation steps are small. This translates into letting the parameters
$\zeta$ and $\varsig$ tend to zero such that the fraction of $\varsig$ and $\zeta$ converges to a non-trivial limit. Note that the mutation parameter $p$ is kept constant. For simplicity, we choose for
every $N\in\mathbb{N}$,
\begin{itemize}
\item {\bf Small distance constant. } $\varsig_N:=\tfrac{1}{N}$.
\item {\bf Small individual mass. } $\zeta_N:=\tfrac{1}{N}$.
\end{itemize}

We also make the following assumption.

\begin{ass}[Mutation operator] Consider a family $\{\alpha_N;\,N\in\mathbb{N}\}$ of stochastic kernels on $K$ such that
the following holds:
\begin{itemize}
\item[(i)] {\bf Mutation processes along one line are tight. }Let for  fixed $N\in\mathbb{N}$, $\smallk^N$ be the $K$-valued Markov process which jumps at rate
$N\cdot\alpha_N(\kappa,\mathrm{d}\tilde{\kappa})$ from $\kappa$ to $\tilde{\kappa}$.
 The family of jump processes $\{\smallk^N;\,N\in\mathbb{N}\}$
is assumed to be tight.
\item[(ii)] {\bf Limit mutation along one line is uniquely characterized. } There is a linear operator  $(A,{\mathcal D}(A))$
on ${\mathcal C}_b(K)$ such that  ${\mathcal D}(A)$ is an algebra, dense in ${\mathcal C}_b(K)$ and the $(A,{\mathcal D}(A))$-martingale problem has a unique solution, and
for all $h\in{\mathcal D}(A)$ and $\kappa\in K$,
\begin{equation}
\label{s:012}
   N\cdot\int_K\big(\alpha_N(\kappa,\mathrm{d}\tilde{\kappa})-\delta_\kappa(\mathrm{d}\tilde{\kappa})\big)h(\tilde{\kappa})\tNo Ah(\kappa),
\end{equation}
uniformly in $\kappa \in K$. Here, as usual, $\delta_\kappa(\cdot)=\delta(\kappa,\cdot)$ denotes the dirac measure.
Note that (\ref{s:012}) implies that $\int_K \alpha_N(\boldsymbol{\cdot},\mathrm{d}\tilde{\kappa})f(\tilde{\kappa}) \tNo f$ whenever $f$ is in ${\mathcal D}(A)$.
\end{itemize}
\label{ass:004}
\end{ass}\sm

We make the following further assumptions on the natural birth and the death-competition rate.
\begin{ass}[Natural birth rate bounded away from zero] There exists a constant $\underline{\beta}>0$ such that
\begin{equation}
\label{e:bds-comp-rates-d1}
\begin{aligned}
  &\inf_{\kappa\in K}\beta(\kappa) \geq \underline{\beta}.
\end{aligned}
\end{equation}
\label{ass:007c}
\end{ass}\sm

\begin{remark} We point out that Assumptions~\ref{ass:007b} and~\ref{ass:007c}
imply that there exist a $C\in(0,\infty)$ such that for all $\kappa_2\in K$,
\begin{equation}\label{e:011}
  \sup_{m\in\R_+,r\in\R_+,\kappa_1\in K}\gamma^{\mathrm{birth}}(m,r,\kappa_1,\kappa_2)\le C\cdot\beta(\kappa_2).
\end{equation}
\label{Rem:003}
\end{remark}

\begin{ass}[Linear bound on the competition death rate] There exists a constant $\overline{\gamma}_d$ such that
\begin{equation}
\label{e:bds-comp-rates-d2}
\begin{aligned}
  &\sup_{m,r\in\R_+,\kappa,\kappa'\in K}\gamma^{\mathrm{death}}(m,r,\kappa,\kappa') \leq (1\vee m) \overline{\gamma}_d.
\end{aligned}
\end{equation}
\label{ass:007a}
\end{ass}\sm

Our first main result is the following:
\begin{theorem}[Tightness] Let for each $N\in\mathbb{N}$, ${\mathcal X}^N$ be the  tree-valued $(\varsig_N,\zeta_N,\alpha_N)$-trait-dependent branching dynamics with mutation and competition
such that
$\{\mathcal X_0^N;\,N\in\mathbb{N}\}$ is a tight family in $\mathbb{M}^K$ with $\sup_{N\in\mathbb{N}}\mathbb{E}[(m({\mathcal X}^N_0))^3]<\infty$.
Under Assumptions~\ref{ass:002},~\ref{ass:007b},~\ref{ass:004},~\ref{ass:007c}  and~\ref{ass:007a},
the family $\{{\mathcal X}_N;\,N\in\mathbb{N}\}$ is tight.
\label{T:tightness}
\end{theorem}\sm

\begin{definition}[Tree-valued trait-dependent branching with mutation and competition] Any limit process is called tree-valued trait-dependent branching dynamics with mutation and competition.
%Denote its law by $\mathbb{P}$.
\label{Def:004}
\end{definition}\sm

\subsection{The martingale problem of the limit dynamics}
\label{Sub:mp}
In this subsection we present an analytic representation of the limit process in terms of a  martingale problem.

We begin by introducing a class of suitable test functions.

\begin{definition}[Polynomials] A polynomial is a function $F=F^{\ups}:\mathbb{M}^K\to\R$ of the form
\begin{equation}
\label{s:009}
   F(\smallx):=\int_{\R_+^{\N\choose 2}\times K^\mathbb{N}}\nu^{\smallx}(\mathrm{d}(\mr,\vk))\,\ups\big(m(\smallx),\mr,\vk\big)
\end{equation}
for a function $\ups\in{\mathcal C}_b(\R_+\times\R_+^{\N\choose 2}\times K^\mathbb{N})$ such that there is a constant $\bar{\ups}_0\in\mathbb{R}$ such that
\begin{equation}
\label{s:011}
   \ups(0,\boldsymbol{\cdot},\boldsymbol{\cdot})\equiv \bar{\ups}_0.
\end{equation}

If the function $\ups\in{\mathcal C}_b(\R_+\times{\R_+^{\N\choose 2}\times K^\mathbb{N}})$ depends on $(m,(r_{i,j})_{1\le i<j},(\kappa_i)_{i\in\mathbb{N}})$ only through
$(m,(r_{i,j})_{1\le i<j\le n},(\kappa_i)_{i\le n})$ for some $n \in \mathbb{N}$, then we refer to $F$ as a polynomial of finite degree.

Denote by $\Pi$ the space of all polynomials, and by $\Pi_{\mathrm{finite}}$ the subspace of polynomials of finite degree.
\label{Def:003}
\end{definition}\sm

Recall the mutation operator $(A,{\mathcal D}(A))$ from Assumption~\ref{ass:004}.
We consider for each  $l_1,l_2\in\mathbb{N}$ the subspace
\begin{equation}
\label{e:004}
\begin{aligned}
   \Pi_{0}^{l_1,l_2,A}
  &:=
    \Big\{F=F^\ups\in\Pi:\,\ups(\boldsymbol{\cdot},\mr,\vk)\in{\mathcal C}_b^{l_1}(\R_+),\,\forall(\mr,\vk);\,\ups(m,\boldsymbol{\cdot},\vk)\in{\mathcal C}_b^{l_2}(\R_+^{\mathbb{N}\choose 2}),\,\forall(m,\vk);
    \\
  &\hspace{3.5cm}\ups(m,\mr,(\kappa_1,...,\kappa_{l-1},\boldsymbol{\cdot},\kappa_{l+1},...))\in{\mathcal D}(A),\,\forall l\in\mathbb{N},(m,\mr,\vk)\Big\}.
\end{aligned}
\end{equation}

Recall the branching rate $\beta:K\to\R_+$ and the competition rates $\gamma^{\mathrm{birth}},\gamma^{\mathrm{death}}:\R_+\times\R_+\times K^2\to\R_+$, and put
\begin{equation}
\label{s:010}
   \Gamma\big(m(\smallx),r_{1,2},\kappa_1,\kappa_2\big):=\gamma^{\mathrm{birth}}\big(m(\smallx),r_{1,2},\kappa_1,\kappa_2\big)-
   \gamma^{\mathrm{death}}\big(m(\smallx),r_{1,2},\kappa_1,\kappa_2\big).
\end{equation}

In Proposition~\ref{e:021} the following operator is obtained as the limit of a sequence of operators corresponding to the suitably rescaled approximating individual based models (cf. Subsections~\ref{Sub:branchmodel} and~\ref{Sub:rescale}, and~(\ref{e:OmeggaN})).
Consider the operator $\Omega$ acting on
\begin{equation}
\label{e:DOmega}
   {\mathcal D}(\Omega):=\big\{F\in \Pi^{2,1,A}:\,\Omega F\mbox{ is well-defined and finite}\big\},
\end{equation}
where
\begin{equation}
\label{y:001}
\begin{aligned}
   \Omega
  &:=
   \Omega^{\beta,\Gamma}_{\mbox{\tiny total mass}}
%  +
%    {\color{green} \Omega^{\beta}_{\mbox{\tiny mass flow}}}
  +
   \Omega^{p,\beta,A}_{\mbox{\tiny trait mutation}}
  +
   \Omega^{p,\beta}_{\mbox{\tiny growth}}
  +
   \Omega^{\Gamma}_{\mbox{\tiny $\Gamma$-reweigh}}
 % +
 %   {\color{green} \Omega^{\beta}_{\mbox{\tiny $\beta$-reweigh}}}
 % +
 %  \Omega^{\beta}_{\mbox{\tiny resample}}
  +
   \Omega^{\beta}_{\mbox{\tiny natural branching}}
\end{aligned}
\end{equation}
is reflecting
\begin{enumerate}
\item $\Omega^{\beta,\Gamma}_{\mbox{\tiny total mass}}$ the changes in the total mass due to competition, fluctuation and a flow due to reweighing the sampling measure with respect to updating $\beta$,
\item $\Omega^{p,\beta,A}_{\mbox{\tiny trait mutation}}$ trait mutation,
\item $\Omega^{p,\beta}_{\mbox{\tiny growth}}$ growth of substitution distances,
\item $\Omega^{\Gamma}_{\mbox{\tiny $\Gamma$-reweigh}}$ a flow due to reweighing the sampling measure with respect to updating $\Gamma$, and
\item $\Omega^{\beta}_{\mbox{\tiny natural branching}}$ the effects on genealogies spanned by a sample of fixed size due to neutral branching without mutation.
\end{enumerate}

If $m(\smallx)=0$, we put $\Omega^{\beta,\Gamma}_{\mbox{\tiny total mass}}F^{\ups}\big(\smallx\big)=0$.
Notice that this implies that also $\Omega F^{\ups}\big(\smallx\big)=0$ as $\ups(0,\boldsymbol{\cdot},\boldsymbol{\cdot})$ is assumed to be constant.\smallskip

Otherwise, if $m(\smallx)>0$ we
introduce the several parts step by step. We use the abbreviations:
 \begin{equation}
\label{hat:beta}
   \widehat{\beta}^{\smallx}:=\int\nu^{\smallx}(\mathrm{d}(\mr,\vk))\,\beta\big(\kappa_1\big),
\end{equation}
and
\begin{equation}
\label{hat:gammad}
   \widehat{\gamma}^{\mathrm{death/birth}}\big(m(\smallx)\big)
 :=
   \int\,\nu^{\smallx}(\mathrm{d}(\mr,\vk))\,
   \gamma^{\mathrm{death/birth}}\big(m(\smallx),r_{1,2},\kappa_1,\kappa_2\big).
\end{equation}

Moreover, we put
\begin{equation}
\label{hat:Gamma}
\begin{aligned}
   \widehat{\Gamma}\big(m(\smallx)\big)
 &:=
   \int\nu^{\smallx}(\mathrm{d}(\mr,\vk))\,\Gamma\big(m(\smallx),r_{1,2},\kappa_1,\kappa_2\big).
\end{aligned}
\end{equation}\smallskip

\noindent{\bf Step~1 (Changes in the total mass). }
%For all $F=F^{\ups}\in\Pi^{2,1,A}\cap\Pi_{\mathrm{finite}}$, put
Put
\begin{equation}
\label{y:002}
\begin{aligned}
   \Omega^{\beta,\Gamma}_{\mbox{\tiny total mass}}F^{\ups}\big(\smallx\big)
   &=:
    \Omega^{\Gamma}_{\mbox{\tiny competition}}F^{\ups}\big(\smallx\big)+\Omega^{\beta}_{\mbox{\tiny total mass fluctuation}}F^{\ups}\big(\smallx\big)
    +\Omega^{\beta}_{\mbox{\tiny mass flow}}F^{\ups}\big(\smallx\big)
    \\
  &=
    \widehat{\Gamma}(m(\smallx))\cdot m(\smallx)\cdot\tfrac{\partial}{\partial m}F^{\ups}\big(\smallx\big)+\widehat{\beta}^{\smallx}\cdot m(\smallx)\cdot\tfrac{\partial^2}{\partial^2 m}F^{\ups}\big(\smallx\big)
    \\
    &\;\;\;+2\cdot\int\nu^{\smallx}(\mathrm{d}(\mr,\vk))\,\tfrac{\partial}{\partial m}\ups\big(m(\smallx),\mr,\vk\big)\cdot\sum_{l\ge 1}\big(\beta(\kappa_l)-\widehat{\beta}^{\smallx}\big).
\end{aligned}
\end{equation}
In words, given the evolution of the sampling measure,
the total mass follows a branching diffusion with branching rate $\widehat{\beta}^{\smallx}$
and state dependent drift $\widehat{\Gamma}(m(\smallx))$.  In addition, as $\beta$ is trait dependent, changes in the sampling measure lead to a flow of mass which we would not see if the (natural) branching $\beta$ were a constant.\smallskip

\noindent{\bf Step~2 (Trait Mutation). }
%For all $F=F^\ups\in\Pi^{2,1,A}\cap\Pi_{\mathrm{finite}}$, put
Put
\begin{equation}
\label{y:004}
\begin{aligned}
   \Omega^{p,\beta,A}_{\mbox{\tiny trait mutation}}F^{\ups}\big(\smallx\big)
  &=
    p\cdot\int\nu^{\smallx}(\mathrm{d}(\mr,\vk))\,\sum_{l\ge 1}\beta(\kappa_{l})\cdot A^{(l)}\ups\big(m(\smallx),\mr,\vk\big),
\end{aligned}
\end{equation}
where $A^{(l)}$ acts on $\ups$ as the mutation operator $A$ on the function of the $l^{\mathrm{th}}$-trait-coordinate of $\ups$ (assuming that all other variables are kept constant). \smallskip

\noindent{\bf Step~3 (Growth of genetic distances). }
%For all $F_0\in\Pi^{2,1,A}\cap\Pi_{\mathrm{finite}}$, put
Put
\begin{equation}
\label{y:005}
\begin{aligned}
   \Omega^{p,\beta}_{\mbox{\tiny growth}}F^{\ups}\big(\smallx\big)
  &=
    p\cdot\int\nu^{\smallx}(\mathrm{d}(\mr,\vk))\,\sum_{1\le l_1<l_2}\big(\beta(\kappa_{l_1})+\beta(\kappa_{l_2})\big)\cdot\tfrac{\partial}{\partial r_{l_1,l_2}}\ups\big(m(\smallx),\mr,\vk\big).
\end{aligned}
\end{equation}
That is, the distance between two individuals, one from the clan $x_1$ and of current type $\kappa_1$ and the other one from the clan $x_2$ and of current type $\kappa_2$ grows at speed  $p\cdot(\beta(\kappa_1)+\beta(\kappa_2))$, which indeed models growing of the substitution distance.
\smallskip

\noindent{\bf Step~4 (Reweigh of the sampling measure with respect to $\Gamma$). }
%For all $F=F^\ups\in\Pi^{2,1,A}\cap\Pi_{\mathrm{finite}}$, put
Put
\begin{equation}
\label{y:006}
\begin{aligned}
   \Omega^{\Gamma}_{\mbox{\tiny $\Gamma$-reweigh}}F^{\ups}\big(\smallx\big)
  &=
   \int\nu^{\smallx}(\mathrm{d}(\mr,\vk))\,\sum_{l\ge 1}\big(\Gamma(m(\smallx),r_{1,l+1},\kappa_1,\kappa_{l+1}))-\widehat{\Gamma}(m(\smallx))\big)\cdot
   \ups\big(m,\tau_1(\mr,\vk)\big),
\end{aligned}
\end{equation}
where for $\ell\ge 1$,
$\tau_{\ell}$ denotes an index shift by $\ell\in\mathbb{N}$, that is,
\begin{equation}
\label{s:tauell}
   \tau_{\ell}\big(\mr,\vk\big):=\big((r_{i,j})_{l+1\le i<j},(\kappa_i)_{l+1\le i}\big).
\end{equation}
\smallskip

\noindent{\bf  Step~5 (Effect of neutral branching without mutation on the genealogy). }
%For all $F=F^\ups\in\Pi^{2,1,A}\cap\Pi_{\mathrm{finite}}$, put
Put
\begin{equation}
\label{y:007}
\begin{aligned}
   &\Omega^{\beta}_{\mbox{\tiny natural branching}}F^{\ups}\big(\smallx\big)
   \\
  &=:
     \Omega^{\beta}_{\mbox{\tiny $\beta$-reweigh}}F^{\ups}\big(\smallx\big)+\Omega^{\beta}_{\mbox{\tiny resample}}F^{\ups}\big(\smallx\big)
     \\
  &=
   \tfrac{1}{m(\smallx)}\int\nu^{\smallx}(\mathrm{d}(\mr,\vk))\,\sum_{1\le l}\big(\widehat{\beta}^{\smallx}-\beta(\kappa_l)\big)\ups\big(m(\smallx),\mr,\vk\big)
   \\
   &\;\;\;+
   \tfrac{1}{m(\smallx)}\int\nu^{\smallx}(\mathrm{d}(\mr,\vk))\sum_{1\le l_1,l_2}\Big\{ \beta(\kappa_{l_1})\Big(\Theta_{l_1,l_2}\ups\big(m(\smallx),\mr,\vk\big)
   -\ups\big(m(\smallx),\mr,\vk\big)\Big)+\Big(\widehat{\beta}^{\smallx}-\beta(\kappa_{l_2})\Big)
   \ups\big(m(\smallx),\mr,\vk\big)\Big\},
   \end{aligned}
\end{equation}
where for $1\le l_1 \neq l_2$, the {\em replacement map} $\Theta_{l_1,l_2}$ on ${\mathcal C}_b(\R_+\times \R_+^{\mathbb{N}\choose 2}\times K^{\mathbb{N}})$ sends a function to a new function by ``replacing its
$l_2^{\mathrm{nd}}$ argument by its $l_1^{\mathrm{st}}$, that is, $\Theta_{l_1,l_2}\ups=\ups\circ\vartheta_{l_1,l_2}$, where
\begin{equation}
\label{e:018_2}
   \vartheta_{l_1,l_2}\big(m(\smallx),(r_{i,j})_{1\le i<j},(\kappa_i)_{i\in\mathbb{N}}\big):=\big(m(\smallx),(r_{\theta_{l_1,l_2}(i),\theta_{l_1,l_2}(j)})_{1\le i<j},(\kappa_{\theta_{l_1,l_2}(i)})_{i\in\mathbb{N}}\big),
\end{equation}
and
\begin{equation}
\label{e:018_3}
   \theta_{l_1,l_2}(i):=\left\{\begin{array}{cc}i, & \mbox{ if }i\not=l_2, \\
   l_1, & \mbox{ if }i=l_2.\end{array}\right.
\end{equation}
\bigskip

\begin{remark}[Consistency] We point out that
\begin{equation}
\label{e:containsfinite}
   \Pi^{2,1,A}\cap\Pi_{\mathrm{finite}}\subseteq{\mathcal D}(\Omega).
\end{equation}
More precisely, summands associated with an index $l$ appearing at any of the right hand sides in (\ref{y:002})--(\ref{y:007}) are zero (and thus do not contribute to the series) if the function $\ups$ does not depend explicitly
on the metric entries $r_{i\wedge l,i\vee l}$, $i\in\mathbb{N}$,  and  the type entry $\kappa_l$.
%Notice that we have defined the operator $\Omega$ acting on polynomials of finite degree without %referring explicitly to the degree.
This is obvious to see for all parts of the generator except maybe $\Omega^\beta_{\mbox{\tiny resample}}$. We therefore want to present a short argument here.

Assume first that $l_1\in\mathbb{N}$ is such that
 the function $\ups$ does not depend on $r_{i\wedge l_1,i\vee l_1}$
for all $i\in\mathbb{N}$ and also does not depend on $\kappa_{l_1}$. In that case,
\begin{equation}
\label{e:014nol1}
\begin{aligned}
   &\int\nu^{\smallx}(\mathrm{d}(\mr,\vk))\,\Big\{\beta(\kappa_{l_1})\cdot\Big(\Theta_{l_1,l_2}\ups\big(m(\smallx),\mr,\vk\big)-\ups\big(m(\smallx),\mr,\vk\big)\Big)+
   \Big(\widehat{\beta}^{\smallx}-\beta(\kappa_{l_2})\Big)\cdot
   \ups\big(m(\smallx),\mr,\vk\big)\Big\}
   \\
   &=\int\nu^{\smallx}(\mathrm{d}(\mr,\vk))\,\Big\{\beta(\kappa_{l_1})\cdot\Theta_{l_1,l_2}\ups\big(m(\smallx),\mr,\vk\big)-\beta(\kappa_{l_2})\cdot
   \ups\big(m(\smallx),\mr,\vk\big)\Big\}
   \\&\;\;\;
   -\widehat{\beta}^{\smallx}\int\nu^{\smallx}(\mathrm{d}(\mr,\vk))\,\ups\big(m(\smallx),\mr,\vk\big)+\widehat{\beta}^{\smallx}\int\nu^{\smallx}(\mathrm{d}(\mr,\vk))\,\ups\big(m(\smallx),\mr,\vk\big)
   \\
   &=0.
\end{aligned}
\end{equation}

On the other hand, assume  that
 $l_2\in\mathbb{N}$ is such that
 the function $\ups$ does not depend on $r_{i\wedge l_2,i\vee l_2}$
for all $i\in\mathbb{N}$ and also does not depend on $\kappa_{l_2}$. Then once more
\begin{equation}
\label{e:014nol2}
\begin{aligned}
   &\int\nu^{\smallx}(\mathrm{d}(\mr,\vk))\,\Big\{\beta(\kappa_{l_1})\cdot\Big(\Theta_{l_1,l_2}\ups\big(m(\smallx),\mr,\vk\big)-\ups\big(m(\smallx),\mr,\vk\big)\Big)+
   \Big(\widehat{\beta}^{\smallx}-\beta(\kappa_{l_2})\Big)\cdot
   \ups\big(m(\smallx),\mr,\vk\big)\Big\}
   \\
   &=\int\nu^{\smallx}(\mathrm{d}(\mr,\vk))\,\beta(\kappa_{l_1})\cdot\Big(\ups\big(m(\smallx),\mr,\vk\big)-\ups\big(m(\smallx),\mr,\vk\big)\Big) 
   \\&\;\;\;+\widehat{\beta}^{\smallx}\int\nu^{\smallx}(\mathrm{d}(\mr,\vk))\,\ups\big(m(\smallx),\mr,\vk\big)-\widehat{\beta}^{\smallx}\int\nu^{\smallx}(\mathrm{d}(\mr,\vk))\,\ups\big(m(\smallx),\mr,\vk\big)
   \\
   &=0.
\end{aligned}
\end{equation}
This shows that indeed $\Omega^\beta_{\mbox{\tiny resample}}F^\ups$
reduces to a sum with finitely many summands and is thus well-defined whenever $\ups$ has finite degree.
\hfill$\qed$
\end{remark}\sm

\iffalse{Consider the domain
\begin{equation}
\label{e:domm}
\begin{aligned}
   &{\mathcal D}(\Omega)
 \\
  &:=\mbox{algebra generated by those $F^\ups\in\Pi^{2,1,A}$ with $\sup_{m>0}\tfrac{\ups(m,\mr,\vk)}{m^n}<\infty$
   for all $n,(\mr,\vk)$}
   \\
   &\hspace{7cm}{\color{red} \mbox{ and }\ups(0,\boldsymbol{\cdot},\boldsymbol{\cdot})=\tfrac{\partial}{\partial m}\ups(m,\boldsymbol{\cdot},\boldsymbol{\cdot})|_{m=0}\equiv 0}.
\end{aligned}
\end{equation}

\todo{check the conditions: we need ${\mathcal D}(\Omega):=\{F^\ups\in\Pi^{2,1,A}\}$ such that
$\ups\in\tilde{mathcal H}$ ensuring uniform convergence of the generators; also every function which can be approximated is good too}
}\fi

We are now in a position to state our second main result.
%The following will be proved in Section~\ref{S:proofs}.
\begin{theorem}[The martingale problem] Let for each $N\in\mathbb{N}$, ${\mathcal X}^N$ be the  tree-valued $(\varsig_N,\zeta_N,\alpha_N)$-trait-dependent branching dynamics with mutation and competition
such that the assumptions of Theorem~\ref{T:tightness} hold. If ${\mathcal X}_0^N$ converges weakly to ${\mathcal X}_0$ in $\mathbb{M}^K$ with $\mathbb{E}[(m({\mathcal X}_0))^3]<\infty$, then
any limit process ${\mathcal X}$ satisfies the $(\Omega,{\mathcal D}(\Omega))$-martingale problem.
\iffalse{
Assume that the mutation operator  $(A,{\mathcal D}(A))$
on ${\mathcal C}_b(K)$ is such that ${\mathcal D}(A)$ is an algebra and the $(A,{\mathcal D}(A))$-martingale problem is well-posed.
Assume further that $\beta:K\to\R_+$ satisfies Assumptions~\ref{ass:002} and~\ref{ass:007c} and in addition that $\beta\in{\mathcal D}(A)$. Moreover, assume that $\Gamma:\R_+\times\R_+\times K\times K\to\R$ satisfies Assumptions~\ref{ass:007b} and~\ref{ass:007a} and that in addition $\Gamma(m,r,\boldsymbol{\cdot},\kappa_2),\Gamma(m,r,\kappa_1,\boldsymbol{\cdot})\in{\mathcal D}(A)$ for all $m,r\ge 0$ and $\kappa_1,\kappa_2\in K$.  Then for all
$P\in{\mathcal M}_1(\mathbb{M}^{K})$ with $\int P(\mathrm{d}\smallx)(m(\smallx))^2<\infty$, the
$(P,\Omega,{\mathcal D}(\Omega))$-martingale problem is well-posed.}\fi
\label{T:mp}
\end{theorem}\sm

\begin{remark}[Relation with other tree-valued dynamics] Several related models of evolving genealogies have been considered so far in the literature.
\begin{enumerate}
\item {\bf Tree-valued Fleming-Viot with mutation and selection.}
In \cite{DGP12} the authors consider the evolving genealogies of a trait-dependent individual based particle system with mutation and selection, where the total mass is fixed, say $m(\smallx_t)\equiv 1$ (and thus $\Omega^{\beta,\Gamma}_{\mbox{\tiny total mass}}\equiv 0$), and the natural branching rate does not depend on the trait, that is, $\beta\equiv\widehat{\beta}$ (and thus $\Omega^{\beta}_{\mbox{\tiny $\beta$-reweigh}}\equiv 0$), and their single mutation is a continuous state and continuous time Markov chain on $K$, that is, of the form
     \begin{equation}
    \label{e:DGPA}
        Af(\kappa):=\vartheta\cdot\int_K \big(a(\kappa,\mathrm{d}\kappa')-\delta_\kappa(\mathrm{d}\kappa')\big)f(\kappa')
   \end{equation}
   (which is the limit of rare mutations).
    Moreover, their genealogical distances grow with time no matter whether mutation occurs or not
    (might be interpreted as $p=1$).
    In this particular case,
  %\begin{equation}
   % \label{e:DGPGamma}
        $\Gamma\big(m,r_{1,2},\kappa_1,\kappa_2\big)\equiv \tilde{\Gamma}\big(r_{1,2},\kappa_1,\kappa_2\big)$,
   %\end{equation}
   for some bounded $\tilde{\Gamma}:\R_+\times K\times K$,
    and (\ref{y:001}) reduces to $
     %\begin{equation}
    %\label{e:DGP}
     %   \begin{aligned}
      \Omega
  %&
    =
      \Omega^{\widehat{\beta},1,{\frac{1}{\widehat{\beta}}A}}_{\mbox{\tiny trait mutation}}
    +
      \Omega^{1,{1}}_{\mbox{\tiny growth}}
    +
      \Omega^{\tilde{\Gamma}}_{\mbox{\tiny $\Gamma$-reweigh}}
    +
      \Omega^{\widehat{\beta}}_{\mbox{\tiny resample}}$.
%\end{aligned}
 %   \end{equation}
\item {\bf Tree-valued (neutral) Fleming-Viot (without mutation). }
    In the particular case without selection and mutation, that is, where $\Gamma\equiv 0$ (and thus $\Omega^{\Gamma}_{\mbox{\tiny $\Gamma$-reweigh}}\equiv 0$) and $A\equiv\mathrm{Id}$ (and thus $\Omega^{p,\beta,A}_{\mbox{\tiny trait mutation}}\equiv 0$), this reduces to
      % \begin{equation}
    %\label{e:GPW}
     %   \begin{aligned}
      $\Omega
  %&:
    =
     \Omega^{1,{1}}_{\mbox{\tiny growth}}
    +
      \Omega^{\widehat{\beta}}_{\mbox{\tiny resample}}$,
%\end{aligned}
 %   \end{equation}
    which is the operator corresponding to the tree-valued Fleming-Viot diffusion constructed in \cite{GPW13}.
\item   {\bf Historical processes with mutation and selection. }
 In \cite{MT12} consider a trait-dependent individual based branching model with mutation and selection is considered, where the selection takes the trait history of an individual into account. If we restrict their model to the case, where the dependence on the historical trait path is only through the current trait, this results in a measure-valued process (and thus $\Omega^{p,\beta}_{\mbox{\tiny growth}}\equiv 0$) with the particular set-up
       \begin{equation}
    \label{e:MTGamma}
       \Gamma\big(m,r_{1,2},\kappa_1,\kappa_2\big)=b(\kappa_2)-d(\kappa_2)-m\cdot U(\kappa_1,\kappa_2)
    \end{equation}
    and
       \begin{equation}
    \label{e:MTA}
       Af(\kappa):=\tfrac{\sigma^2}{2}f''(\kappa).
    \end{equation}
  Thus (\ref{y:001}) reads as follows: $\Omega :=
      \Omega^{\beta,\Gamma}_{\mbox{\tiny total mass}}
    +
      \Omega^{p,\beta,A}_{\mbox{\tiny trait mutation}}
   +
      \Omega^{\Gamma}_{\mbox{\tiny $\Gamma$-reweigh}}
   +
      \Omega^{\beta}_{\mbox{\tiny natural branching}}$,
     where the generator acts for instance on functions $F^{\ups}$ of the particular form
    \begin{equation}
    \label{e:MTups}
     \ups(m,(r_{i,j})_{1\le i<j\le n},(\kappa_i)_{1\le i\le n}):=m^n\cdot f\big(\vk\big).
  \end{equation}

\iffalse{ The measure-valued process can be characterized as the unique solution $(X_t)_{t\ge 0}$ of the following square integrable martingale problem:
\begin{equation}
   M_s:=\langle\mu_t,f\rangle-\langle\mu_0,f\rangle-\int^t_0\mathrm{d}s\big\langle\mu_s,\big(p\beta(\boldsymbol{\cdot})Af(\boldsymbol{\cdot})+
   \Gamma(m(\smallx_s),\boldsymbol{\cdot},\bar{\mu}_s)
   f\big)\big\rangle
\end{equation}
with quadratic variation
\begin{equation}
\label{e:MTvariation}
   \langle M\rangle_t=2\int^t_0\mathrm{d}s\,\langle\mu_s,\beta(\boldsymbol{\cdot})f^2(\boldsymbol{\cdot})\rangle.
\end{equation}}\fi
The {latter} model has been extended in \cite{Kli14} by considering general mutation operators on Polish trait spaces.\hfill$\qed$
\end{enumerate}
\end{remark}\sm

% ======================================================================

 \section{Uniform moment bounds for the discrete models}
 \label{S:proofdiscrete}
 Throughout this section we assume that Assumptions~\ref{ass:002}, \ref{ass:007b} and \ref{ass:003} hold.
 Recall the  tree-valued $(\varsig,\zeta,\alpha)$-trait-dependent branching dynamics with mutation and competition, ${\mathcal X}^{(\varsig,\zeta,\alpha)}$,
 as defined
 in Subsection~\ref{Sub:branchmodel} (Definition~\ref{Def:002}) by means of the natural branching rate $\beta:K\to\R_+$, the competition birth rate $\gamma^{\mathrm{death}}:\,\R_+\times\R_+\times K^2\to\R_+$,
 the competition death rate $\gamma^{\mathrm{death}}:\,\R_+\times\R_+\times K^2\to\R_+$, and the mutation kernel $\alpha$.
 Its generator $\Omega^{(\varsig,\zeta,\alpha)}$ acts on bounded functions $F$ on $\mathbb{M}^{K,(\ell,\zeta)}$ as follows (cf. (\ref{s:002aa})--(\ref{s:004})):
 \begin{equation}
\label{e:OmegaN}
\begin{aligned}
  &\Omega^{(\varsig,\zeta,\alpha)} F\big(\smallx\big)
 \\
  &=:\Big(\Omega^{(\ell,\zeta,\alpha)}_{\mbox{\tiny death}}+\Omega^{(\ell,\zeta,\alpha)}_{\mbox{\tiny birth}}\Big)F\big(\smallx\big)
  \\
  &=
   \tfrac{m(\smallx)}{\zeta}\int_{(X\times K)^2}\bar{\mu}^{\otimes 2}(\mathrm{d}((x_1,\kappa_1),(x_2,\kappa_2)))\big\{\tfrac{\beta(\kappa_2)}{\zeta}+\gamma^{\mathrm{death}}(m(\smallx),r(x_1,x_2),\kappa_1,\kappa_2)\big\}
   \\
   &\hspace{2cm}\cdot
   \Big(F\big((X,r,\mu-\zeta\cdot\delta_{(x_2,\kappa_2)})\big)-F\big(\smallx\big)\Big)
    \\
  &+\tfrac{m(\smallx)}{\zeta}\int_{(X\times K)^2}\bar{\mu}^{\otimes 2}
   (\mathrm{d}((x_1,\kappa_1),(x_2,\kappa_2)))\big\{\tfrac{\beta(\kappa_2)}{\zeta}+\gamma^{\mathrm{birth}}(m(\smallx),r(x_1,x_2),\kappa_1,\kappa_2)\big\}
   \int_K\widehat{\alpha}\big(\kappa_2,\mathrm{d}\tilde{\kappa}_2\big)
     \\
   &\hspace{2cm}\cdot
   \Big(F\big((X\uplus\{z\},r^{(x_2,z),\varsig},\mu+\zeta\cdot\delta_{(z,\tilde{\kappa}_2)})\big)-
   F\big(\smallx\big)\Big).
\end{aligned}
\end{equation}
\smallskip

In the proof of tightness we will make use of the following uniform moment bounds.
\begin{proposition}[Uniform moment bounds]
Let ${\mathcal X}:=({\mathcal X}^{(\ell,\zeta,\alpha)})_{t\in[0,\infty)}$
be the $(\ell,\zeta,\alpha)$-trait-dependent branching dynamics with mutation and competition. Fix $q\in\mathbb{N}$.
\begin{itemize}
\item[(i)]
Then for all $(\ell,\zeta,\alpha)$ with $\zeta\le 1$, and all $t>0$,
\begin{equation}
\label{e:Gronwall}
\begin{aligned}
  \sup_{s\in[0,t]}\E\big[\big(m({\mathcal X}^{(\ell,\zeta,\alpha)}_s)\big)^q\big]
   &\le
     \big(1+\E\big[\big(m({\mathcal X}^{(\ell,\zeta,\alpha)}_0)\big)^q\big]\big)\cdot e^{C_q\cdot(2\bar{\beta}+\bar{\gamma}_b)\cdot t}
\end{aligned}
\end{equation}
for some $C_q>0$.
\item[(ii)] Assume that in addition Assumption~\ref{ass:007a} holds.
If the initial masses satisfy
\begin{equation}
\label{ass-ini-cond}
 \sup_{\zeta\in(0,1]}\E\big[\big(m({\mathcal X}^{(\ell,\zeta,\alpha)}_0)\big)^{2q+1}\big] < \infty,
\end{equation}
then for all $t>0$, there is a constant $C_{q,\bar{\beta},\bar{\gamma_b},\bar{\gamma}_d,t}\in(0,\infty)$ such that for all $(\ell,\zeta,\alpha)$ with $\zeta\le 1$,
\begin{equation}
\label{e:t-moment-bounds-1}
  \E\big[\sup_{s\in[0,t]} \big(m({\mathcal X}^{(\ell,\zeta,\alpha)}_s)\big)^q\big]\le C_{q,\bar{\beta},\bar{\gamma_b},\bar{\gamma}_d,t}.
\end{equation}
\item[(iii)]
Let $q=1$ and $\gamma^{\mathrm{death}}(m,r,\kappa,\kappa') \equiv 0$, that is $\overline{\gamma}_d=0$. Under the assumptions of (ii), including (\ref{ass-ini-cond}), for all $t>0$ and $\delta, \epsilon>0$ there exists $m_0>0$ small enough such that for all $(\ell,\zeta,\alpha)$ with $\zeta\le 1$,
\begin{equation}
\label{e:t-moment-bounds-1iii}
  \mathbb{P}\big( \sup_{s \in [0,t]} m({\mathcal X}^{(\ell,\zeta,\alpha)}_s) \geq \delta \big)
  \le \epsilon
\end{equation}
if $\E\big[m({\mathcal X}^{(\ell,\zeta,\alpha)}_0\big] \leq m_0$.
\end{itemize}
\label{P:005}
\end{proposition}\sm

\begin{proof}
The proof is inspired by \cite[Step~2 of the proof of Theorem~5.6]{FM04}.
Fix $q\ge 1$, and put for $L\ge 0$,
\begin{equation}
\label{e:Fq}
   F^{q,L}\big(\smallx\big)
 :=
   \big(m(\smallx)\wedge L\big)^q.
\end{equation}

Then $F^{q,L}\in{\mathcal D}(\Omega^{(\varsig,\zeta,\alpha)})$, and
\begin{equation}
\label{e:omegaN_g}
\begin{aligned}
  &\Omega^{(\varsig,\zeta,\alpha)} F^{q,L}\big(\smallx\big)
   \\
   &=
   \zeta^{-1}\cdot m(\smallx)\cdot\big(\zeta^{-1}\cdot\widehat{\beta}^{\smallx}+\widehat{\gamma}^{\mathrm{death}}(m(\smallx))\big)\cdot\Big(\big((m(\smallx)-\zeta)\wedge L\big)^q-\big(m(\smallx)\wedge L\big)^q\Big)
   \\
  &+
    \zeta^{-1}\cdot m(\smallx)\cdot\big(\zeta^{-1}\cdot\widehat{\beta}^{\smallx}+\widehat{\gamma}^{\mathrm{birth}}(m(\smallx))\big)\cdot\Big(\big((m(\smallx)+\zeta)\wedge L\big)^q-\big(m(\smallx)\wedge L\big)^q\Big).
\end{aligned}
\end{equation}
(recall $\widehat{\beta}^{\smallx}$ from (\ref{hat:beta}) and $\widehat{\gamma}^{\mathrm{death/birth}}$ from (\ref{hat:gammad})).

Thus, for all $t\ge 0$,
\begin{equation}
\label{e:mp-g}
   F^{q,L}\big({\mathcal X}^{(\ell,\zeta,\alpha)}_t\big)=F^{q,L}\big({\mathcal X}^{(\ell,\zeta,\alpha)}_0\big)+\int_0^t \mathrm{d}s\, \Omega^{(\varsig,\zeta,\alpha)}F^{q,L}\big({\mathcal X}^{(\ell,\zeta,\alpha)}_s\big)+M_t^{q,L},
\end{equation}
where $(M^{q,L}_{t\wedge \tau_M})_{t\ge 0}$ with
\begin{equation}
\label{e:tauM}
   \tau_M:=\inf\big\{t\ge 0:\,m({\mathcal X}^{(\ell,\zeta,\alpha)}_t)\ge M\big\},\hspace{.2cm}M>0
\end{equation}
is a c\`adl\`ag $L^2$-martingale starting from $0$ and with quadratic variation
\begin{equation}
\label{e:mart-g-0}
  \langle M\rangle_{t\wedge\tau_M}
 =
   \int_0^{t\wedge\tau_M} \mathrm{d}s \Big(\Omega^{(\varsig,\zeta,\alpha)}(F^2)-2F\Omega^{(\varsig,\zeta,\alpha)}F\Big)
  \big({\mathcal X}^{(\ell,\zeta,\alpha)}_s\big).
\end{equation}

{\em (i).} To obtain an upper bound for $F^{q,L}\big({\mathcal X}^{(\varsig,\zeta,\alpha)}_{t\wedge\tau_M}\big)$, we can drop the competition death term in (\ref{e:omegaN_g}) respectively (\ref{e:mp-g}), and use that with $C_q:=2^q-1$ for all $m>0$ and $\zeta\in(0,1]$ (note that the case $q=1$ is trivial),
\begin{equation}
\label{e:mq}
\begin{aligned}
   m\cdot\big|(m+\zeta)^q - 2 m^q +(m-\zeta)^q \big|
   &=2m\zeta^2\cdot\big|\sum_{k=0,q-k~\mbox{\tiny even}}^{q-2}{q\choose k}m^k\zeta^{q-k-2}\big|
   \\
   &\leq 2{C_q}\cdot\zeta^2\cdot(1+m^q)
   \end{aligned}
\end{equation}
as well as
\begin{equation}
\label{e:mq1}
   m\cdot\big|(m+\zeta)^q - m^q\big| \leq C_q\cdot\zeta\cdot \big(1+m^{q}\big).
\end{equation}

Therefore by Assumptions~\ref{ass:002} and~\ref{ass:007b}, for all $M>0$ and
 all $L\ge M+\zeta$,
\begin{equation}
\begin{aligned}
\label{e:bound-E-m_q0}
  &F^{q,L}\big({\mathcal X}^{(\varsig,\zeta,\alpha)}_{t\wedge\tau_M}\big)
  \\
  &\le
     F^{q,L}\big({\mathcal X}^{(\varsig,\zeta,\alpha)}_0\big)+M_{t\wedge\tau_M}^{q,L}
   \\
   &\,+
    \zeta^{-2}\cdot\int_0^{t\wedge\tau_M}\mathrm{d}s\,m({\mathcal X}^{(\ell,\zeta,\alpha)}_s)\cdot\widehat{\beta}^{{\mathcal X}^{(\ell,\zeta,\alpha)}_s}\cdot\Big(\big(m({\mathcal X}^{(\ell,\zeta,\alpha)}_s)-\zeta\big)^q+\big(m({\mathcal X}^{(\ell,\zeta,\alpha)}_s)+\zeta\big)^q-2\big(m({\mathcal X}^{(\ell,\zeta,\alpha)}_s)\big)^q\Big)
   \\
   &\,+
   \zeta^{-1}\cdot\int_0^{t\wedge\tau_M}\mathrm{d}s\,m({\mathcal X}^{(\ell,\zeta,\alpha)}_s)\cdot\widehat{\gamma}^{\mathrm{birth}}(m({\mathcal X}^{(\ell,\zeta,\alpha)}_s))\cdot
   \Big(\big(m({\mathcal X}^{(\ell,\zeta,\alpha)}_s)+\zeta)\big)^q-\big(m({\mathcal X}^{(\ell,\zeta,\alpha)}_s)\big)^q\Big)
   \\
   &\le
     F^{q,L}\big({\mathcal X}^{(\varsig,\zeta,\alpha)}_0\big)+M_{t\wedge\tau_M}^{q,L}
     \\
   &\;+
     C_q\cdot\int^{t\wedge\tau_M}_0\mathrm{d}s\,\Big\{2\bar{\beta}\cdot\big((m({\mathcal X}^{(\ell,\zeta,\alpha)}_s))^q+1\big)+
     \bar{\gamma}_b\cdot\big((m({\mathcal X}^{(\ell,\zeta,\alpha)}_s))^q+1\big)\Big\}.
\end{aligned}
\end{equation}

Hence by monotone convergence as $L\uparrow\infty$, for all $M$,
\begin{equation}
\begin{aligned}
\label{e:bound-E-m_q}
  &\E\big[\big(m({\mathcal X}^{(\ell,\zeta,\alpha)}_{t\wedge\tau_M})\big)^q\big]
  \\
   &\le
       \E\big[\big(m({\mathcal X}^{(\ell,\zeta,\alpha)}_0)\big)^q\big]
   +
     C_q\cdot\big(2\bar{\beta}+\bar{\gamma}_b\big)\cdot\int^{t}_0\mathrm{d}s\,\big\{\E\big[\big(m({\mathcal X}^{(\ell,\zeta,\alpha)}_{s\wedge\tau_M})\big)^{q}\big]+1\big\}.
\end{aligned}
\end{equation}

We therefore obtain (\ref{e:Gronwall}) from Gronwall's lemma and subsequently taking $M \rightarrow \infty$. \smallskip

{\em (ii).} Now we move the supremum over time inside the expectation.
Using once more (\ref{e:bound-E-m_q0}),  for all $L\ge M+\zeta$,
\begin{equation}
\begin{aligned}
\label{e:bound-E-m_q2}
  \sup_{s\in[0,t\wedge\tau_M]}F^{q,L}\big({\mathcal X}^{(\ell,\zeta,\alpha)}_s\big)
   &\le
     F^{q,L}\big({\mathcal X}_0\big)+\sup_{s\in[0,t \wedge\tau_M]}\big|M_s^{q,L}\big|
     \\
   &\;+
    C_q\cdot\big(2\bar{\beta}+\bar{\gamma}_b\big)\cdot
    \int^{t\wedge\tau_M}_0\mathrm{d}s\,\big((m({\mathcal X}^{(\ell,\zeta,\alpha)}_s))^q+1\big).
\end{aligned}
\end{equation}

Taking expectations,
\begin{equation}
\begin{aligned}
\label{e:bound-E-m_q4}
   &\E\big[\sup_{s\in[0,t\wedge\tau_M]}\big(m({\mathcal X}^{(\ell,\zeta,\alpha)}_s)\big)^q\big]
   \\
   &\le
   \E\big[\big(m({\mathcal X}^{(\ell,\zeta,\alpha)}_0)\big)^q\big]+\E\big[\sup_{s\in[0,t\wedge\tau_M]}\big| M^{q,L}_s\big|\big]+C_q\cdot\big(2\bar{\beta}+\bar{\gamma}_b\big)\cdot
   \int^t_0\mathrm{d}s\,\mathbb{E}\big[\big((m({\mathcal X}^{(\ell,\zeta,\alpha)}_s))^q+1\big)\cdot 1\{s\le\tau_M\}\big]
   \\
  &\le
    \E\big[\big(m({\mathcal X}^{(\ell,\zeta,\alpha)}_0)\big)^q\big]+\E\big[\sup_{s\in[0,t\wedge\tau_M]}\big| M^{q,L}_s\big|\big]+C_q\cdot\big(2\bar{\beta}+\bar{\gamma}_b\big)\cdot
   \int^t_0\mathrm{d}s\,\mathbb{E}\big[\big((m({\mathcal X}^{(\ell,\zeta,\alpha)}_{s\wedge\tau_M}))^q+1\big)\big].
 \end{aligned}
\end{equation}

By a Burkholder-Davis-Gundy inequality and Cauchy-Schwarz's inequality, there is a $C<\infty$ such that
\begin{equation}
\begin{aligned}
\label{e:Doob}
   \E\big[\sup_{s \in [0,t\wedge\tau_M]}\big|M_s^{q,L}\big|\big]
  &\le
     C\cdot\E\big[\langle M_{\boldsymbol{\cdot}}^{q,L}\rangle^{\frac{1}{2}}_{t\wedge\tau_M}\big]
   \le
     C\cdot\big(\E\big[\langle M_{\boldsymbol{\cdot}}^{q,L}\rangle_{t\wedge\tau_M}\big]\big)^{\frac{1}{2}}.
\end{aligned}
\end{equation}

We conclude with (\ref{e:mart-g-0}) that for all $L\ge M+\zeta$,
\begin{equation}
\label{e:006}
\begin{aligned}
   &\E\big[\langle M_{\boldsymbol{\cdot}}^{q,L}\rangle_{t\wedge\tau_M}\big]
   \\
   &\leq
    \int^t_0\mathrm{d}s\,\E\big[\tfrac{m({\mathcal X}^{(\ell,\zeta,\alpha)}_s)}{\zeta}\Big(\big(m({\mathcal X}^{(\ell,\zeta,\alpha)}_s)-\zeta\big)^{q}+\big(m({\mathcal X}^{(\ell,\zeta,\alpha)}_s)\big)^q\Big)^2
   \cdot\big(\tfrac{\widehat{\beta}^{{\mathcal X}^{(\ell,\zeta,\alpha)}_s}}{\zeta}+\widehat{\gamma}^{\mathrm{death}}(m({\mathcal X}^{(\ell,\zeta,\alpha)}_s))\big);\,s\le\tau_M
   \big]
   \\
  &\;\;+ \int^t_0\mathrm{d}s\,\E\big[\tfrac{m({\mathcal X}^{(\ell,\zeta,\alpha)}_s)}{\zeta}\Big(\big(m({\mathcal X}^{(\ell,\zeta,\alpha)}_s)+\zeta\big)^{q}+\big(m({\mathcal X}^{(\ell,\zeta,\alpha)}_s)\big)^q\Big)^2
   \cdot\big(\tfrac{\widehat{\beta}^{{\mathcal X}^{(\ell,\zeta,\alpha)}_s}}{\zeta}+\widehat{\gamma}^{\mathrm{birth}}(m({\mathcal X}^{(\ell,\zeta,\alpha)}_s))\big);\,s\le\tau_M
   \big]
  \\
  &\leq
     C\cdot\big(2\bar{\beta}+\bar{\gamma}_b+\bar{\gamma}_d\big)\cdot\int^t_0\mathrm{d}s\,\big(1+ \E\big[(m({\mathcal X}^{(\ell,\zeta,\alpha)}_{s}))^{2q+1}\big]\big),
\end{aligned}
\end{equation}
which gives the claim after combining (\ref{e:Gronwall}), (\ref{e:bound-E-m_q4}), (\ref{e:Doob}) and (\ref{e:006}). \smallskip

{\em (iii).} By (ii) we already know that the integrals in question are well-defined. For $q=1$, instead of (\ref{e:mq}) and (\ref{e:mq1}), use that
\begin{equation}
\label{e:mqiii}
  m\cdot\big|(m+\zeta) - 2 m +(m-\zeta) \big| = 0
  \ \mbox{ and } \ 
   m\cdot\big|(m+\zeta) - m\big| = m \zeta
\end{equation}
to conclude as in (\ref{e:bound-E-m_q0}) that 
\begin{equation}
\label{e:bound-E-m_qiii}
  m({\mathcal X}^{(\ell,\zeta,\alpha)}_{t\wedge\tau_M}) \wedge L
  =
       m({\mathcal X}^{(\ell,\zeta,\alpha)}_0) \wedge L
   +
       M_{t\wedge\tau_M}^{1,L}
   +
     \int^{t\wedge\tau_M}_0\mathrm{d}s\,m({\mathcal X}^{(\ell,\zeta,\alpha)}_{s}) \cdot \widehat{\gamma}^{\mathrm{birth}}(m({\mathcal X}^{(\ell,\zeta,\alpha)}_s)).
\end{equation}
Hence, $m({\mathcal X}^{(\ell,\zeta,\alpha)}_{t\wedge\tau_M}) \wedge L$ is a non-negative submartingale and Gronwall's lemma yields $\E\big[ m({\mathcal X}^{(\ell,\zeta,\alpha)}_t) \big] \leq \E\big[ m({\mathcal X}^{(\ell,\zeta,\alpha)}_0) \big] \cdot e^{\overline{\gamma}_b t}$. In combination with Doob's inquality we obtain for all $\delta>0$,
\begin{equation}
\label{e:dmi}
  \mathbb{P}\big( \sup_{s \in [0,t]} m({\mathcal X}^{(\ell,\zeta,\alpha)}_s) \geq \delta \big)
  \le \delta^{-2} \E\big[ \big( m({\mathcal X}^{(\ell,\zeta,\alpha)}_t) \big)^2 \big].
\end{equation}
We use (\ref{e:mart-g-0}) again, together with $m \big( (m\pm\zeta)^2-m^2-2m((m\pm\zeta)-m) \big) = m \zeta^2$, to get 
\begin{equation}
\label{e:006iii}
  \E\big[\langle M_{\boldsymbol{\cdot}}^{q,1}\rangle_{t\wedge\tau_M}\big]
  \le
  \big(2\bar{\beta}+\bar{\gamma}_b\big)\cdot\int^t_0\mathrm{d}s\,\mathbb{E}\big[m({\mathcal X}^{(\ell,\zeta,\alpha)}_s)\big]
  \le 
  \big(2\bar{\beta}+\bar{\gamma}_b\big)\cdot m_0\cdot\int^t_0\mathrm{d}s\,e^{\overline{\gamma}_b s}.
\end{equation}
Now apply It\^o's formula to see that
\begin{equation}
\label{e:dmi-ito}
  \E\big[ \big( m({\mathcal X}^{(\ell,\zeta,\alpha)}_t) \big)^2 \big]
  \le m_0^2
  +
  2\overline{\gamma}_b\cdot\int_0^t\mathrm{d}s\,\E\big[ \big( m({\mathcal X}^{(\ell,\zeta,\alpha)}_t) \big)^2 \big] 
  +
  \big(2\bar{\beta}+\bar{\gamma}_b\big)\cdot m_0\cdot\int^t_0\mathrm{d}s\,e^{\overline{\gamma}_b s}.
\end{equation}
Another application of Gronwall's lemma yields
\begin{equation}
  \E\big[ \big( m({\mathcal X}^{(\ell,\zeta,\alpha)}_t) \big)^2 \big]
  \le
  m_0\cdot\Big( m_0 + \big(2\bar{\beta}+\bar{\gamma}_b\big)\cdot\int^t_0\mathrm{d}s\,e^{\overline{\gamma}_b s} \Big) \cdot e^{2 \bar{\gamma}_b t}.
\end{equation}
Use $m_0$ small enough in (\ref{e:dmi}) to conclude the claim.
\end{proof}\sm

 \section{Uniform convergence of generators}
 \label{S:generators}
 Recall the $(\ell,\zeta,\alpha)$-trait-dependent branching with mutation and competition from Definition~\ref{Def:002}, its state space $\mathbb{M}^{K,(\ell,\zeta)}$ from (\ref{s:001}), and its generator {$\Omega^{(\ell,\zeta,\alpha)}=\Omega^{(\ell,\zeta,\alpha)}_{\mbox{\tiny birth}}
 +\Omega^{(\ell,\zeta,\alpha)}_{\mbox{\tiny death}}$} acting on bounded measurable functions from (\ref{e:OmegaN}).
 In the following we abbreviate
 \begin{equation}
\label{e:OmeggaN}
    \Omega_N:=\Omega^{(\frac{1}{N},\frac{1}{N},\alpha_N)}
    \end{equation}
    and
    \begin{equation}
    \label{e:OmmegaN}
    \mathbb{M}^K_N:=\mathbb{M}^{K,(\frac{1}{N},\frac{1}{N})}
 \end{equation}
 with a family $\{\alpha_N ;\,N\in\mathbb{N}\}$ of mutation operators satisfying Assumption~\ref{ass:004}.
Furthermore, recall the tree-valued trait-dependent branching with mutation and competition from Definition~\ref{Def:004}, its state space $\mathbb{M}^K$ from (\ref{eq:defM}), and its generator $(\Omega,{\mathcal D}(\Omega))$ from (\ref{e:DOmega}) and (\ref{y:001}).

As a first step in proving tightness we are proving the uniform convergence of the generators in this section.
Put
\begin{equation}
\label{e:tildePi}
\begin{aligned}
   \widetilde{\Pi}
   &:=\big\{F^{\ups}\in\Pi^{2,1,A}\cap\Pi_{\mathrm{finite}}\mbox{ of the form }\ups(m,\mr,\vk)=g(m)\cdot \phi\big((r_{i,j})_{1\le i<j\le n}\big)\cdot f\big((\kappa_i)_{1\le i\le n}\big)\mbox{ with}
   \\
   &\qquad n\in\mathbb{N}_0, g\in{\mathcal C}_b^3(\R_+)\mbox{ such that }g'(0)=0,
   \limsup_{m\to\infty}m\tilde{\gamma}(m)|g'(m)|<\infty,
   \\
   &\qquad
   \limsup_{m\to\infty}m|g''(m)|<\infty,
     \limsup_{m\to\infty}(1 \vee \tilde{\gamma}(m))|g(m)|<\infty \mbox{ and }
     \\
    &\qquad \lim_{\epsilon \downarrow 0} \limsup_{m\to\infty}m\cdot (1\vee\tilde{\gamma}(m)) \max_{\xi \in [m-\epsilon,m+\epsilon]} |g'''(\xi)|<\infty\big\},
\end{aligned}
\end{equation}
where $\tilde{\gamma}$ is as in Assumption~\ref{ass:003}.

\begin{definition}[Degree of $\ups$ and degree of $F$] Let $\ups\in{\mathcal C}_b(\R_+\times{\R_+^{\N\choose 2}\times K^\mathbb{N}})$ depend on $(m,(r_{i,j})_{1\le i<j},(\kappa_i)_{i\in\mathbb{N}})$ only through
$(m,(r_{i,j})_{1\le i<j: i,j \in I},(\kappa_i)_{i\in I})$ for a finite set $I \subseteq \mathbb{N}$. The smallest number $\# I\in\mathbb{N}$ with this property is referred to as the degree of $\ups$.  For a polynomial $F:\,\mathbb{M}^K\to\R$ of finite degree
we denote by  $\mathrm{deg}(F)$ the smallest number for which $F$ is of the form (\ref{s:009}) for a function $\ups\in{\mathcal C}_b(\R_+\times{\R_+^{\N\choose 2}\times K^\mathbb{N}})$
of degree $\mathrm{deg}(h)$.
\label{Def:006}
\end{definition}\sm

\begin{nota}
For $\ups \in \tilde{\Pi}$ with $\ups(\boldsymbol{\cdot},\mr,\vk)=g(m)\phi(\mr)f(\vk)$ and $n:=\mathrm{deg}(\ups)$ write $F^{g,(n,\phi,f)}:=F^\ups$.
\end{nota}\sm

The main result reads as follows:

{
\begin{proposition}[Convergence of generators] For all $F\in\tilde{\Pi}$, under Assumptions~1,~2,~3 and~4
\begin{equation}
\label{e:021}
   \lim_{N\to\infty}\sup_{\smallx\in\mathbb{M}^{K}_N}\big|\Omega_NF(\smallx)-\Omega F(\smallx)\big|=0.
\end{equation}
\label{P:generators}
\end{proposition}\sm

Throughout this section we will fix  a sequence $(a_N)_{N\in\mathbb{N}}$, $a_N\to\infty$ for $N \rightarrow \infty$ but slow enough such that
\begin{equation}
\label{e:a_N}
   \tfrac{a_N}{N}\tNo 0.
\end{equation}\smallskip

To prepare the proof of  Proposition~\ref{P:generators} we first show the following lemma.
\begin{lemma} If $F\in\widetilde{\Pi}$, then
\begin{equation}
\label{e:010y}
   \sup_{\smallx:\,m(\smallx)\le\frac{a_N}{N}}\big|\Omega F(\smallx)\big|\tNo 0.
\end{equation}
\label{L:001}
\end{lemma}

\begin{proof} Recall that $\Omega F(\smallx)=0$ if $m(\smallx)=0$
or if $F$ is constant. Fix therefore $\smallx\in\mathbb{M}^K$ with $0<m(\smallx)\le\frac{\alpha_N}{N}$. As $\Omega$ is linear, we may assume without loss of generality that $g(0)=0$. Then recalling (\ref{y:001}),
\begin{equation}
\label{e:010}
\begin{aligned}
   &\big|\Omega F(\smallx)\big|
   \\
   &\le
   \big|\Omega^{\beta,\Gamma}_{\mbox{\tiny total mass}} F(\smallx)\big|+\big|\Omega^{p,\beta,A}_{\mbox{\tiny trait mutation}}F(\smallx)\big|+\big|\Omega^{p,\beta}_{\mbox{\tiny growth}} F(\smallx)\big|+\big|\Omega^{\Gamma}_{\mbox{\tiny $\Gamma$-reweigh}} F(\smallx)\big|+\big|\Omega^\beta_{\mbox{\tiny natural branching}} F(\smallx)\big|
   \\
   &\le
   \|\phi f\|\cdot\sup_{m\in [0,\frac{a_N}{N}]}\big(\widehat{\Gamma}(m)m|g'(m)|+\overline{\beta}m|g''(m)|+2\mathrm{deg}(\ups)\overline{\beta}|g'(m)|\big)
   \\
   &\;+
   p\|\phi\|\cdot \overline{\beta} \sum\nolimits_{l=1}^{\mathrm{deg}(\ups)}\|A^{(l)}f\|\cdot\sup_{m\in [0,\frac{a_N}{N}]}|g(m)|
%   \\
%   &\;
+
   2p\overline{\beta}\|f\|\cdot \sum\nolimits_{1\le l_1<l_2\le\mathrm{deg}(\ups)} \| \tfrac{\partial\phi}{\partial r_{l_1,l_2}}\|\cdot\sup_{m\in [0,\frac{a_N}{N}]}|g(m)|
   \\
   &\,+2\mathrm{deg}(\ups)\cdot\|\phi f\|\cdot\sup_{m\in [0,\frac{a_N}{N}]}(\tilde{\gamma}(m)+\overline{\gamma}_b)|g(m)|
   \\
   &\;+2\|\phi f\|\overline{\beta}\big(\mathrm{deg}(\ups)+2(\mathrm{deg}(\ups))^2\big)\cdot\sup_{m\in [0,\frac{a_N}{N}]}\tfrac{|g(m)|}{m}.
\end{aligned}
\end{equation}
Thus (\ref{e:010y}) follows as $g$ is continuously differentiable with $g'(0)=0$ and $\tilde{\gamma}$ is continuous.
\end{proof}\smallskip

We continue with the proof of Proposition~\ref{P:generators} and turn
to the approximating operators.
\begin{proof}[Proof of Proposition~\ref{P:generators}]
{Fix $F=F^{g,(n,\phi,f)}\in\tilde{\Pi}$, $N\in\mathbb{N}$, and  $\smallx=\overline{(X,r,\mu)}\in\mathbb{M}^K_N$. We assume once more without loss of generality that $g(0)=0$.

Notice that if $m(\smallx)=0$, then $\Omega_NF(\smallx)=0$ for all $N\in\mathbb{N}$ by definition. Also, if $m(\smallx)=0$, then $\Omega F(\smallx)=0$, cf. below (\ref{y:001}).

If  $\smallx$ is
such that $m(\smallx)=\tfrac{1}{N}$, that is, $\mu=\tfrac{1}{N}\delta_{(x,k)}$ say, then
\begin{equation}
\label{e:009}
\begin{aligned}
   &\big|\Omega_NF(\smallx)\big|
   \\
   &\le
   \big|\Omega^{\mathrm{death}}_NF(\smallx)\big|+\big|\Omega^{\mathrm{birth}}_NF(\smallx)\big|
   \\
   &\le\big(N\cdot\beta(k)+\gamma_{\mathrm{death}}(\tfrac{1}{N},0,k,k)\big)\big|g\big(\tfrac{1}{N}\big)
   \phi\big(\underline{\underline{0}}\big)f\big(\underline{k}\big)\big|
   \\
   &\;\;\;+
    \big(N\cdot\beta(k)+\gamma_{\mathrm{birth}}(\tfrac{1}{N},0,k,k)\big)
    \\
    &\;\hspace{.8cm}
    \cdot\Big|\int\widehat{\alpha}_N(k,\mathrm{d}\tilde{k})\big\{g\big(\tfrac{2}{N}\big)\int\big(\tfrac{1}{2}\delta_{(x,k)}+\tfrac{1}{2}\delta_{(z,\tilde{k})}\big)^{\otimes n}(\mathrm{d}(\underline{x}',\underline{k}'))\,\phi\big(\mr^{(x,z),\frac{1}{N}}(\underline{x}')\big)f\big(\underline{k}'\big)
    -g\big(\tfrac{1}{N}\big)\phi\big(\underline{\underline{0}}\big)f\big(\underline{k}\big)\big\}\Big|
    \\
    &\le
    %\big(\overline{\beta}+\tfrac{1}{N}\overline{\gamma}_b\big)\cdot
    %\big|\int N(\alpha_N(k,\boldsymbol{\cdot})-\delta_k)(\mathrm{d}\tilde{k})\big\{
    %g\big(\tfrac{2}{N}\big)\int\big(\tfrac{1}{2}\delta_{(x,k)}+\tfrac{1}{2}\delta_{(z,\tilde{k})}\big)^{\otimes %n}(\mathrm{d}(\underline{x}',\underline{k}'))\,\phi\big(\mr^{(x,z),\frac{1}{N}}\big)f\big(\underline{k}'\big)
    %\\
    %&\;\hspace{1cm}
    %-g\big(\tfrac{1}{N}\big)\phi\big(\underline{\underline{0}}\big)f\big(\underline{k}\big)\big\}\big|
    %\\
    %&\;
    \big(\overline{\beta}+\tfrac{1}{N}\overline{\gamma}_d\big)\cdot\|\phi f\|\cdot N \big|g\big(\tfrac{1}{N}\big)\big|
    \\
    &\;\;\;+\big(\overline{\beta}+\tfrac{1}{N}\overline{\gamma}_b\big)\cdot\|\phi f\|\cdot N\big| g\big(\tfrac{2}{N}\big)-g\big(\tfrac{1}{N}\big)\big|+2\big(\overline{\beta}+\tfrac{1}{N}\overline{\gamma}_b\big)\cdot\|\phi f\|\cdot N\big|g\big(\tfrac{1}{N}\big)\big|
    \\
    &\tNo 0,
\end{aligned}
\end{equation}
where in the last line we made use once more of the assumption that $g'(0)=0$.\smallskip

Assume next that $\smallx=\overline{(X,r,\mu)}\in\mathbb{M}^K_N$ is such that $ m(\smallx)\ge \tfrac{2}{N}$.
}
%Fix $F=F^{g,(n,\phi,f)}\in\tilde{\Pi}$, $N\in\mathbb{N}$, and $\smallx=\overline{(X,r,\mu)}\in\mathbb{M}^K_N$ with $m(\smallx)\cdot N\ge 2$.
Then (recall (\ref{e:sample}), (\ref{e:OmegaN}) and (\ref{e:OmeggaN}))
\begin{equation}
\label{e:022}
\begin{aligned}
  &\Omega_NF\big(\smallx\big)
  \\
  & =:
   \Omega^{\mathrm{death}}_NF\big(\smallx\big)+\Omega^{\mathrm{birth}}_NF\big(\smallx\big)
   \\
   &=
   Nm(\smallx)\int_{(X\times K)^2}\bar{\mu}^{\otimes 2}(\mathrm{d}((x_1,\kappa_1),(x_2,\kappa_2)))\big\{N\beta(\kappa_2)+\gamma^{\mathrm{death}}(m(\smallx),r(x_1,x_2),\kappa_1,\kappa_2)\big\}
   \\
   &\hspace{.2cm}\cdot
   \Big(g\big(m(\smallx)-\tfrac{1}{N}\big)\int_{(X\times K)^n}\big(\tfrac{Nm(\smallx)\cdot\bar{\mu}-\delta_{(x_2,\kappa_2)}}{Nm(\smallx)-1}\big)^{\otimes n}(\mathrm{d}(\vx,\vk))\,\phi(\mr(\vx))\cdot f(\vk)-g(m(\smallx))\cdot F^{1,(n,\phi,f)}\big(\smallx\big)\Big)
    \\
    &+Nm(\smallx)\int_{(X\times K)^2}\bar{\mu}^{\otimes 2}
   (\mathrm{d}((x_1,\kappa_1),(x_2,\kappa_2)))\big\{N\beta(\kappa_2)+\gamma^{\mathrm{birth}}(m(\smallx),r(x_1,x_2),\kappa_1,\kappa_2)\big\}
   \int_K\widehat{\alpha}_N\big(\kappa_2,\mathrm{d}\tilde{\kappa}_2\big)
     \\
   &\hspace{.2cm}\cdot
   \Big(g\big(m(\smallx)+\tfrac{1}{N}\big)\int_{(X\uplus\{z\}\times K)^n}\big(\tfrac{Nm(\smallx)\cdot\bar{\mu}+\delta_{(z,\tilde{\kappa}_2)}}{Nm(\smallx)+1}\big)^{\otimes n}(\mathrm{d}(\vx,\vk))\,\phi(\underline{\underline{r}}^{(x_2,z),\frac{1}{N}}(\vx))\cdot f(\vk)-g(m(\smallx))\cdot F^{1,(n,\phi,f)}\big(\smallx\big)\Big),
\end{aligned}
\end{equation}
with $r^{(x_2,z),\frac{1}{N}}$ as defined in (\ref{s:006}).
We make use of a Taylor expansion. Namely, for all $g\in{\mathcal C}_b^3(\R_+)$, {$m\ge\tfrac{1}{N}$},
\begin{equation}
\label{exp-g}
   g\big(m\pm\tfrac{1}{N}\big)=g(m)\pm\tfrac{1}{N}\cdot g'(m)+\tfrac{1}{2 N^2}\cdot g''(m)+{\mathcal O}(N^{-3}) C_g(m,N).
\end{equation}
with
\begin{equation}
\label{Cg}
   C_g(m,N)
 :=
   %\begin{cases}
   \max_{\xi \in [m-1/N,m+1/N]}\big|g'''(\xi)\big|.
    %\end{cases}
\end{equation}

 We thereby obtain by Assumptions~\ref{ass:002} and \ref{ass:003},
 \begin{equation}
\label{e:Omegadeath1}
\begin{aligned}
  &\Omega^{\mathrm{death}}_N F^{g,(n,\phi,f)}\big(\smallx\big)
   \\
  &=
    g(m(\smallx))\cdot\Omega^{\mathrm{death}}_N F^{1,(n,\phi,f)}\big(\smallx\big)
   \\
  &\;\;
    +m(\smallx)\big(-g'(m(\smallx))+\tfrac{1}{2N}g''(m(\smallx))\big)\int\nolimits_{(X\times K)^2}\bar{\mu}^{\otimes 2}
   (\mathrm{d}((x_1,\kappa_1),(x_2,\kappa_2)))\,
   \\
   &\hspace{1cm}\cdot\big\{N\beta(\kappa_2)+\gamma^{\mathrm{death}}(m(\smallx),r(x_1,x_2),\kappa_1,\kappa_2)\big\}
     \cdot
   \int_{(X\times K)^n}\big(\tfrac{Nm(\smallx)\cdot\bar{\mu}-\delta_{(x_2,\kappa_2)}}{Nm(\smallx)-1}\big)^{\otimes n}(\mathrm{d}(\vx,\vk))\,\phi(\mr(\vx))\cdot f(\vk)
   \\
   &\;\;
   +m(\smallx)\cdot\big(\overline{\beta}+\frac{1}{N}\tilde{\gamma}(m(\smallx)) \big)\cdot C_g(m(\smallx),N)\cdot {\mathcal O}\big(\tfrac{1}{N}\big).
 \end{aligned}
\end{equation}

Moreover, by Assumptions~\ref{ass:002} and~\ref{ass:007b} (distinguishing between clones and mutants),
\begin{equation}
\label{e:OmegaNbirth1}
\begin{aligned}
  &\Omega^{\mathrm{birth}}_N F^{g,(n,\phi,f)}\big(\smallx\big)
   \\
  &=
    g(m(\smallx))\cdot\Omega^{\mathrm{birth}}_N F^{1,(n,\phi,f)}\big(\smallx\big)
   \\
  &\;\;\;
    +m(\smallx)\big(g'(m(\smallx))+\tfrac{1}{2N}g''(m(\smallx))\big)\int_{(X\times K)^2}\bar{\mu}^{\otimes 2}
   (\mathrm{d}((x_1,\kappa_1),(x_2,\kappa_2)))\big\{N\beta(\kappa_2)+\gamma^{\mathrm{birth}}(m(\smallx),r(x_1,x_2),\kappa_1,\kappa_2)\big\}
     \\
   &\;\;\;\cdot
   \Big((1-p)\int_{(X\times K)^n}\big(\tfrac{Nm(\smallx)\cdot\bar{\mu}+\delta_{(x_2,\kappa_2)}}{Nm(\smallx)+1}\big)^{\otimes n}(\mathrm{d}(\vx,\vk))\,\phi(\mr(\vx))\cdot f(\vk)
   \\
    &\;\;\;+p\int_K\alpha_N\big(\kappa_2,\mathrm{d}\tilde{\kappa}_2\big)
   \int_{(X\uplus\{z\}\times K)^n}\big(\tfrac{Nm(\smallx)\cdot\bar{\mu}+\delta_{(z,\tilde{\kappa}_2)}}{Nm(\smallx)+1}\big)^{\otimes n}(\mathrm{d}(\vx,\vk))\,\phi(\underline{\underline{r}}^{(x_2,z),\frac{1}{N}}(\vx))\cdot f(\vk)\Big)
   \\
   &\;\;\;+m(\smallx)\cdot \big(\overline{\beta}+\frac{1}{N}\gamma_b \big) \cdot C_g(m(\smallx),N)\cdot {\mathcal O}\big(\tfrac{1}{N}\big).
\end{aligned}
\end{equation}

We abbreviate
\begin{equation}
\label{e:abABCDE}
\begin{aligned}
    a&=a_N(m(\smallx))
 :=
    \tfrac{1}{N}g'(m(\smallx))
  \\
     b&=b_N(m(\smallx))
  :=
    \tfrac{1}{2N^2}\cdot  g''(m(\smallx))
    \\
     A&=A_N(m(\smallx),\kappa_2)
  :=
     N^2m(\smallx)\cdot\beta(\kappa_2)
    \\
    B&=B_N(m(\smallx),x_2,\kappa_2)
  :=
     Nm(\smallx)\cdot\int\bar{\mu}(\mathrm{d}(x_1,\kappa_1))\,\gamma^{\mathrm{birth}}\big(m(\smallx), r(x_1,x_2),\kappa_1,\kappa_2\big)
    \\
     C&=C_N(m(\smallx),x_2,\kappa_2)
  :=
     Nm(\smallx)\cdot\int\bar{\mu}(\mathrm{d}(x_1,\kappa_1))\,\gamma^{\mathrm{death}}\big(m(\smallx), r(x_1,x_2),\kappa_1,\kappa_2\big)
    \\
      D&=D_N(m(\smallx),x_2,\kappa_2)
   :=
      (1-p)\int_{(X\times K)^n}\big(\tfrac{Nm(\smallx)\cdot\bar{\mu}+\delta_{(x_2,\kappa_2)}}{Nm(\smallx)+1}\big)^{\otimes n}(\mathrm{d}(\vx,\vk))\,\phi(\mr(\vx))\cdot f(\vk)
   \\
    &\;+p\int_K\alpha_N\big(\kappa_2,\mathrm{d}\tilde{\kappa}_2\big)
   \int_{(X\uplus\{z\}\times K)^n}\big(\tfrac{Nm(\smallx)\cdot\bar{\mu}+\delta_{(z,\tilde{\kappa}_2)}}{Nm(\smallx)+1}\big)^{\otimes n}(\mathrm{d}(\vx,\vk))\,\phi(\underline{\underline{r}}^{(x_2,z),\frac{1}{N}}(\vx))\cdot f(\vk)
   \\
      E&=E_N(m(\smallx),x_2,\kappa_2)
   :=
      \int_{(X\times K)^n}\big(\tfrac{Nm(\smallx)\cdot\bar{\mu}-\delta_{(x_2,\kappa_2)}}{Nm(\smallx)-1}\big)^{\otimes n}
      (\mathrm{d}(\vx,\vk))\,\phi(\mr(\vx))\cdot f(\vk).
\end{aligned}
\end{equation}

Rewrite
\begin{equation}
\label{e:OmegaNdeath1_short}
\begin{aligned}
  &\Omega^{\mathrm{death}}_N F^{g,(n,\phi,f)}\big(\smallx\big)
   \\
  &=
    g(m(\smallx))\cdot \Omega^{\mathrm{death}}_N F^{1,(n,\phi,f)}\big(\smallx\big)
    +\big(-a+b\big)\cdot\int_{X\times K}\bar{\mu}
   (\mathrm{d}(x_2,\kappa_2))\,\big\{A+C \big\}\cdot E
    \\
   &\;+ m(\smallx) \cdot \big(\overline{\beta}+\tfrac{1}{N}\tilde{\gamma}(m(\smallx)) \big) C_g(m(\smallx),N) {\mathcal O}\big(\tfrac{1}{N}\big),
\end{aligned}
\end{equation}
and
\begin{equation}
\label{e:OmegaNbirth1_short}
\begin{aligned}
  &\Omega^{\mathrm{birth}}_N F^{g,(n,\phi,f)}\big(\smallx\big)
   \\
  &=
    g(m(\smallx))\cdot \Omega^{\mathrm{birth}}_N F^{1,(n,\phi,f)}\big(\smallx\big)+\big(a+b\big)\cdot\int_{X\times K}\bar{\mu}(\mathrm{d}(x_2,\kappa_2))\,\big\{A+B \big\}\cdot D
    \\
   &\;+ m(\smallx) \cdot  \big(\overline{\beta}+\frac{1}{N}\gamma_b \big) C_g(m(\smallx),N) {\mathcal O}\big(\tfrac{1}{N}\big).
\end{aligned}
\end{equation}

As
\begin{equation}
\label{e:arithmetic}
\begin{aligned}
   &(a+b)\{A+B\}D-(a-b)\{A+C\}E
   \\
 &=
   aA(D-E)+a(B-C)E{+bA(D+E)+aB(D-E)+b\{ B\cdot D+C\cdot E\}} ,
\end{aligned}
\end{equation}
this yields
\begin{equation}
\label{e:OmegaN1}
\begin{aligned}
  &\Omega_N F^{g,(n,\phi,f)}\big(\smallx\big)
   \\
  &= g(m(\smallx))\cdot \Omega_N F^{1,(n,\phi,f)}\big(\smallx\big)
   \\
  &\; + a\int\mathrm{d}\bar{\mu}\, A (D-E)  + a\int\mathrm{d}\bar{\mu} (B-C) E + b\cdot \int\mathrm{d}\bar{\mu}\, A (D+E)
  \\
  &\;+ a\int\mathrm{d}\bar{\mu}\, B (D-E) +b\int\mathrm{d}\bar{\mu}\,\big\{ B\cdot D+C\cdot E\big\} + m(\smallx)\cdot\big(1\vee \tfrac{\tilde{\gamma}(m(\smallx))}{N}\big) C_g(m(\smallx),N)\cdot {\mathcal O}\big(\tfrac{1}{N}\big)
  \\
  &=: g(m(\smallx))\cdot T_1 + T_2  + E_3.
\end{aligned}
\end{equation}

In what follows we analyse each of the three terms separately.\bigskip

% ----------------------------------------------------------------------

\noindent{\em Step~0 (Preparatory calculations) }
For all $n\ge 2$ fixed, $N\cdot m \geq 2$,
\begin{equation}
\label{e:exp-mun}
\begin{aligned}
  &\big(\tfrac{Nm\bar{\mu}\pm\delta_{(z,\kappa)}}{Nm \pm 1}\big)^{\otimes n}
  \\
  &=
   \big(\bar{\mu}\pm\tfrac{(\delta_{(z,\kappa)}-\bar{\mu})}{Nm \pm 1} \big)^{\otimes n}
   \\
  &=
    \bar{\mu}^{\otimes n} \pm \sum_{l=1}^n \bar{\mu}^{\otimes (l-1)} \otimes \big( \tfrac{\delta_{(z,\kappa)}-\bar{\mu}}{Nm \pm 1} \big) \otimes \bar{\mu}^{\otimes (n-l)}
   \\
  &\;+
    \sum_{1 \leq l_1 < l_2 \leq n} \bar{\mu}^{\otimes (l_1-1)} \otimes \big( \tfrac{\delta_{(z,\kappa)}-\bar{\mu}}{Nm \pm 1} \big) \otimes \bar{\mu}^{\otimes (l_2-l_1-1)} \otimes \big( \tfrac{\delta_{(z,\kappa)}-\bar{\mu}}{Nm \pm 1} \big) \otimes \bar{\mu}^{\otimes (n-l_2)} + {\mathcal O}\big((N m)^{-3}\big),
\end{aligned}
\end{equation}
which we rewrite to
\begin{equation}
\label{e:exp-diff}
\begin{aligned}
   &\big(\tfrac{Nm\bar{\mu}\pm\delta_{(z,\kappa)}}{Nm \pm 1}\big)^{\otimes n}-\bar{\mu}^{\otimes n}
   \\
  &=
    \tfrac{n}{Nm \pm 1}  \Big( \mp 1 + \tfrac{(n-1)}{2 (Nm \pm 1)} \Big) \bar{\mu}^{\otimes n}
   \\
  &\;+
    \tfrac{1}{Nm \pm 1} \big( \pm 1 - \tfrac{n-1}{(Nm \pm 1)} \big)
    \sum_{l=1}^n \bar{\mu}^{\otimes (l-1)} \otimes \delta_{(z,\kappa)} \otimes \bar{\mu}^{\otimes (n-l)}
   \\
  &\;+
    \tfrac{1}{(Nm \pm 1)^2} \sum_{1 \leq l_1 < l_2 \leq n} \bar{\mu}^{\otimes (l_1-1)} \otimes \delta_{(z,\kappa)} \otimes \bar{\mu}^{\otimes (l_2-l_1-1)} \otimes
    \delta_{(z,\kappa)} \otimes \bar{\mu}^{\otimes (n-l_2)} + {\mathcal O}\big((N m)^{-3}\big)
    \\
   &=
   -\tfrac{n(n-1)}{2(Nm\pm 1)^2}\bar{\mu}^{\otimes n}
   +\tfrac{1}{Nm \pm 1} \big( \pm 1 - \tfrac{n-1}{(Nm \pm 1)} \big)
    \sum_{l=1}^n \bar{\mu}^{\otimes (l-1)} \otimes \big(\delta_{(z,\kappa)}-\bar{\mu}\big) \otimes \bar{\mu}^{\otimes (n-l)}
   \\
  &\;+
    \tfrac{1}{(Nm \pm 1)^2} \sum_{1 \leq l_1 < l_2 \leq n} \bar{\mu}^{\otimes (l_1-1)} \otimes \delta_{(z,\kappa)} \otimes \bar{\mu}^{\otimes (l_2-l_1-1)} \otimes
    \delta_{(z,\kappa)} \otimes \bar{\mu}^{\otimes (n-l_2)} + {\mathcal O}\big((N m)^{-3}\big).
\end{aligned}
\end{equation}

This implies that
\begin{equation}
\label{e:exp-diffoh0}
\begin{aligned}
   &\big(\tfrac{Nm\bar{\mu}+\delta_{(z,\kappa)}}{Nm+ 1}\big)^{\otimes n}+\big(\tfrac{Nm\bar{\mu}-\delta_{(z,\kappa)}}{Nm- 1}\big)^{\otimes n}-2\bar{\mu}^{\otimes n}
   \\
   &=
   -\tfrac{n(n-1)}{2}\Big(\tfrac{1}{(Nm+1)^2}+\tfrac{1}{(Nm-1)^2}\Big)\bar{\mu}^{\otimes n}
   \\
   &\;\;\;+\Big(\tfrac{1}{Nm+1} \big(1-\tfrac{n-1}{(Nm+1)}\big)-\tfrac{1}{Nm- 1} \big(1+\tfrac{n-1}{(Nm -1)} \big)\Big)\sum_{l=1}^n \bar{\mu}^{\otimes (l-1)} \otimes \big(\delta_{(z,\kappa)}-\bar{\mu}\big) \otimes \bar{\mu}^{\otimes (n-l)}
   \\
   &\;\;\;+\Big(\tfrac{1}{(Nm+ 1)^2}+\tfrac{1}{(Nm- 1)^2}\Big)\sum_{1 \leq l_1 < l_2 \leq n} \bar{\mu}^{\otimes (l_1-1)} \otimes \delta_{(z,\kappa)} \otimes \bar{\mu}^{\otimes (l_2-l_1-1)} \otimes
    \delta_{(z,\kappa)} \otimes \bar{\mu}^{\otimes (n-l_2)} + {\mathcal O}\big((N m)^{-3}\big)
   \\
   &=
   -\tfrac{n(n-1)((Nm)^2+1)}{((Nm)^2-1)^2}\bar{\mu}^{\otimes n}
    \\
   &\;\;\;+\tfrac{1}{((Nm)^2-1)^2}\Big(-2\big((Nm)^2-1\big)-2(n-1)\big((Nm)^2+1\big)\Big)
   \sum_{l=1}^n \bar{\mu}^{\otimes (l-1)} \otimes \big(\delta_{(z,\kappa)}-\bar{\mu}\big) \otimes \bar{\mu}^{\otimes (n-l)}
   \\
   &\;\;\;+\tfrac{(Nm-1)^2+(Nm+1)^2}{((Nm)^2-1)^2}\sum_{1 \leq l_1 < l_2 \leq n} \bar{\mu}^{\otimes (l_1-1)} \otimes \delta_{(z,\kappa)} \otimes \bar{\mu}^{\otimes (l_2-l_1-1)} \otimes
    \delta_{(z,\kappa)} \otimes \bar{\mu}^{\otimes (n-l_2)} + {\mathcal O}\big((N m)^{-3}\big).
\end{aligned}
\end{equation}

Thus
\begin{equation}
\label{e:exp-diffoh}
\begin{aligned}
   &\big(\tfrac{Nm\bar{\mu}+\delta_{(z,\kappa)}}{Nm+ 1}\big)^{\otimes n}+\big(\tfrac{Nm\bar{\mu}-\delta_{(z,\kappa)}}{Nm- 1}\big)^{\otimes n}-2\bar{\mu}^{\otimes n}
      \\
   &=
   -\tfrac{2n}{(Nm)^2}\sum_{l=1}^n \bar{\mu}^{\otimes (l-1)} \otimes \big(\delta_{(z,\kappa)}-\bar{\mu}\big) \otimes \bar{\mu}^{\otimes (n-l)}
   \\
  &\;\;\;+\tfrac{2}{(Nm)^2}\sum_{1 \leq l_1 < l_2 \leq n} \Big(\bar{\mu}^{\otimes (l_1-1)} \otimes \delta_{(z,\kappa)} \otimes \bar{\mu}^{\otimes (l_2-l_1-1)} \otimes
    \delta_{(z,\kappa)} \otimes \bar{\mu}^{\otimes (n-l_2)}-\bar{\mu}^{\otimes n}\Big) + {\mathcal O}\big((N m)^{-3}\big).
%\\
%  &= \tfrac{n(n+1)(Nm)^2}{((Nm)^2- 1)^2}\bar{\mu}^{\otimes n}
%    -\tfrac{2n(Nm)^2}{((Nm)^2-1)^2}
%    \sum_{l=1}^n \bar{\mu}^{\otimes (l-1)} \otimes \delta_{(z,\kappa)} \otimes \bar{\mu}^{\otimes (n-l)}
%   \\
%  &\;+
%    \tfrac{2(Nm)^2}{((Nm)^2 - 1)^2} \sum_{1 \leq l_1 < l_2 \leq n} \bar{\mu}^{\otimes (l_1-1)} \otimes %\delta_{(z,\kappa)} \otimes \bar{\mu}^{\otimes (l_2-l_1-1)} \otimes
%    \delta_{(z,\kappa)} \otimes \bar{\mu}^{\otimes (n-l_2)} + {\mathcal O}\big((N m)^{-3}\big) + %{\mathcal O}\big(((N m)^2-1)^{-2}\big)
%    \\
%%  &= -\tfrac{2n}{(Nm)^2} \sum_{l=1}^n \bar{\mu}^{\otimes (l-1)} \otimes %\big(\delta_{(z,\kappa)}-\bar{\mu}\big) \otimes \bar{\mu}^{\otimes (n-l)}
%%    \\
%%    &\;+\tfrac{2}{(Nm)^2} \sum_{1 \leq l_1 < l_2 \leq n} \big(\bar{\mu}^{\otimes (l_1-1)} \otimes %\delta_{(z,\kappa)} \otimes \bar{\mu}^{\otimes (l_2-l_1-1)} \otimes
%%    \delta_{(z,\kappa)} \otimes \bar{\mu}^{\otimes (n-l_2)} -\bar{\mu}^{\otimes n}\big)
%%    + {\mathcal O}\big((N m)^{-3}\big)+{\mathcal O}\big(((N m)^2-1)^{-2}\big)+{\mathcal O}\big(\big)
%%    \\
%    &= -\tfrac{2n}{(Nm)^2} \sum_{l=1}^n \bar{\mu}^{\otimes (l-1)} \otimes %\big(\delta_{(z,\kappa)}-\bar{\mu}\big) \otimes \bar{\mu}^{\otimes (n-l)}
%    \\
%    &\;+\tfrac{2}{(Nm)^2} \sum_{1 \leq l_1 < l_2 \leq n} \big(\bar{\mu}^{\otimes (l_1-1)} \otimes %\delta_{(z,\kappa)} \otimes \bar{\mu}^{\otimes (l_2-l_1-1)} \otimes
%    \delta_{(z,\kappa)} \otimes \bar{\mu}^{\otimes (n-l_2)} -\bar{\mu}^{\otimes n}\big)
%    + {\mathcal O}\big((N m)^{-3}\big).
\end{aligned}
\end{equation}
%Recall that had assumed that $Nm\ge 2$. We could therefore identify all above error terms by a term of %order ${\mathcal O}\big((N m)^{-3}\big)$.

Formula (\ref{e:exp-mun}) also yields 
\begin{equation}
\label{e:exp-diffq}
\begin{aligned}
   &\big(\tfrac{Nm\bar{\mu}+\delta_{(z,\kappa)}}{Nm + 1}\big)^{\otimes n}-\big(\tfrac{Nm\bar{\mu}-\delta_{(\tilde{z},\tilde{\kappa})}}{Nm - 1}\big)^{\otimes n}
   \\
    &=
    \sum_{l=1}^n \bar{\mu}^{\otimes (l-1)} \otimes \big( \tfrac{\delta_{(z,\kappa)}-\bar{\mu}}{Nm+ 1} \big) \otimes \bar{\mu}^{\otimes (n-l)}
    +\sum_{l=1}^n \bar{\mu}^{\otimes (l-1)} \otimes \big( \tfrac{\delta_{(\tilde{z},\tilde{\kappa})}-\bar{\mu}}{Nm- 1} \big) \otimes \bar{\mu}^{\otimes (n-l)}
   \\
  &\;+
    \sum_{1 \leq l_1 < l_2 \leq n} \bar{\mu}^{\otimes (l_1-1)} \otimes \big( \tfrac{\delta_{(z,\kappa)}-\bar{\mu}}{Nm+ 1} \big) \otimes \bar{\mu}^{\otimes (l_2-l_1-1)} \otimes \big( \tfrac{\delta_{(z,\kappa)}-\bar{\mu}}{Nm+ 1} \big) \otimes \bar{\mu}^{\otimes (n-l_2)}
    \\
    &\;\;\;-\sum_{1 \leq l_1 < l_2 \leq n} \bar{\mu}^{\otimes (l_1-1)} \otimes \big( \tfrac{\delta_{(\tilde{z},\tilde{\kappa})}-\bar{\mu}}{Nm- 1} \big) \otimes \bar{\mu}^{\otimes (l_2-l_1-1)} \otimes \big( \tfrac{\delta_{(\tilde{z},\tilde{\kappa})}-\bar{\mu}}{Nm- 1} \big) \otimes \bar{\mu}^{\otimes (n-l_2)} + {\mathcal O}\big((N m)^{-3}\big),
    \\
   &=
     \tfrac{1}{Nm}\sum_{l=1}^n \bar{\mu}^{\otimes (l-1)}\otimes\big(\delta_{(z,\kappa)}+\delta_{(\tilde{z},\tilde{\kappa})}-2\bar{\mu}\big)\otimes \bar{\mu}^{\otimes (n-l)}
     +{\mathcal O}\big((N m)^{-2}\big).
\end{aligned}
\end{equation}\bigskip

%.......................................................................

\noindent{\em Step~1 (The term $T_1$) } Recall from (\ref{e:OmegaN1}) the term $T_1=\Omega_N F^{1,(n,\phi,f)}\big(\smallx\big)$.
This term describes the changes we see once we force the total mass to be constant. Recall from (\ref{e:022}) that $\Omega_N=\Omega^{\mathrm{death}}_N+\Omega^{\mathrm{birth}}_N$.
We start with the death part which we split in {\em natural death} which happens on a faster time scale and death due to {\em competition} which occurs more rarely:
\begin{equation}
\label{e:possible4}
\begin{aligned}
    &\Omega^{\mathrm{death}}_N F^{1,(n,\phi,f)}\big(\smallx\big)
   \\
  &=:
    \Omega^{\mathrm{natural~death}}_N F^{1,(n,\phi,f)}\big(\smallx\big)+\Omega^{\mathrm{competition}}_N F^{1,(n,\phi,f)}\big(\smallx\big)
    \\
  &=
    N^2m(\smallx)\int_{(X\times K)^2}\bar{\mu}^{\otimes 2}(\mathrm{d}((x_1,\kappa_1),(x_2,\kappa_2)))\,\beta(\kappa_2)
    \\
    &\hspace{1cm}\cdot\int_{(X\times K)^n}\Big(\big(\tfrac{Nm(\smallx)\cdot\bar{\mu}-\delta_{(x_2,\kappa_2)}}{Nm(\smallx)-1}\big)^{\otimes n}- \bar{\mu}^{\otimes n}\Big)
   (\mathrm{d}(\vx,\vk))\,\phi(\mr(\vx))\cdot f(\vk)
    \\
  &\;+
    Nm(\smallx)\int_{(X\times K)^2}\bar{\mu}^{\otimes 2}(\mathrm{d}((x_1,\kappa_1),(x_2,\kappa_2)))\,\gamma^{\mathrm{death}}(m(\smallx),r(x_1,x_2),
    \kappa_1,\kappa_2)
    \\
   &\;
    \hspace{1cm}\cdot\int_{(X\times K)^n}\Big(\big(\tfrac{Nm(\smallx)\cdot\bar{\mu}-\delta_{(x_2,\kappa_2)}}{Nm(\smallx)-1}\big)^{\otimes n}- \bar{\mu}^{\otimes n}\Big)
   (\mathrm{d}(\vx,\vk))\,\phi(\mr(\vx))\cdot f(\vk).
\end{aligned}
\end{equation}

For the birth-part, we have to take into consideration possible {\em mutation} events. Thus
\begin{equation}
\label{e:possible1}
\begin{aligned}
    &\Omega^{\mathrm{birth}}_N F^{1,(n,\phi,f)}\big(\smallx\big)
   \\
  &=:
    (1-p)\cdot\Omega^{\mathrm{birth; no~mutation}}_N F^{1,(n,\phi,f)}\big(\smallx\big)+p\cdot\Omega^{\mathrm{birth; mutation}}_N F^{1,(n,\phi,f)}\big(\smallx\big)
    \\
  &=
    Nm(\smallx)\int_{(X\times K)^2}\bar{\mu}^{\otimes 2}(\mathrm{d}((x_1,\kappa_1),(x_2,\kappa_2)))\,\big\{N\beta(\kappa_2)+\gamma^{\mathrm{birth}}(m(\smallx),r(x_1,x_2),\kappa_1,\kappa_2)\big\}
    \\
   &\;
    \cdot\Big\{(1-p) \cdot\int_{(X\times K)^n}\Big(\big(\tfrac{Nm(\smallx)\cdot\bar{\mu}+\delta_{(x_2,\kappa_2)}}{Nm(\smallx)+1}\big)^{\otimes n}- \bar{\mu}^{\otimes n}\Big)
   (\mathrm{d}(\vx,\vk))\,\phi(\mr(\vx))\cdot f(\vk)
    \\
  & \quad  +p\int_K\alpha_N\big(\kappa_2,\mathrm{d}\tilde{\kappa}_2\big)
   \hspace{-.3cm}\int\limits_{(X\uplus\{z\}\times K)^n}\hspace{-.3cm}\Big(\big(\tfrac{Nm(\smallx)\cdot\bar{\mu}+\delta_{(z,\tilde{\kappa}_2)}}{Nm(\smallx)+1}\big)^{\otimes n} -\bar{\mu}^{\otimes n}\Big)(\mathrm{d}(\vx,\vk))\,\phi\big(\underline{\underline{r}}^{(x_2,z),\frac{1}{N}}(\vx)\big)\cdot f(\vk)\Big)\Big\}.
\end{aligned}
\end{equation}

We also distinguish for each of the birth events (with and without mutation) between {\em natural birth events} which happen at a fast time scale and {\em enhancement}  which occurs
more rarely. That is, we use the notation (with the terms obtained similar as in (\ref{e:possible4}))
\begin{equation}
\label{e:possible2}
\begin{aligned}
    \Omega^{\mathrm{birth}}_N
  &=:
    (1-p)\cdot\Omega^{\mathrm{natural~birth;~no~mutation}}_N+(1-p)\cdot\Omega^{\mathrm{enhancement;~no~mutation}}_N
     \\
   &\;\;\;+
    p\cdot\Omega^{\mathrm{natural~birth;~mutation}}_N+p\cdot\Omega^{\mathrm{enhancement;~mutation}}_N.
 \end{aligned}
\end{equation}

We therefore base our study of $T_1$ on the following decomposition:
\begin{equation}
\begin{aligned}
\label{e:024}
   T_1
  &=
     \big(\Omega^{\mathrm{enhancement;~no~mutation}}_N+\Omega^{\mathrm{competition}}_N\big)F^{1,(n,\phi,f)}(\smallx)
    \\
     &\;\;\;+\big(\Omega^{\mathrm{natural~birth;~no~mutation}}_N+\Omega^{\mathrm{natural~death}}_N\big)F^{1,(n,\phi,f)}(\smallx)
     \\
    &\;\;\;
      +p\cdot\big(\Omega^{\mathrm{natural~birth;~mutation}}_N-\Omega^{\mathrm{natural~birth;~no~mutation}}_N\big)F^{1,(n,\phi,f)}(\smallx)
           \\
    &\;\;\; +p\cdot\big(\Omega^{\mathrm{enhancement;~mutation}}_N-\Omega^{\mathrm{enhancement;~no~mutation}}_N\big)F^{1,(n,\phi,f)}(\smallx)
       \\
   &=:
      T_{1.1}+T_{1.2}+T_{1.3}+E_{1.4}.
\end{aligned}
\end{equation}
Once more we are handling the different terms separately.
Recall the sequence $(a_N)$ from (\ref{e:a_N}). \smallskip

\noindent{\em Step~1.1 (Reweighing the sampling measure with respect to $\Gamma$) } Recall $\Gamma$ from (\ref{s:010}), and $\Omega^{\Gamma}_{\mbox{\tiny $\Gamma$-reweigh}}$ from (\ref{y:006}).
By (\ref{e:possible4}) and (\ref{e:possible1}),
\begin{equation}
\label{e:025x}
\begin{aligned}
   T_{1.1}
   &=
   Nm(\smallx)\int_{(X\times K)^2}\bar{\mu}^{\otimes 2}(\mathrm{d}((x_1,\kappa_1),(x_2,\kappa_2)))\,\gamma^{\mathrm{death}}(m(\smallx),r(x_1,x_2),
    \kappa_1,\kappa_2)
    \\
   &\;
    \hspace{2cm}\cdot\int_{(X\times K)^n}\Big(\big(\tfrac{Nm(\smallx)\cdot\bar{\mu}-\delta_{(x_2,\kappa_2)}}{Nm(\smallx)-1}\big)^{\otimes n}- \bar{\mu}^{\otimes n}\Big)
   (\mathrm{d}(\vx,\vk))\,\phi(\mr(\vx))\cdot f(\vk)
   \\
   &\;\;\;+ Nm(\smallx)\int_{(X\times K)^2}\bar{\mu}^{\otimes 2}(\mathrm{d}((x_1,\kappa_1),(x_2,\kappa_2)))\,
\gamma^{\mathrm{birth}}(m(\smallx),r(x_1,x_2),\kappa_1,\kappa_2)
    \\
   &\;
    \hspace{2cm}\cdot\int_{(X\times K)^n}\Big(\big(\tfrac{Nm(\smallx)\cdot\bar{\mu}+\delta_{(x_2,\kappa_2)}}{Nm(\smallx)+1}\big)^{\otimes n}- \bar{\mu}^{\otimes n}\Big)
   (\mathrm{d}(\vx,\vk))\,\phi(\mr(\vx))\cdot f(\vk)
   \\
   &=Nm(\smallx)\int_{(X\times K)^2}\bar{\mu}^{\otimes 2}(\mathrm{d}((x_1,\kappa_1),(x_2,\kappa_2)))\,\gamma^{\mathrm{death}}(m(\smallx),r(x_1,x_2),
    \kappa_1,\kappa_2)
    \\
   &\;
    \hspace{1.5cm}\cdot\int_{(X\times K)^n}\Big(\big(\tfrac{Nm(\smallx)\cdot\bar{\mu}-\delta_{(x_2,\kappa_2)}}{Nm(\smallx)-1}\big)^{\otimes n}+\big(\tfrac{Nm(\smallx)\cdot\bar{\mu}+\delta_{(x_2,\kappa_2)}}{Nm(\smallx)+1}\big)^{\otimes n} - 2\bar{\mu}^{\otimes n}\Big)
   (\mathrm{d}(\vx,\vk))\,\phi(\mr(\vx))\cdot f(\vk)
   \\
   &\;\;\;+
Nm(\smallx)\int_{(X\times K)^2}\bar{\mu}^{\otimes 2}(\mathrm{d}((x_1,\kappa_1),(x_2,\kappa_2)))\,
\Gamma(m(\smallx),r(x_1,x_2),\kappa_1,\kappa_2)
    \\
   &\;
    \hspace{1.5cm}\cdot\int_{(X\times K)^n}\Big(\big(\tfrac{Nm(\smallx)\cdot\bar{\mu}+\delta_{(x_2,\kappa_2)}}{Nm(\smallx)+1}\big)^{\otimes n}- \bar{\mu}^{\otimes n}\Big)
   (\mathrm{d}(\vx,\vk))\,\phi(\mr(\vx))\cdot f(\vk).
 \end{aligned}
\end{equation}

Thus by (\ref{e:exp-diffoh}) and (\ref{e:exp-diff}) together with Assumption~\ref{ass:003},
 \begin{equation}
\label{e:025}
\begin{aligned}
   T_{1.1}
   &=\Big(1\vee \tilde{\gamma}\big(m(\smallx)\big)\Big)\cdot{\mathcal O}\big((Nm(\smallx))^{-1}\big)
   \\
   &\;\;+
     \int_{(X\times K)^2}\bar{\mu}^{\otimes 2}(\mathrm{d}((x_1,\kappa_1),(x_2,\kappa_2)))\,
     \Gamma(m(\smallx),r(x_1,x_2),\kappa_1,\kappa_2)
     \\
     &\hspace{1cm}\cdot\sum_{l=1}^n \int_{(X\times K)^n} \bar{\mu}^{\otimes (l-1)}\otimes\big(\delta_{(x_2,\kappa_2)}-\bar{\mu}\big)\otimes\bar{\mu}^{\otimes (n-l)}(\mathrm{d}(\vx,\vk))\,\phi(\mr(\vx))\cdot f(\vk)
   \\
    &=
      \Omega^{\Gamma}_{\mbox{\tiny $\Gamma$-reweigh}}F^{1,(n,\phi,f)}(\smallx) +\Big(1\vee \tilde{\gamma}\big(m(\smallx)\big)\Big)\cdot{\mathcal O}\big((Nm(\smallx))^{-1}\big).
\end{aligned}
\end{equation}

Hence, for some $C\in(0,\infty)$,
\begin{equation}
\label{e:leT11a}
\begin{aligned}
   &\sup_{\smallx\in\mathbb{M}^K_N;\,\frac{2}{N}\le m(\smallx)\le\frac{a_N}{N}}\big|g(m(\smallx))\cdot T_{1.1}-\Omega^{\Gamma}_{\mbox{\tiny $\Gamma$-reweigh}}F^{g,(n,\phi,f)}(\smallx)\big|
   \\
   &\le
    C\cdot\sup_{m\in[\frac{2}{N},\frac{a_N}{N}]}\tfrac{|g(m)|}{m}\tfrac{(1\vee \tilde{\gamma}(m))}{N}\tNo 0,
\end{aligned}
\end{equation}
as $g$ is bounded differentiable at $m=0$ and $\tilde{g}$ is continuous. Moreover,
\begin{equation}
\label{e:geT11b}
\begin{aligned}
   &\sup_{\smallx\in\mathbb{M}^K_N;\,m(\smallx)\ge\frac{a_N}{N}}\big|g(m(\smallx))\cdot T_{1.1}-\Omega^{\Gamma}_{\mbox{\tiny $\Gamma$-reweigh}}F^{g,(n,\phi,f)}(\smallx)\big|
   \\
   &\le
    C\cdot a^{-1}_N\cdot\sup_{m\in[\frac{a_N}{N},\infty)}|g(m)|(1\vee \tilde{\gamma}(m))\tNo 0,
\end{aligned}
\end{equation}
as $\limsup_{m\to\infty}(1 \vee \tilde{\gamma}(m))|g(m)|<\infty$ (compare (\ref{e:tildePi})) and $a_N \rightarrow \infty$ for $N \rightarrow \infty$.
\smallskip

\noindent{\em Step~1.2 (Effect of natural branching on the genealogy) }
Recall $\widehat{\beta}^{\smallx}$ from (\ref{hat:beta}), and the replacement map
$\Theta_{l_1,l_2}$ from (\ref{e:018_2}) and (\ref{e:018_3}). Recall further $\Omega^{\beta}_{\mbox{\tiny $\beta$-natural branching}}$, $\Omega^{\beta}_{\mbox{\tiny $\beta$-reweigh}}$ and $\Omega^{\beta}_{\mbox{\tiny resample}}$
from (\ref{y:007}).
By (\ref{e:possible4}) and (\ref{e:possible1}),
\begin{equation}
\label{e:026pre}
\begin{aligned}
   &T_{1.2}
   \\
   &=
   N^2m(\smallx)\int_{X\times K}\bar{\mu}(\mathrm{d}(x',\kappa'))\,\beta(\kappa')
   \\
   &\hspace{1cm}\cdot
   \int_{(X\times K)^n}\Big(\big(\tfrac{Nm(\smallx)\cdot\bar{\mu}-\delta_{(x',\kappa')}}{Nm(\smallx)-1}\big)^{\otimes n}+\big(\tfrac{Nm(\smallx)\cdot\bar{\mu}+\delta_{(x',\kappa')}}{Nm(\smallx)+1}\big)^{\otimes n} - 2\bar{\mu}^{\otimes n}\Big)
   (\mathrm{d}(\vx,\vk))\,\phi(\mr(\vx))\cdot f(\vk).
\end{aligned}
\end{equation}

Using once more (\ref{e:exp-diffoh}), we find
that
\begin{equation}
\label{e:026}
\begin{aligned}
   T_{1.2}
    &=-\tfrac{2n}{m(\smallx)}
    \int\sum_{l=1}^n\bar{\mu}^{\otimes n}(\mathrm{d}(\vx,\vk))\,\big(\beta(\kappa_l)-\widehat{\beta}^{\smallx}\big)\cdot\phi(\mr(\vx))\cdot f(\vk)
    \\
    &\;\;
    +\tfrac{1}{m(\smallx)}\int\bar{\mu}^{\otimes n}(\mathrm{d}(\vx,\vk))\,\sum_{1\le l_1\not=l_2\le n}\big\{\beta(\kappa_{l_1})\Theta_{l_1,l_2}´(\phi\cdot f)-\widehat{\beta}^{\smallx}\phi\cdot f\big\}(\mr,\vk)\big\}
    + N^2 m(\smallx)\cdot {\mathcal O}\big((Nm(\smallx))^{-3}\big)
    \\
    &=
    -\tfrac{1}{m(\smallx)}
    \int\sum_{l=1}^n\bar{\mu}^{\otimes n}(\mathrm{d}(\vx,\vk))\,\big(\beta(\kappa_l)-\widehat{\beta}^{\smallx}\big)\cdot\phi(\mr(\vx))\cdot f(\vk)+{\mathcal O}\big(N^{-1}(m(\smallx))^{-2}\big)
   \\
   &\;\;+
    \tfrac{1}{m(\smallx)}\int\bar{\mu}^{\otimes n}(\mathrm{d}(\vx,\vk))\,\sum_{1\le l_1,l_2\le n}\Big\{\beta(\kappa_{l_1})\cdot\big(\Theta_{l_1,l_2}(\phi\cdot f)-\phi\cdot f\big)+\big(\widehat{\beta}^{\smallx}-\beta(\kappa_{l_2})\big) \phi\cdot f \Big\}(m(\smallx),\mr,\vk) \\
    &=
     \big(\Omega^{\beta}_{\mbox{\tiny $\beta$-reweigh}}+\Omega^{\beta}_{\mbox{\tiny resample}}\big)F^{1,(n,\phi,f)}(\smallx)+{\mathcal O}\big(N^{-1}(m(\smallx))^{-2}\big)
     \\
     &=\Omega^{\beta}_{\mbox{\tiny $\beta$-natural branching}}F^{1,(n,\phi,f)}(\smallx)+{\mathcal O}\big(N^{-1}(m(\smallx))^{-2}\big).
\end{aligned}
\end{equation}

Thus for some $C\in(0,\infty)$,
\begin{equation}
\label{e:leT12}
\begin{aligned}
   &\sup_{\smallx\in\mathbb{M}^K_N;\,\frac{2}{N}\le m(\smallx)\le\frac{a_N}{N}}\big|g(m(\smallx))\cdot T_{1.2}-\Omega^{\beta}_{\mbox{\tiny natural branching}}F^{g,(n,\phi,f)}(\smallx)\big|
   \\
   &\le
    C\cdot\sup_{m\in[\frac{2}{N},\frac{a_N}{N}]}\tfrac{|g(m)|}{m} \tfrac{1}{N \cdot m}\tNo 0,
\end{aligned}
\end{equation}
as $g$ is bounded differentiable at $m=0$ and $g'(0)=0$. Moreover,
\begin{equation}
\label{e:geT12}
\begin{aligned}
   &\sup_{\smallx\in\mathbb{M}^K_N;\,m(\smallx)\ge\frac{a_N}{N}}\big|g(m(\smallx))\cdot T_{1.2}-\Omega^{\beta}_{\mbox{\tiny natural branching}}F^{g,(n,\phi,f)}(\smallx)\big|
   \\
   &\le
    C\cdot a^{-1}_N\cdot\sup_{m\in[\frac{a_N}{N},\infty)}\tfrac{|g(m)|}{m}\tNo 0,
\end{aligned}
\end{equation}
as $\sup_{m>0}\tfrac{g(m)}{m}<\infty$ and $a_N\tNo\infty$.
\smallskip

\noindent{\em Step~1.3 (Trait mutation and substitution distance growth) } Recall $r^{(x_2,z),l}$ from (\ref{s:006}) as well as $\Omega^{p,\beta,A}_{\mbox{\tiny trait~mutation}}$ and $\Omega^{p,\beta}_{\mbox{\tiny growth}}$ from (\ref{y:004}) and (\ref{y:005}). Then by (\ref{e:possible1}),
\begin{equation}
\label{e:027}
\begin{aligned}
 &T_{1.3}
 \\
 &=
 p\cdot N^2m(\smallx)\cdot\int_{X\times K}\bar{\mu}(\mathrm{d}(x_2,\kappa_2))\,\beta(\kappa_2)
 \\
   &\;\hspace{.8cm}
     \cdot\Big\{\int_K
    \alpha_N\big(\kappa_2,\mathrm{d}\tilde{\kappa}_2\big)\int_{(X\uplus\{z\}\times K)^n}\Big(\big(\tfrac{Nm(\smallx)\cdot\bar{\mu}+\delta_{(z,\tilde{\kappa}_2)}}{Nm(\smallx)+1}\big)^{\otimes n}-\bar{\mu}^{\otimes n}\Big)
     (\mathrm{d}(\vx,\vk))\,\phi(\mr^{(x_2,z),\frac{1}{N}}(\vx))\cdot f\big(\vk\big)
     \\
     &\;\;\;\hspace{.8cm}-\int_{(X\uplus\{z\}\times K)^n}\Big(\big(\tfrac{Nm(\smallx)\cdot\bar{\mu}+\delta_{(x_2,\kappa_2)}}{Nm(\smallx)+1}\big)^{\otimes n}-\bar{\mu}^{\otimes n}\Big)
     (\mathrm{d}(\vx,\vk))\,\phi(\mr^{(x_2,z),\tfrac{1}{N}}(\vx))\cdot f\big(\vk\big)\Big\}
     \\
  &=
   p\cdot\int_{X\times K}\bar{\mu}(\mathrm{d}(x_2,\kappa_2))\,\beta(\kappa_2)\int_K
    N\Big(\alpha_N\big(\kappa_2,\mathrm{d}\tilde{\kappa}_2\big)
    -\delta\big(\kappa_2,\mathrm{d}\tilde{\kappa}_2\big)\Big)
    \\
   &\;\hspace{0.3cm}
     \cdot\int_{(X\uplus\{z\}\times K)^n} Nm(\smallx)\Big(\big(\tfrac{Nm(\smallx)\cdot\bar{\mu}+\delta_{( x_2,\tilde{\kappa}_2)}}{Nm(\smallx)+1}\big)^{\otimes n}-\bar{\mu}^{\otimes n}\Big)(\mathrm{d}(\vx,\vk))\,\phi(\mr^{(x_2,z),\tfrac{1}{N}}(\vx))\cdot f\big(\vk\big)
   \\
  &\;\;\;+
     p\cdot\int_{X\times K}\bar{\mu}(\mathrm{d}(x_2,\kappa_2))\,\beta(\kappa_2)\int_K
     \alpha_N\big(\kappa_2,\mathrm{d}\tilde{\kappa}_2\big)
    \\
   &\;\hspace{0.3cm}
     \cdot\int_{(X\uplus\{z\}\times K)^n} Nm(\smallx)\cdot\Big(\big(\tfrac{Nm(\smallx)\cdot\bar{\mu}+\delta_{(z,\tilde{\kappa}_2)}}{Nm(\smallx)+1}\big)^{\otimes n}-\big(\tfrac{Nm(\smallx)\cdot\bar{\mu}+\delta_{(x_2,\tilde{\kappa}_2)}}{Nm(\smallx)+1}\big)^{\otimes n}\Big)
   (\mathrm{d}(\vx,\vk))\,N\phi(\mr^{(x_2,z),\tfrac{1}{N}}(\vx))\cdot f\big(\vk\big).
\iffalse{   \\
  &= p\cdot Nm(\smallx)\cdot\int_{X\times K}\bar{\mu}(\mathrm{d}(x_2,\kappa_2))\,\beta(\kappa_2)\int_K
    N\Big(\alpha_N\big(\kappa_2,\mathrm{d}\tilde{\kappa}_2\big)-\delta\big(\kappa_2,\mathrm{d}\tilde{\kappa}_2\big)\Big)
    \\
   &\;\hspace{2cm}
     \cdot\int_{(X\uplus\{z\}\times K)^n}\Big(\big(\tfrac{Nm(\smallx)\cdot\bar{\mu}+\delta_{(z,\tilde{\kappa}_2)}}{Nm(\smallx)+1}\big)^{\otimes n}-\big(\tfrac{Nm(\smallx)\cdot\bar{\mu}+\delta_{(x_2,\kappa_2)}}{Nm(\smallx)+1}\big)^{\otimes n}\Big)(\mathrm{d}(\vx,\vk))\,\phi(\mr^{(x_2,z),0}(\vx))\cdot f\big(\vk\big)
  \\
  &\;\;\;+
     p\cdot Nm(\smallx)\cdot\int_{X\times K}\bar{\mu}(\mathrm{d}(x_2,\kappa_2))\,\beta(\kappa_2)\int_K
     N\Big(\alpha_N\big(\kappa_2,\mathrm{d}\tilde{\kappa}_2\big)-\delta\big(\kappa_2,\mathrm{d}\tilde{\kappa}_2\big)\Big)
    \\
   &\;\hspace{2cm}
     \cdot\int_{(X\uplus\{z\}\times K)^n}\big(\tfrac{Nm(\smallx)\cdot\bar{\mu}+\delta_{(z,\tilde{\kappa}_2)}}{Nm(\smallx)+1}\big)^{\otimes n}
   (\mathrm{d}(\vx,\vk))\,\Big\{\phi(\mr^{(x_2,z),\tfrac{1}{N}}(\vx))-\phi(\mr^{(x_2,z),0}(\vx))\Big\}\cdot f\big(\vk\big).}\fi
 %  &{\color{red}\mbox{hier fehlen noch Terme}}
\end{aligned}
\end{equation}

Applying
(\ref{e:exp-diff}) we find that
\begin{equation}
\label{e:027b}
\begin{aligned}
 &T_{1.3}
 \\
 &=
 p\int_{X\times K}\bar{\mu}(\mathrm{d}(x_2,\kappa_2))\,\beta(\kappa_2)\int_K
    N\Big(\alpha_N\big(\kappa_2,\mathrm{d}\tilde{\kappa}_2\big)
    -\delta\big(\kappa_2,\mathrm{d}\tilde{\kappa}_2\big)\Big)
    \\
   &\;\hspace{1cm}
     \cdot\Big\{\int_{(X\uplus\{z\}\times K)^n}\sum_{l=1}^n \bar{\mu}^{\otimes (l-1)}\otimes\big(\delta_{(x_2,\tilde{\kappa}_2)}-\bar{\mu}\big)\otimes \bar{\mu}^{\otimes (n-l)}(\mathrm{d}(\vx,\vk))\,\phi(\mr^{(x_2,z),\tfrac{1}{N}}(\vx)) \cdot f\big(\vk\big)
     +{\mathcal O}\big((Nm(\smallx))^{-1}\big)\Big\}
     \\
     &\;\;\;+p\cdot\int_{X\times K}\bar{\mu}(\mathrm{d}(x_2,\kappa_2))\,\beta(\kappa_2)\int_K\alpha_N\big(\kappa_2,
     \mathrm{d}\tilde{\kappa}_2\big)
    \\
   &\;\hspace{1cm}
     \cdot\Big\{\int_{(X\uplus\{z\}\times K)^n}\sum_{l=1}^n \bar{\mu}^{\otimes (l-1)} \otimes \big(\delta_{(z,\tilde{\kappa}_2)}-\delta_{(x_2,\tilde{\kappa}_2} \big) \otimes \bar{\mu}^{\otimes (n-l)}
   (\mathrm{d}(\vx,\vk))\,\tfrac{\phi(\mr^{(x_2,z),\tfrac{1}{N}}(\vx)) }{\frac{1}{N}}f\big(\vk\big)
   \\
   &\hspace{3cm}+{\mathcal O}\big((Nm(\smallx))^{-1}\big)\Big\}.
\end{aligned}
\end{equation}

Recall the mutation operators $A^{(l)};\,l\in\{1,...,n\}$ from (\ref{y:004}) which act on a function as the mutation operator $A$ on the function of the $l^{\mathrm{th}}$-trait-coordinate (assuming that all other variables are kept constant).
Then by Assumption~\ref{ass:004} on the rescaling of the mutation kernel, and as $\phi\in{\mathcal C}_b^1(\R_+^{n\choose 2})$,
\begin{equation}
\label{e:027c}
\begin{aligned}
  &T_{1.3}
  \\
  &=p\int_{X\times K}\bar{\mu}(\mathrm{d}(x_2,\kappa_2))\,\beta(\kappa_2)\int_K
    N\Big(\alpha_N\big(\kappa_2,\mathrm{d}\tilde{\kappa}_2\big)
    -\delta\big(\kappa_2,\mathrm{d}\tilde{\kappa}_2\big)\Big)
    \\
   &\;\hspace{1cm}
     \cdot\Big\{\int_{(X\uplus\{z\}\times K)^n}\sum_{l=1}^n \bar{\mu}^{\otimes (l-1)}\otimes\delta_{(x_2,\tilde{\kappa}_2)}\otimes \bar{\mu}^{\otimes (n-l)}(\mathrm{d}(\vx,\vk))\,\phi(\mr^{(x_2,z),\tfrac{1}{N}}(\vx)) \cdot f\big(\vk\big)
     +{\mathcal O}\big((Nm(\smallx))^{-1}\big)\Big\}
     \\
     &\;\;\;+p\cdot\int_{X\times K}\bar{\mu}(\mathrm{d}(x_2,\kappa_2))\,\beta(\kappa_2)\int_K\alpha_N\big(\kappa_2,
     \mathrm{d}\tilde{\kappa}_2\big)
    \\
   &\;\hspace{1cm}
     \cdot\Big\{\int_{(X\uplus\{z\}\times K)^n}\sum_{l=1}^n \bar{\mu}^{\otimes (l-1)} \otimes \big( \delta_{(z,\tilde{\kappa}_2)} - \delta_{(x_2,\tilde{\kappa}_2)} \big) \otimes \bar{\mu}^{\otimes (n-l)}
   (\mathrm{d}(\vx,\vk))\,\tfrac{\phi(\mr^{(x_2,z),\tfrac{1}{N}}(\vx))}{\frac{1}{N}}f\big(\vk\big)
   +{\mathcal O}\big((Nm(\smallx))^{-1}\big)\Big\}
  \\
  &=
    p\cdot\int_{(X\times K)^n}\bar{\mu}^{\otimes n}(\mathrm{d}(\vx,\vk))\,
     \phi(\mr(\vx))\cdot \sum_{l=1}^n\beta(\kappa_l)\cdot A^{(l)}f\big(\vk\big)+o_N(1)+{\mathcal O}\big((Nm(\smallx))^{-1}\big)
          \\
     &\;\;\;+
      p\cdot\int_{(X\times K)^n}\bar{\mu}^{\otimes n}(\mathrm{d}(\vx,\vk))\,\sum_{1\le l_1<l_2\le n}\Big\{\beta(\kappa_{l_1})+\beta(\kappa_{l_2})\Big\}\cdot\tfrac{\partial \phi}{\partial r_{l_1,l_2}}(\mr(\vx)\big)\cdot f\big(\vk\big)
     \\
    &= \Omega^{p,\beta,A}_{\mbox{\tiny trait~mutation}}F^{1,(n,\phi,f)}(\smallx)+\Omega^{p,\beta}_{\mbox{\tiny growth}}F^{1,(n,\phi,f)}(\smallx)+o_N(1)+{\mathcal O}\big( \tfrac{1}{Nm(\smallx)}\big).
 \end{aligned}
\end{equation}

{
Thus for some $C\in(0,\infty)$,
\begin{equation}
\label{e:leT13}
\begin{aligned}
   &\sup_{\smallx\in\mathbb{M}^K_N;\,\frac{2}{N}\le m(\smallx)\le\frac{a_N}{N}}\big|g(m(\smallx))\cdot T_{1.3}-(\Omega^{p,\beta,A}_{\mbox{\tiny trait mutation}}+\Omega^{p,\beta}_{\mbox{\tiny growth}})F^{g,(n,\phi,f)}(\smallx)\big|
   \\
   &\le
    C\cdot\sup_{m\in[\frac{2}{N},\frac{a_N}{N}]}|g(m)|\tNo 0,
\end{aligned}
\end{equation}
as $g(0)=0$. Moreover, for some sequence $c_N \downarrow 0$ for $N \rightarrow \infty$,
\begin{equation}
\label{e:geT13}
\begin{aligned}
   &\sup_{\smallx\in\mathbb{M}^K_N;\,m(\smallx)\ge\frac{a_N}{N}}\big|g(m(\smallx))\cdot T_{1.3}-(\Omega^{p,\beta,A}_{\mbox{\tiny trait mutation}}+\Omega^{p,\beta}_{\mbox{\tiny growth}})F^{g,(n,\phi,f)}(\smallx)\big|
   \\
   &\le
    C\cdot (c_N+a^{-1}_N)\cdot\sup_{m\in[\frac{a_N}{N},\infty)}|g(m)|\tNo 0,
\end{aligned}
\end{equation}
as $g$ is bounded and $a_N\tNo\infty$.
}
\smallskip

\noindent{\em Step~1.4 (The term $E_{1.4}$) } It turns out that this last term is negligible in the limit as $N\to\infty$. Indeed, analogous calculations to Step~1.3 show that
$E_{1,4}$ goes to zero as $N\to\infty$ uniformly over all $m\ge\frac{2}{N}$.
%\begin{equation}
%\label{e:027b}
%\begin{aligned}
% N\cdot E_{1.4}
%  &=
%    p\cdot\int_{(X\times K)^n}\bar{\mu}^{\otimes n}(\mathrm{d}(\vx,\vk))\,
%     \phi(\mr^{z,0}(\vx))\cdot \sum_{l=1}^n\int\bar{\mu}(\mathrm{d}(x_1,\kappa_1))\gamma^{\mathrm{birth}}(m(\smallx),r(x_1,x_l),\kappa_1,\kappa_l)\cdot A^{(l)}f\big(\vk\big)
%           \\
%     &\;\;\;+
%      p\cdot\int_{(X\times K)^n}\bar{\mu}^{\otimes n}(\mathrm{d}(\vx,\vk))\,\sum_{1\le l_1<l_2\le n}\int\bar{\mu}(\mathrm{d}(x_1,\kappa_1))\Big\{\gamma^{\mathrm{birth}}(m(\smallx),r(x_1,x_{l_1}),\kappa_1,\kappa_{l_1})+
%      \\
%     &\hspace{3cm}
%      \gamma^{\mathrm{birth}}(m(\smallx),r(x_1,x_{l_2}),\kappa_1,\kappa_{l_2})\Big\}
%      \cdot\tfrac{\partial \phi}{\partial r_{l_1,l_2}}(\mr(\vx)\big)\cdot f\big(\vk\big)
%     \\
%     &\;\;\;+o_N(1),
% \end{aligned}
%\end{equation}
%which implies that $E_{1.4}={\mathcal O}()$.
\bigskip

%.......................................................................

\noindent{\em Step~2 (The term $T_2$) }
Recall $a$, $b$, $A$, $B$, $C$, $D$ and $E$ from (\ref{e:abABCDE}),
and $T_2=a\int\mathrm{d}\bar{\mu}\, A (D-E)  + a\int\mathrm{d}\bar{\mu} (B-C) E +
b\cdot \int\mathrm{d}\bar{\mu}\, A (D+E)$ from (\ref{e:OmegaN1}). This term describes the evolution of the total mass. Recall $\Gamma(m(\smallx),r_{1,2},\kappa_1,\kappa_2)$ from (\ref{s:010}) and $\widehat{\Gamma}(m(\smallx))$ from (\ref{hat:Gamma}), and $\Omega^{\beta,\Gamma}_{\mbox{\tiny total mass}}:=\Omega^{\Gamma}_{\mbox{\tiny competition}}+\Omega^{\beta}_{\mbox{\tiny mass flow}}+\Omega^{\beta}_{\mbox{\tiny total mass fluctuation}}$ from (\ref{y:002}).

Abbreviate
\begin{equation}
\label{e:Gamma_int}
\begin{aligned}
   \Gamma\big(m(\smallx),r(\bar{\mu},x_2),\kappa(\bar{\mu}),\kappa_2\big)
   &:=
   \int\bar{\mu}(\mathrm{d}(x_1,\kappa_1))\,\Gamma\big(m(\smallx),r(x_1,x_2),\kappa_1,\kappa_2\big)
   \\
   &=
   \tfrac{1}{Nm(\smallx)}(B-C).
\end{aligned}
\end{equation}
Then by (\ref{e:exp-diff}),
\begin{equation}
\label{e:a(B-C)E}
\begin{aligned}
   &a\int\mathrm{d}\bar{\mu} (B-C) E
   \\
   &=g'(m(\smallx))\cdot m(\smallx)\cdot\int_{X \times K}\bar{\mu}(\mathrm{d}(x_2,\kappa_2))\,
   \Gamma\big(m(\smallx),r(\bar{\mu},x_2),\kappa(\bar{\mu}),\kappa_2\big)
   \\
   &\hspace{4cm}\cdot\int_{(X \times K)^n} \big(\tfrac{Nm(\smallx)\cdot \bar{\mu}-\delta_{(x_2,\kappa_2)}}{Nm(\smallx)-1}\big)^{\otimes n}
      (\mathrm{d}(\vx,\vk))\,\phi(\mr(\vx))\cdot f(\vk)
      \\
    &=g'(m(\smallx))\cdot m(\smallx)\cdot\hat{\Gamma}\big(m(\smallx)\big)\cdot\int_{(X \times K)^n}\bar{\mu}
    ^{\otimes n}
      (\mathrm{d}(\vx,\vk))\,\phi(\mr(\vx))\cdot f(\vk)
      \\
      &\;\;\;+g'(m(\smallx)) m(\smallx) \big(1 \vee \tilde{\gamma}\big(m(\smallx)\big)\big){\mathcal O}\big((Nm(\smallx))^{-1}\big)
      \\
    &=\Omega^{\Gamma}_{\mbox{\tiny competition}}F^{g,(n,\phi,f)}(\smallx) +{\mathcal O}\big((Nm(\smallx))^{-1}\big),
\end{aligned}
\end{equation}
as we have assumed that $m\tilde{\gamma}(m)g'(m)$ is bounded.

Further by (\ref{e:exp-diffq}),
\begin{equation}
\label{e:aA(D-E)}
\begin{aligned}
  &=a\int\mathrm{d}\bar{\mu} A(D-E)
  \\
  &=g'(m(\smallx))\cdot p\int_{X \times K}\bar{\mu}(\mathrm{d}(x_2,\kappa_2))\,\beta(\kappa_2)
  \int_K\alpha_N\big(\kappa_2,\mathrm{d}\tilde{\kappa}_2\big)
  \\
  &\;\;\cdot\int_{(X\uplus\{z\} \times K)^n}N m(\smallx)\cdot\Big(\big(\tfrac{Nm(\smallx)\bar{\mu}+\delta_{(z,\tilde{\kappa}_2)}}{Nm(\smallx) + 1}\big)^{\otimes n}-\big(\tfrac{Nm(\smallx)\bar{\mu}-\delta_{(x_2,\kappa_2)}}{Nm(\smallx) - 1}\big)^{\otimes n}\Big)
      (\mathrm{d}(\vx,\vk))\,\phi\big(\mr^{(x_2,z),\frac{1}{N}}(\vx)\big)\cdot f(\vk)
\\
  &\;\;\;+g'(m(\smallx))\cdot (1-p)\int_{X \times K}\bar{\mu}(\mathrm{d}(x_2,\kappa_2))\,\beta(\kappa_2)
   \\
  &\;\;\cdot\int_{(X \times K)^n}Nm(\smallx)\cdot \Big(\big(\tfrac{Nm(\smallx)\bar{\mu}+\delta_{(x_2,\kappa_2)}}{Nm(\smallx) + 1}\big)^{\otimes n}-\big(\tfrac{Nm(\smallx)\bar{\mu}-\delta_{(x_2,\kappa_2)}}{Nm(\smallx) - 1}\big)^{\otimes n}\Big)
      (\mathrm{d}(\vx,\vk))\,\phi\big(\mr(\vx)\big)\cdot f(\vk)
      \\
   &=2g'(m(\smallx))\int_{X \times K}\bar{\mu}(\mathrm{d}(x_2,\kappa_2))\,\beta(\kappa_2)\int_{(X\times K)^n}\sum_{l=1}^n\bar{\mu}^{\otimes (l-1)}\otimes\big(\delta_{(x_2,\kappa_2)}-\bar{\mu}\big)\otimes\bar{\mu}^{n-l}(\mathrm{d}(\vx,\vk))\,\phi\big(\mr(\vx)\big)\cdot f(\vk)
   \\
   &\;\;\;+g'(m(\smallx))o_N(1)+g'(m(\smallx)){\mathcal O}\big((Nm(\smallx))^{-1}\big)
       \\
   &=2g'(m(\smallx))\int_{(X\times K)^n}\bar{\mu}^{\otimes n}(\mathrm{d}(\vx,\vk))\, \sum_{l=1}^n\big(\beta(\kappa_l)-\hat{\beta}^{\smallx}\big) \phi\big(\mr(\vx)\big)\cdot f(\vk)+
   o_N(1)+{\mathcal O}\big((Nm(\smallx))^{-1}\big)
   \\
   &=\Omega^{\beta}_{\mbox{\tiny mass flow}}F^{g,(n,\phi,f)}(\smallx)+o_N(1)+{\mathcal O}\big((Nm(\smallx))^{-1}\big),
\end{aligned}
\end{equation}
as we have assumed $g$ is bounded differentiable.

Moreover,
\begin{equation}
\label{e:bA(D+E)}
\begin{aligned}
  &=b\int\mathrm{d}\bar{\mu} A(D+E)
  \\
  &=\tfrac{1}{2}g''(m(\smallx))\cdot m(\smallx)\cdot\int_{X\times K}\bar{\mu}(\mathrm{d}(x_2,\kappa_2))\,\beta(\kappa_2)
  \int_K\alpha_N\big(\kappa_2,\mathrm{d}\tilde{\kappa}_2\big)
  \\
  &\;\;\cdot \Big\{p\int_{(X\uplus\{z\} \times K)^n}\Big(\big(\tfrac{Nm(\smallx)\bar{\mu}+\delta_{(z,\tilde{\kappa}_2)}}{Nm(\smallx) + 1}\big)^{\otimes n}+\big(\tfrac{Nm(\smallx)\bar{\mu}-\delta_{(x_2,\kappa_2)}}{Nm(\smallx) - 1}\big)^{\otimes n}\Big)
      (\mathrm{d}(\vx,\vk))\,\phi\big(\mr^{(x_2,z),\frac{1}{N}}(\vx)\big)\cdot f(\vk)
      \\
  &\;\;\hspace{.5cm}+(1-p)\int_{(X \times K)^n}\Big(\big(\tfrac{Nm(\smallx)\bar{\mu}+\delta_{(x_2,\kappa_2)}}{Nm(\smallx) + 1}\big)^{\otimes n}+\big(\tfrac{Nm(\smallx)\bar{\mu}-\delta_{(x_2,\kappa_2)}}{Nm(\smallx) - 1}\big)^{\otimes n}\Big)
      (\mathrm{d}(\vx,\vk))\,\phi\big(\mr(\vx)\big)\cdot f(\vk)\Big\}
  \\
  &=g''(m(\smallx))\cdot m(\smallx)\cdot \hat{\beta}^{\smallx} \int_{(X\times K)^n}\bar{\mu}^{\otimes n} (\mathrm{d}(\vx,\vk))\,\phi\big(\mr(\vx)\big)\cdot f(\vk)+g''(m(\smallx))\cdot m(\smallx)\cdot o_N(1)
  \\
  &=\Omega^{\beta}_{\mbox{\tiny total mass fluctuation}}F^{g,(n,\phi,f)}(\smallx)+ o_N(1),
\end{aligned}
\end{equation}
as we have assumed that $mg''(m)$ is bounded.

Repeating the same arguments as in (\ref{e:leT11a}) and (\ref{e:geT11b}), (\ref{e:leT12}) and
(\ref{e:geT12}), as well as (\ref{e:leT13}) and (\ref{e:geT13}) we can conclude once more that
\begin{equation}
\label{e:leT2}
\begin{aligned}
   &\sup_{\smallx\in\mathbb{M}^K_N}\big|T_{2}-\Omega^{\beta,\Gamma}_{\mbox{\tiny total mass}}F^{g,(n,\phi,f)}(\smallx)\big|
   \tNo 0.
\end{aligned}
\end{equation}

\iffalse{
{Thus for some $C\in(0,\infty)$,
\begin{equation}
\label{e:leT2}
\begin{aligned}
   &\sup_{\smallx\in\mathbb{M}^K_N;\,m(\smallx)\le\frac{a_N}{N}}\big|T_{2}-\Omega^{\beta,\Gamma}_{\mbox{\tiny total mass}}F^{g,(n,\phi,f)}(\smallx)\big|
   \\
   &\le
    C\cdot\sup_{m\in[\frac{2}{N},\frac{a_N}{N}]}\tfrac{g''(m)}{N}\vee g'(m)\tNo 0,
\end{aligned}
\end{equation}
as $g''$ is bounded, and $g'(m)=0$. Moreover,
\begin{equation}
\label{e:geT2}
\begin{aligned}
   &\sup_{\smallx\in\mathbb{M}^K_N;\,m(\smallx)\ge\frac{a_N}{N}}\big|T_{2}-\Omega^{\beta,\Gamma}_{\mbox{\tiny total mass}}F^{g,(n,\phi,f)}(\smallx)\big|
   \\
   &\le
    C\cdot \frac{1}{a_N}\cdot\sup_{m\in[\frac{a_N}{N},\infty)}(m\cdot g'(m)+m\cdot g''(m))\tNo 0,
\end{aligned}
\end{equation}
as $\Omega F^{g,(n,\phi,f)}$ is bounded.
}}\fi
\smallskip

\bigskip

%.......................................................................

\noindent{\em Step~3 (The error-term $E_3$) }
Recall $E_3$ from (\ref{e:OmegaN1}) and $C_g$ from (\ref{Cg}). It turns out that the term $E_3$ is negligible in the limit as $N\to\infty$. Indeed, by (\ref{e:exp-diffq}) and Assumption~\ref{ass:004},
\begin{equation}
\label{e:E4bound}
\begin{aligned}
  E_3
 &=
   g'(m(\smallx)) \cdot\bar{\gamma}_b \cdot {\mathcal O}\big(\tfrac{1}{N}\big)
   +
    {\mathcal O}\big( \tfrac{1}{N} \big) g''(m(\smallx))m(\smallx)\big(\bar{\gamma}_b +\tilde{\gamma}(m(\smallx))\big)
   \\
   &+  m(\smallx)\cdot C_g(m(\smallx),N)\cdot \big( 1\vee{\tfrac{\tilde{\gamma}(m(\smallx))}{N}} \big)\cdot{\mathcal O}(\tfrac{1}{N}).
\end{aligned}
\end{equation}
Hence $E_3$ goes to zero as $N\to\infty$ uniformly over all $m\ge\frac{2}{N}$.
\end{proof}\sm

% ======================================================================

 \section{The compact containment condition}
 \label{S:Tightness}
Let for each $N\in\mathbb{N}$, ${\mathcal X}^N$ be the  tree-valued $(\varsig_N,\zeta_N,\alpha_N)$-trait-dependent branching dynamics with mutation and competition
rescaled as given in Subsection~\ref{Sub:rescale}.
In this section we  verify the following compact containment condition.

\begin{proposition}[Compact containment]
Assume that
$\{\mathcal X_0^N;\,N\in\mathbb{N}\}$ is a tight family of random elements in $\mathbb{M}^K$.
Suppose Assumptions~\ref{ass:002},~\ref{ass:007b},~\ref{ass:004},~\ref{ass:007c} and ~\ref{ass:007a} hold as well as (\ref{ass-ini-cond}) with $q=1$. Then for every $\epsilon, T>0$ there exists a compact $\ccset_{\epsilon,T} \subset \M^K$ such that
\begin{equation}
\label{equ:cc}
  \inf_{N \in \N} \PP( {\mathcal X}_t^N \in \ccset_{T,\epsilon} \mbox{ for all } t \in [0,T] ) \geq 1-\epsilon.
\end{equation}
\label{prop:cc}
\end{proposition}\sm

Recall from Proposition~\ref{pro-rel-comp} the characterization of (pre-)compact sets in $\mathbb{M}^K$.
To prove the claim of this proposition it is therefore enough to show the following.

\begin{lemma}[Compact containment modified]
Under the assumptions of Proposition~\ref{prop:cc}, for $T, \epsilon_0>0$ arbitrarily fixed, for all $k \in \N$ there exist sets $\ccset_{T,\epsilon_0,k} \subset \M^K$ such that conditions (i)--(iii) of Proposition~\ref{pro-rel-comp} with $\epsilon = 2^{-k}$ are satisfied for all $\smallx \in \ccset_{T,\epsilon_0,k}$ and
$$ \inf_{N \in \N} \PP( \{ {\mathcal X}_t^N \in \ccset_{T,\epsilon_0,k} \mbox{ for all } t \in [0,T] \} ) \ge 1-\epsilon_0 2^{-k}. $$
\label{lem-cc}
%
%\todo{Define $(\Omega,{\mathcal F},\mathbb{P})$ somewhere.}
%
\end{lemma}\sm

Indeed, set $\mathcal{K}_{T,\epsilon_0} := \bigcap_{k \in \N} \ccset_{T,\epsilon_0,k}$. Then the latter is relatively compact by the characterization of relatively compact sets in Proposition \ref{pro-rel-comp} and its closure satisfies (\ref{equ:cc}).

From now on assume that the assumptions of Proposition~\ref{prop:cc} hold. Let $T>0$ be arbitrarily fixed. We proceed to verify the compact containment condition from Lemma~\ref{lem-cc} in three steps. Let $\epsilon_0>0$ be arbitrarily fixed. Item~(i) of Proposition \ref{pro-rel-comp} (control of total mass) is established first. Items~(ii) and~(iii-a)--(iii-b) (control of traits and the diameter) follow in Step~2. The proof of Item (iii-c) (coverage number by balls of radius $3\epsilon$) is the most involved and concludes the proof in a third step. As a preparation we will start to provide a comparison argument which will allow for reductions to technically simpler situations.
\begin{lemma}[Comparison]
Suppose that $(\beta,\gamma^{\mathrm{birth}}_1,\gamma^{\mathrm{death}}_1)$ and $(\beta,\gamma^{\mathrm{birth}}_2,\gamma^{\mathrm{death}}_2)$ satisfy the assumptions of Proposition~\ref{prop:cc}, $\gamma^{\mathrm{birth}}_2(m,r,\kappa_1,\kappa_2)=\gamma^{\mathrm{birth}}_2(\kappa_2)$, $\gamma^{\mathrm{death}}_2(m,r,\kappa_1,\kappa_2)=\gamma^{\mathrm{death}}_2(\kappa_2)$ and suppose that $\gamma^{\mathrm{birth}}_1\le \gamma^{\mathrm{birth}}_2$ and that $\gamma^{\mathrm{death}}_2 \le \gamma^{\mathrm{death}}_1$. Assume that Lemma~\ref{lem-cc} holds
with $(\beta,\gamma^{\mathrm{birth}}_2,\gamma^{\mathrm{death}}_2)$. Then it also holds with
$(\beta,\gamma^{\mathrm{birth}}_1,\gamma^{\mathrm{death}}_1)$.
\label{lem-cc-red}
\end{lemma}

\begin{proof} Fix $p\in[0,1]$, and couple for each $N\in\mathbb{N}$, $\mathcal{X}^{N,1}$
with competition rate $\gamma^{\mathrm{death}}_1$ and  birth-enhancement rate $\gamma^{\mathrm{birth}}_1$ to  $\mathcal{X}^{N,2}$
with competition rate $\gamma^{\mathrm{death}}_2$ and  birth-enhancement rate $\gamma^{\mathrm{birth}}_2$
both starting in the same state
in such a way that for all realizations and all $t\ge 0$, $X_t^{N,1} \subseteq X_t^{N,2}$, $r^{N,1}_t=r^{N,2}_t |_{X^{N,1}_t}$ and $\mu^{N,1}_t\le \mu^{N,2}_t$.
Such a coupling is always possible as $\mathcal{X}^{N,2}$ can only result in additional births or omissions of deaths that are captured by means of $\mu_t^{N,2}-\mu_t^{N,1}$ and have no impact on the death- respectively birth-rates in (\ref{s:002a}) respectively (\ref{s:002d}).
\end{proof}\sm

Recall that by Remark~\ref{Rem:003} there exists a finite constant $C>0$ such that $\gamma^{\mathrm{birth}}(m,r,\kappa_1,\kappa_2)\le C\beta(\kappa_2)$ for all $m\ge 0$, $r\ge 0$, $\kappa_1,\kappa_2\in K$.
For the remainder of this section, assume that without loss of generality for  $m\ge 0$, $r\ge 0$ and $\kappa,\kappa'\in K$,
\begin{equation}
\label{e:dominating-rates}
  \gamma^{\mathrm{death}}(m,r,k,k') \equiv 0 \mbox{ and } \gamma^{\mathrm{birth}}(m,r,k,k') \equiv C\beta(\kappa').
\end{equation}
In particular, the model now satisfies the additional assumptions of Proposition~\ref{P:005}(iii). \sm

%%%%%%%%%%%%%%%%%%%%%%%%%%%%%%%%%%%%%%%%%%%%%%%%%%%%%%%

\noindent{\bf Step 1: The mass [Item (i) of Proposition \ref{pro-rel-comp}].}
For all $T,\epsilon_0>0$, $k \in \N$ there exists $M_k>0$ big enough such that
\begin{equation}
\label{e:mass-tight}
  \sup_{N \in \N} \PP(\sup_{t \in [0,T]} m(\smallx_t^N) > M_k) < \tfrac{1}{3} \epsilon_0 2^{-k}.
\end{equation}
This follows from (\ref{e:t-moment-bounds-1}) and Markov's inequality. \bigskip

%%%%%%%%%%%%%%%%%%%%%%%%%%%%%%%%%%%%%%%%%%%%%%%%%%%%%%%

\begin{remark}[Excluding empty clans]
Recall the remarks in and around (\ref{s:002aa})--(\ref{s:005}) with regards to empty clans, that is, clans $x \in X$ that satisfy $n_x=0=\mu(\{x\} \times K)$. For the remainder of this section, we assume without loss of generality that every $x \in X_t^N$ satisfies $\mu(\{x\} \times K)>0$. Indeed, it remains to verify Items (ii) and (iii) of Proposition~\ref{pro-rel-comp}. Empty clans are without effect on Item (ii) and can be included into the exceptional sets $X_\epsilon^c$ in Item (iii).
\end{remark}
 
For $(X,r,\mu)\in\mathbb{M}^K$ let
\begin{equation}
\label{e:cvernr}
\CS_\epsilon(\smallx) := \min\big\{ \# \mbox{ balls of radius } \epsilon \mbox{ needed to cover } \mathrm{supp}\big(\mu(\boldsymbol{\cdot}\times K)\big) \big\}.
\end{equation}

As we assume that $\{{\mathcal X}_0^N;\,N\in\mathbb{N}\}$ is a tight family of random elements in $\mathbb{M}^K$, Proposition~\ref{P:005}(iii) yields for each $\epsilon_0>0$ and $k\in\mathbb{N}$ the existence of
 $m_k=m_k(\epsilon_0)>0$ such that
\begin{equation}
\label{e:2aaaiii}
   \sup_{n\in\mathbb{N}}\mathbb{P}\big(\sup_{t\in[0,T]}m({\mathcal X}_t^N)>2^{-k}/4\big|m({\mathcal X}_0^N) \leq m_k\big)\le\tfrac{1}{24}\epsilon_02^{-k}.
\end{equation}
Moreover, 
Proposition~\ref{pro-rel-comp} yields the existence of $L_k=L_k(\epsilon_0), N_k=N_k(\epsilon_0)>0$ big enough such that 
\begin{equation}
\label{e:2aa}
\begin{aligned}
  & \inf_{N \in \N} \PP\Big(\exists (X_0^N)_k\subseteq X_0^N: \mu_0^N((X_0^N)_k^c \times K) \leq m_k, \ \mathrm{diam}((X^N_0)_k) \leq \tfrac{L_k}{2} \ \mbox{ and } \ \CS_{2^{-k}}((\smallx_0^N)_k) \leq\tfrac{N_k}{2} \Big) \\
  &
  \;> 1- \tfrac{1}{24} \epsilon_0 2^{-k}.
\end{aligned}
\end{equation}

%%%%%%%%%%%%%%%%%%%%%%%%%%%%%%%%%%%%%%%%%%%%%%%%%%%%%%%

\noindent{\bf Step 2: The trait and the diameter [Item (ii), (iii-a,b) of Proposition \ref{pro-rel-comp}].}
In this step we prove the following (recall that we consider $\epsilon=2^{-k}$ in (i)--(iii) of Proposition \ref{pro-rel-comp}):

\begin{lemma}
\label{lem-tight-step-2}
For all $k \in \N$ there exists $K_k \subseteq K$ compact and $L_k>0$ big enough such that
\begin{equation}
\label{e:tight-ii}
  \inf_{N \in \N} \PP_N\Big( \sup_{t \in [0,T]} \mu_t^N(X_t^N \times K_k^c) \leq 2^{-k} \Big) > 1- \tfrac{1}{6} \epsilon_0 2^{-k}
\end{equation}
and
\begin{equation}
\begin{aligned}
  \inf_{N \in \N} \PP\big( \forall t \in [0,T], \ \exists (X_t^N)_k \subseteq X_t^N:\;\;
  \sup_{t \in [0,T]} \mu_t^N(((X_t^N)_k)^c \times K) \leq 2^{-k}/2, &
  \\
  \sup_{t \in [0,T]} \mathrm{diam}((X_t^N)_k) \leq L_k \big) &> 1- \tfrac{1}{6} \epsilon_0 2^{-k}.
\label{e:tight-iiib}
\end{aligned}
\end{equation}
\end{lemma}

\begin{proof}[Proof of Lemma~\ref{lem-tight-step-2}.] Consider the individual-based description of our
trait-dependent branching particle model with mutation and competition given in Subsection~\ref{Sub:branchmodel}. Notice that due to the branching mechanism, each particle $x$ alive at time $t\ge 0$ is the descendant of a well-defined ancestor, $A_{t-s}(x)$, at time $0\le s\le t$.
If
$A_{t}(x)=A_{t}(y)$, let
$T(x,y):=\sup\{s\in(0,t]:\,A_{t-s}(x)=A_{t-s}(y)\}$
denote the time instant at which the ancestral lines of these particles had split.
Extend this definition to letting $T(x,y):=0$ in case that $A_{t}(x)\not=A_{t}(y)$.
If we define the {\em (genetic) age} $a^N_t(\{x\})$ of an individual $x$ living at time $t$  as $\zeta_N \big(=\tfrac{1}{N} \big)$ times the number of mutations which separate $A_t(x)$ and $x$,
then
the genetic distance of particles $x$ and $y$ satisfies
\begin{equation}
\label{rtat}
\begin{aligned}
   &r^N_t\big(x,y\big)
   \\
   &=a^N_t\big(\{x\}\big)-a^N_{T(x,y)}\big(\{A_{(t-T(x,y))}(x)\}\big)
   +a^N_t\big(\{y\}\big)-a^N_{T(x,y)}\big(\{A_{(t-T(x,y))}(y)\}\big)+r^N_0\big(A_{t}(x),A_{t}(y)\big)
   \\
   &\le
   a^N_t\big(\{x\}\big)+a^N_t\big(\{y\}\big)+r^N_0\big(A_{t}(x),A_{t}(y)\big).
   \end{aligned}
\end{equation}

\iffalse{By (\ref{e:006}), we can choose $m=m(\epsilon_0,k)>0$ suitably small such that
\begin{equation}
\label{e:2aaa}
   \sup_{n\in\mathbb{N}}\mathbb{P}\big(\sup_{t\in[0,T]}m(\smallx_t^N)>\varepsilon_02^{-k}|m(\small_0^N)=m\big)\le 1-\epsilon_02^{-k}.
\end{equation}
As by assumption of Proposition~\ref{prop:cc}, $\{\mathcal X_0^N;\,N\in\mathbb{N}\}$ is a tight family of random elements in $\mathbb{M}^K$,  Proposition~\ref{pro-rel-comp} yields for $\epsilon=2^{-k}$ the existence of $(X_0^N)_k\subseteq X_0^N$ with $\mu_0^N(((X_0^N)_k)^c)\le m$ and $L_k=L_k(\epsilon_0)>0$ big enough such that
\begin{equation}
\label{e:2aa}
  \inf_{N \in \N} \PP\Big( \mathrm{diam}((X_0^N)_k) > L_k/2 \Big) \leq \tfrac{1}{12} \epsilon_0 2^{-k}.
\end{equation}}\fi

Fix
%$T>0$, $\epsilon_0>0$ and
$k\in\mathbb{N}$. Recall the notation $n_x$ for the number of clones in a clan $x$ from (\ref{e:003a}). We consider the following two auxiliary measure-valued processes $\chi, \eta$ which are Markov jump processes taking values in the space of finite measures on $K$ respectively $K\times\R_+$ and are given for each $t\ge 0$ by
\begin{equation}
\label{etat}
  \chi^N_t:=\tfrac{1}{N}\sum_{x\in X_t^N}n_x\delta_{\kappa_t(x)}
  \quad \mbox{ and } \quad
  \eta^N_t:=\tfrac{1}{N}\sum_{x\in X_t^N; \,A_t(x) \in (X_0^N)_k} n_x\delta_{(\kappa_t(x),a^N_t(\{x\}))}.
\end{equation}
with $(X_0^N)_k$ as in (\ref{e:2aa}).

As by Assumption (\ref{e:dominating-rates}) none of the rates depend on the genetic distance, we obtain the following:
\begin{lemma} Assume that the assumptions of Proposition~\ref{prop:cc} hold.
If the families $\{\chi^N;\,N\in\mathbb{N}\}$ and $\{\eta^N;\,N\in\mathbb{N}\}$ satisfy respective compact containment conditions, then the statement of Lemma~\ref{lem-tight-step-2} follows.
\label{lem-measure-valued}
\end{lemma}

\begin{proof}[Proof of Lemma~\ref{lem-measure-valued}] Fix $T>0$, and $\epsilon_0>0$. By the tightness assumptions and \cite[Theorem~3.6.3 and Remark~3.6.4]{EK86}, for any $k\in\mathbb{N}$, we can find a compact subset $K_k\subseteq K$ and a number $L_k>0$ such that
\begin{equation}
\label{eta-cc-1}
  \sup_{N\in\mathbb{N}}\mathbb{P}\big(\sup_{t\in[0,T]}\chi^N_t\big({K}_k^c\big)> 2^{-k}\big)\le \tfrac{1}{12}\epsilon_02^{-k}.
\end{equation}
and
\begin{equation}
\label{eta-cc}
  \sup_{N\in\mathbb{N}}\mathbb{P}\big(\sup_{t\in[0,T]}\eta^N_t\big(({K}_k\times[0,\tfrac{L_k}{4}])^c\big)> 2^{-k}/2\big)\le \tfrac{1}{12}\epsilon_02^{-k}.
\end{equation}
By (\ref{eta-cc-1}) and the definition of $\chi$,
\begin{equation}
\label{xtn-cc-1}
  \sup_{N\in\mathbb{N}}\mathbb{P}\big(\sup_{t\in[0,T]}\mu_t^N(X_t^N \times K^c_k)> 2^{-k}\big)
  \le \tfrac{1}{12}\epsilon_0 2^{-k}
\end{equation}
and (\ref{e:tight-ii}) follows.

Next set for each $t\in[0,T]$, 
\begin{equation}
\label{e:XtNk-Step2}
   (X_t^N)_k
 := 
   \big\{ x \in X_t^N:\, A_t(x) \in (X_0^N)_k \mbox{ and } a_t^N(\{x\}) \leq L_k/4 \big\}.
\end{equation}
Then $\mathrm{diam}((X^N_t)_k) \leq L_k$ and by (\ref{e:2aaaiii}), (\ref{e:2aa}), (\ref{eta-cc}) and the definition of $\eta$,
\begin{equation}
\label{xtn-cc}
\begin{aligned}
  &\sup_{N\in\mathbb{N}}\mathbb{P}\big(\sup_{t\in[0,T]}\mu_t^N((X_t^N)^c_k\times K) > 2^{-k}/2 \ \mbox{ or }\sup_{t\in[0,T]}\mathrm{diam}((X^N_t)_k)> L_k\big)
  \\
  &\le \tfrac{1}{24}\epsilon_0 2^{-k} + \tfrac{1}{24}\epsilon_0 2^{-k} + \sup_{N\in\mathbb{N}}\mathbb{P}\big(\sup_{t\in[0,T]}\eta^N_t\big(({K}_k\times[0,\tfrac{L_k}{2}])^c\big)> 2^{-k}/2\big)
   \\
   &\le \tfrac{1}{6}\epsilon_0 2^{-k},
\end{aligned}
\end{equation}
and the second claim (\ref{e:tight-iiib}) follows.
\end{proof}

To finish the proof of  Lemma~\ref{lem-tight-step-2}, we need to verify the compact containment for the measure-valued processes.
Recall (\ref{e:dominating-rates}).
In this particular setting the claim is covered by \cite[Theorem~3.4]{Kli14} (applied to $r(t,y):=\beta(\pi_K(y_t))$ with $\pi_K:K\times\R_+\to K$ denoting the projection map on to the trait, $b(t,y):=C\beta(\pi_K(y_t))$, $D(t,y)\equiv 0$, $U(t,y,y')\equiv 0$ and $\alpha_N((\kappa,a),\mathrm{d}(\kappa',a')):=\alpha_N(\kappa,\mathrm{d}\kappa')
\delta_{a+\frac{1}{N}\mathbf{1}\{\kappa\not =\kappa'\}}(da')$. Here the left hand sides refer to the set-up used in \cite{Kli14} and the right hand sides to our set-up.).
\end{proof}

\begin{remark}[Non-ultrametric setup]
In the present article and \cite{Kli14}, because of the use of exponential times in the modelling of birth- and death-events, the analysis of the modulus of continuity of the trait-history of a particle (cf. \cite[Definition~5.3]{Kli14}) plays a major role in obtaining appropriate bounds. In the present article, the need for such an analysis arises due to the non-ultrametric setup, where instead of defining genetic distance to be twice the time to the most recent common ancestor (MRCA) (cf., for example, \cite[(2.27)]{GPW13} in the context of tree-valued Moran dynamics), genetic distance is taken to be the sum of the number of mutations of two individuals back in time to their MRCA. The control of the modulus of continuity along the path of an individual now allows us to relate genetic distance to time by considering age as part of the type-space in \cite{Kli14}. \qed
\end{remark}

\noindent{\bf Step~3: Coverage number of balls of radius $3\epsilon$ [Item (iii-c) of Proposition \ref{pro-rel-comp}].}
In this step we show the following.
\begin{lemma}
For all $k \in \N$ there exists $N_k>0$ such that
\begin{equation}
\begin{aligned}
  & \inf_{N \in \N} \PP\Big( \forall t \in [0,T], \ \exists (X_t^N)'_k \subseteq X_t^N: \sup_{t \in [0,T]} \mu_t^N(((X_t^N)'_k)^c \times K) < 2^{-k}/2,
  \\
  & \qquad\quad\ \sup_{t \in [0,T]} \mathcal{S}_{3 \cdot 2^{-k}}((\smallx_t^N)'_k) < N_k \Big) > 1- \tfrac{1}{3} \epsilon_0 2^{-k},
\end{aligned}
\end{equation}
where $(\smallx_t^N)_k' = \overline{\big((X_t^N)_k', r_t^N |_{(X_t^N)_k' \times (X_t^N)_k'}, \mu_t^N( \cdot |_{(X_t^N)_k' \times K}) \big)}$.
\label{lem-cc-coverage}
\end{lemma}

\begin{proof} 

Fix $k\in\N$. Choose $m_k, (X_0^N)_k$ and $N_k$ according to (\ref{e:2aaaiii}) and (\ref{e:2aa}).

We first argue that there exists $t_0=t_0(\epsilon_0,k)>0$ small enough such that for all $0 \leq t \leq t_0$ the genetic distances of a large enough proportion of particles alive at time $t$ to their ancestors at time $0$ are at most $\epsilon=2^{-k}$ with high probability, that is,  such that
\begin{equation}
\begin{aligned}
  & \inf_{N \in \N} \PP\Big( \forall t \in [0,t_0], \ \exists (\tilde{X}_t^N)_k \subseteq X_t^N: \sup_{t \in [0,t_0]} \mu_t^N(((\tilde{X}_t^N)_k)^c \times K) < 2^{-k}/4,
  \\
  & \qquad\quad\ \sup_{t \in [0,t_0]}\max_{x\in \tilde{X}_t^N} a^N_t(\{x\})<2^{-k} \Big) > 1- \tfrac{1}{12} \epsilon_0 2^{-k}.
  \label{e:tight-iv}
\end{aligned}
\end{equation}
This follows indeed immediately from \cite[Lemma~3.9]{Kli14}.

With this choice of $t_0$, we can conclude from (\ref{e:2aaaiii}) (recall that $m(\smallx)=\mu(X \times K)$) and (\ref{e:2aa}) by the triangle inequality that
\begin{equation}
\begin{aligned}
  & \inf_{N \in \N} \PP\Big( \forall t \in [0,t_0], \ \exists (\tilde{X}_t^N)_k \subseteq X_t^N: \sup_{t \in [0,t_0]} \mu_t^N(((\tilde{X}_t^N)_k)^c \times K) < 2^{-k}/2,
  \\
  & \qquad\quad\ \sup_{t \in [0,t_0]} \CS_{3 \cdot 2^{-k}}((\tilde{\smallx}_t^N)_k) \leq N_k/2 \Big) > 1-\tfrac{1}{6} \epsilon_0 2^{-k}.
\label{e:time-zero-coverage-claim}
\end{aligned}
\end{equation}

It then remains to show that
\begin{equation}
\begin{aligned}
  & \inf_{N \in \N} \PP\Big( \forall t \in [t_0,T], \ \exists (\tilde{X}_t^N)_k \subseteq X_t^N: \sup_{t \in [t_0,T]} \mu_t^N(((\tilde{X}_t^N)_k)^c \times K) < 2^{-k}/2,
  \\
  & \qquad\quad\ \sup_{t \in [t_0,T]} \CS_{3 \cdot 2^{-k}}((\tilde{\smallx}_t^N)_k) \leq N_k/2 \Big) > 1-\tfrac{1}{6} \epsilon_0 2^{-k}.
\label{e:tight-v}
\end{aligned}
\end{equation}

%As we have assumed (\ref{e:dominating-rates}), if we see a branching event, the probability that it leads to birth rather than death of an individual does not depend on the individual's type. Thus in this step we can assume without loss of generality that
%\begin{equation}\label{e:betacons}
%  \beta\equiv\underline{\beta}.
%\end{equation}
%Indeed we can couple both models such that the birth events and death events match, while the times between two successive branching events are longer for the lower bound $\underline{\beta}$ than for the type depending $\beta(\boldsymbol{\cdot})$.

We will first prove the statement in case the genetic distance $r$ (counting the number of substitutions due to mutation) is replaced by the genealogical distance $\tilde{r}$ (which is twice the time to the most recent ancestor), that is, in contrast to (\ref{s:006}) distances
(between any two distinct individuals) grow deterministically with time at speed $2$.
In the set-up of genealogical distance our model fits into the class of {\em finite population models} (compare, \cite[Definition~2.18]{GPW13}).
Put for each $t\in[0,T]$, %and $k\in\mathbb{N}$,
\begin{equation}
\label{e:XtNk}
   (X_t^N)_k
 :=
   \big\{x\in X_t^N:\,A_t(x)\in (X_0^N)_k\big\},
\end{equation}
the set of descendants of $(X_0^N)_k$ at time $t$. It then remains to show
the following:
\begin{lemma} For all $0 < \xi \leq t < T$, the family $\{\CS_{2\xi}(((X_t^N)_k,\tilde{r},\mu_t^N( \cdot |_{(X_t^N)_k \times K})));\,N\in\mathbb{N}\}$ is tight.
\label{lemma-ancestors-tight-single-time}
\end{lemma}\smallskip

\begin{remark} Assume that $\{\CS_{2\xi}(((X_t^N)_k,\tilde{r},\mu_t^N( \cdot |_{(X_t^N)_k \times K})));\,N\in\mathbb{N}\}$ is tight for all $0 < \xi \leq t < T$, then we can literally copy the proof of \cite[Lemma~6.7a)]{GPW13} to obtain for all $0 < \xi \leq t < T$, the existence of $C_\xi$
such that
\begin{equation}\label{e:suptight}
  \sup_{N \in \N} \PP\big( \sup_{t \in [\xi,T)} \CS_{2 \xi}(((X_t^N)_k
  ,\tilde{r},\mu_t^N( \cdot |_{(X_t^N)_k \times K}))) > C_\xi \big) \leq 2 \xi.
\end{equation}
\label{Rem:004}
\end{remark}\smallskip

\begin{proof}[Proof of Lemma~\ref{lemma-ancestors-tight-single-time}]
Until the end of this proof we set $(\tilde{\smallx}_t^N)_k := \overline{((X_t^N)_k,\tilde{r},\mu_t^N( \cdot |_{(X_t^N)_k \times K}))}$. Fix $\xi \leq t < T$. It is enough to show that for all $\delta>0$ one can find $C=C(\delta)>0$ such that
\begin{equation}
\label{e:bd-0-dom-anc}
  \sup_{N \in \N} \PP( \CS_{2\xi}((\tilde{\smallx}_t^N)_k) > C ) < \delta.
\end{equation}

By reasoning as in the derivation of (\ref{e:mass-tight}) we can find a constant $\tilde{M}_\xi$ such that
\begin{equation}
\label{e:tightmass}
   \sup_{N \in \N} \PP(\sup_{t \in [0,T]} m((\tilde{\smallx}_t^N)_k) > \tilde{M}_\xi) < \delta/2
\end{equation}
holds and thus for $t \in [\xi,T)$ arbitrary, using Markov's inequality, we obtain
\begin{equation}
\begin{aligned}
  \PP( \CS_{2\xi}((\tilde{\smallx}_t^N)_k) > C )
  &\leq \tfrac{\delta}{2} + \PP\big( \CS_{2\xi}((\tilde{\smallx}_t^N)_k) > C, m((\tilde{\smallx}_{t-\xi}^N)_k) \leq \tilde{M}_\xi \big)
  \\
  &\leq \tfrac{\delta}{2} + \tfrac{1}{C} \mathbb{E}\big[ \CS_{2\xi}((\tilde{\smallx}_t^N)_k) \1_{\{ m((\tilde{\smallx}_{t-\xi}^N)_k) \leq \tilde{M}_\xi \}} \big].
  \label{e:bd-2-dom-anc}
\end{aligned}
\end{equation}

As we have assumed (\ref{e:dominating-rates}), if we see a branching event, the probability that it leads to birth rather than death of an individual does not depend on the individual's type. Moreover, by (\ref{s:002a}) respectively (\ref{s:002d}) the total death- respectively birth-rate of a particle only depends on its clan's trait. As a result, the progenies of different clans evolve independently from each other and we obtain with (\ref{e:003b}),
\begin{equation}
\label{e:bd-3-dom-anc}
\begin{aligned}
  & \mathbb{E}\big[ \CS_{2\xi}((\tilde{\smallx}_t^N)_k) \1_{\{ m((\tilde{\smallx}_{t-\xi}^N)_k) \leq \tilde{M}_\xi \}} \big] \\
  & \leq N\tilde{M}_\xi \cdot \sup_{\kappa \in K} \mathbb{P}\big( \mbox{a particle starting with trait } \kappa \mbox{ has a child present after a period of time } \xi \big).
\end{aligned}
\end{equation}
We can now assume without loss of generality that $\beta(\cdot) \equiv \underline{\beta}$ with $\underline{\beta}>0$ as in Assumption~\ref{ass:007c}. Indeed, we can couple both models such that the birth events and death events match, while the times between two successive branching events are longer for the lower bound $\underline{\beta}$ than for the type dependent $\beta(\cdot)$. As a result, the progeny of a particle $x$ gets extinct later in the model using $\underline{\beta}$. Apply the theory of continuous-time birth- and death-processes to the $N^{\mathrm{th}}$-model with birth- respectively death-rates of $N \underline{\beta}$ respectively $N \underline{\beta} (1+ C/N)$ to conclude that
\begin{equation}
\label{e:bd-3-dom-anc-2}
  \mathbb{E}\big[ \CS_{2\xi}((\tilde{\smallx}_t^N)_k) \1_{\{ m((\tilde{\smallx}_{t-\xi}^N)_k) \leq \tilde{M}_\xi \}} \big] 
  \leq N \tilde{M}_\xi \cdot \mbox{const.} N^{-1}.
\end{equation}

Now choose $C$ big enough in (\ref{e:bd-2-dom-anc}) to obtain (\ref{e:bd-0-dom-anc}).
\end{proof}\smallskip

We next relate time to genetic distance in order to obtain a bound on $\CS_{2\xi}(((X_t^N)_k,r,\mu_t^N( \cdot |_{(X_t^N)_k \times K})))$
from the bound on $\CS_{2\xi}(((X_t^N)_k,\tilde{r},\mu_t^N( \cdot |_{(X_t^N)_k \times K})))$. In what follows we estimate the mass of the exceptional subset of $(X_t^N)_k$ of those individuals whose ancestors at time $t-\xi$ are in genetic distance greater than $\epsilon= 2^{-k}$.

 \begin{lemma}[Too fast evolving individuals]
For all $0<t_0<T, \sigma>0$ and $k \in \N$ there exists $0<s_0=s_0(t_0,T,\sigma,k)<t_0/2$ small enough such that
\begin{equation}
  \sup_{N \geq 3 \cdot 2^k} \PP\big( \exists t \in [t_0,T]: \mu_t^N \big( x \in X_t^N: a_t^N(\{x\})-a_{t-s_0}^N(A_{t-s_0}(\{x\},t)) \in [2^{-k},\infty) \big) \geq 2^{-k}/4 \big) \leq \sigma,
\label{e:too-fast-is-J-set}
\end{equation}
where we denote by $A_s(\{x\},t) \in X_s^N$ the clan that contains the ancestors at time $0 \leq s \leq t$ of all the individuals in the clan $x \in X_t^N$ at time $t$.
\label{lemma-too-fast-is-J-set}
\end{lemma}\smallskip

\begin{proof} This is covered by \cite[Lemma~3.9]{Kli14}, which gives the analogous result phrased in terms of modulus of continuity of the historical process associated with $\eta^N$ from (\ref{etat}).
As the maximal jump-size in age is $1/N$, \eqref{e:too-fast-is-J-set} follows.
\end{proof}\smallskip

 We now combine Lemma~\ref{lemma-too-fast-is-J-set} with Remark~\ref{Rem:004} to finish the proof of Lemma~\ref{lem-cc-coverage}.
 %Fix  $T>0$, $\epsilon_0>0$ and $k\in\mathbb{N}$, and put $\epsilon:=2^{-k}$.

{\it 1. Case: $N < 3 \cdot 2^k$.}

Recall that $\epsilon=2^{-k}$. Choose $(X_t^N)_k = X_t^N$. By reasoning as in (\ref{e:mass-tight}), we can choose $M'_k > 0$ big enough such that
\begin{equation}
\label{e:case1-of-step3}
  \sup_{N \in \N} \PP(\sup_{t \in [0,T]} m(\smallx_t^N) > M'_k) < \tfrac{1}{12} \epsilon_0 2^{-k}. 
\end{equation}
For $N < 3 \cdot 2^k$ the number of individuals alive at time $t$ is $N m(\smallx_t^N)$ and thus
\begin{equation}
\label{e:case1-of-step3-2}
  \sup_{t \in [t_0,T]} \mathcal{S}_{2^{-k}}(\smallx_t^N) < 3 \cdot 2^k M'_k 
\end{equation}
with probability greater or equal to $1-\tfrac{1}{12} \epsilon_0 2^{-k}$. Now (\ref{e:tight-v}) holds with $N_k/2 = 3 \cdot 2^k M'_k$. It remains to investigate:

{\it 2. Case: $N \geq 3 \cdot 2^k$.}

 Recall the constants $t_0=t_0(k)$ from (\ref{e:tight-iv}).
 For $t \in [0,T]$, let
\begin{equation}
\label{e:def-X-k}
  (X_t^N)'_k := ((X_t^N)_k \setminus \big\{x \in X_t^N: a_t^N(\{x\}) - a_{t-s_0}^N(A_{t-s_0}(\{x\},t)) \in [2^{-k},\infty) \big\}.
\end{equation}

By Remark~\ref{Rem:004} we obtain for $\xi=s_0 \wedge \tfrac{1}{48} \epsilon_0 2^{-k}$,
\begin{equation}
\label{e:ancestor-bounds}
  \sup_{N \in \N} \PP^N\big( \sup_{t \in [t_0,T)} \CS_{2 s_0} (((X_t^N)'_k,\tilde{r},\mu_t^N( \cdot |_{(X_t^N)'_k \times K}))) > C_\xi \big) \leq \tfrac{1}{24} \epsilon_0 2^{-k}
\end{equation}
for some constant $C_{\xi}$. An $\epsilon$-coverage with $\epsilon=2^{-k}$ of $(X_t^N)'_k$ for $t \in [t_0,T]$ can now be constructed as follows: Go back $0<s_0<t_0/2$ in time, determine the number of ancestors of $(X_t^N)'_k$, which is smaller or equal to $C_{\xi}$ with probability at least $1-\tfrac{1}{24} \epsilon_0 2^{-k}$ by (\ref{e:ancestor-bounds}). Then the genetic distance of the particles in $(X_t^N)'_k$ to their ancestors at time $t-s_0$ is indeed bounded by $\epsilon=2^{-k}$ by the definition of $(X_t^N)'_k$ in (\ref{e:def-X-k}). Lemma~\ref{lemma-too-fast-is-J-set} for the choice $\sigma=\tfrac{1}{24} \epsilon_0 2^{-k}$ further yields that  \eqref{e:tight-v} holds
with probability at least $1-\tfrac{1}{6} \epsilon_0 2^{-k}$. We obtain (\ref{e:tight-v}) with $N_k/2=\max\{ 3 \cdot 2^k M'_k, C_\xi \}$.
\end{proof}\smallskip

{\it Completion of the proof of Lemma~\ref{lem-cc} and thereby of Proposition~\ref{prop:cc}.}
Combine (\ref{e:mass-tight}), Lemma~\ref{lem-tight-step-2} and Lemma~\ref{lem-cc-coverage}, replacing $(X_t^N)_k$ by the intersection of the respective sets and choosing $N_\epsilon = \max\{ M_k, L_K, N_k \}$ for $\epsilon=3 \cdot 2^{-k}$ in Proposition~\ref{pro-rel-comp}. \qed

\section{Finishing up (Proofs of Theorems~\ref{T:tightness} and~\ref{T:mp})}
\label{S:proofs}
Let for each $N\in\mathbb{N}$, ${\mathcal X}^N$ be the  tree-valued $(\varsig_N,\zeta_N,\alpha_N)$-trait-dependent branching dynamics with mutation and competition
rescaled as given in Subsection~\ref{Sub:rescale}.
In this section we want to collect the results from the previous sections to present the proof of tightness (Theorem~\ref{T:tightness}) and the fact that any limit process satisfies the $(\Omega,{\mathcal D}(\Omega))$-martingale problem (Theorem~\ref{T:mp}) with $(\Omega,{\mathcal D}(\Omega))$ as in (\ref{e:DOmega})--(\ref{y:001}).

We will apply \cite[Remark~4.5.2]{EK86}. This remark states that the sequence $({\mathcal X}^N)_{N\in\mathbb{N}}$ is relatively compact if
\begin{itemize}
\item[(i)] the (strong) compact containment condition (\ref{equ:cc}) holds,
\item[(ii)] the closure of $\widetilde{\Pi}$ (cf. (\ref{e:tildePi})) contains an algebra which separates points in $\mathbb{M}^K$ and vanishes nowhere, and
\item[(iii)] for all $F\in\widetilde{\Pi}$, (\ref{e:021}) holds.
\end{itemize}
In addition, any limit process satisfies the $(\Omega,\widetilde{\Pi})$-martingale problem.

As (i) is covered by Proposition~\ref{prop:cc} and (iii) is covered by Proposition~\ref{P:generators} ($\tilde{\gamma}(m)=(1 \vee m) \overline{\gamma}_d$ by Assumption~\ref{ass:007a}) it therefore remains to show the following two facts.
 \begin{proposition} The subspace $\widetilde{\Pi}\subseteq{\mathcal D}(\Omega)$ separates points and vanishes nowhere. Moreover, its closure contains an algebra.
\label{P:006}
\end{proposition}\sm

\begin{proof} Recall that
a class $\smallx:=\overline{(X,r,\mu)}\in\mathbb{M}^K$ is uniquely characterized by
the total mass $m(\smallx)$
and the {\em marked distance matrix distribution}
$\nu^{\smallx}\in
\mathcal{M}_1\big(\R_+^{\binom{\N}{2}}\times K^{\mathbb{N}}\big)$ (cf. (\ref{eq:spaces:distance-matrix-distribution})). Thus a separating set for
$\mathbb{M}^K$ are clearly the bounded continuous  functions of the form
$e^{-\lambda m(\smallx)}\int\nu^{\smallx}(\mathrm{d}(\mr,\vk))\,\phi\big(\mr\big)f(\vk)$ for some $n\in\mathbb{N}$, $\lambda\ge 0$, $\phi\in{\mathcal C}_b(\R^{n\choose 2})$ and $f\in{\mathcal C}_b(K^n)$.
As the subspaces ${\mathcal C}^1_b(\R^{n\choose 2})$ of ${\mathcal C}_b(\R^{n\choose 2})$ and ${\mathcal D}(A)$ of ${\mathcal C}_b(K)$ are dense (also see Assumption~\ref{ass:004}), the subspace
\begin{equation}\label{e:separate}
\begin{aligned}
   &\big\{\big(g(m) + e^{-\lambda m(\smallx)} \big) \int\nu^{\smallx}(\mathrm{d}(\mr,\vk))\,\phi\big(\mr\big)f\big(\vk\big);\,\lambda\ge 0, \epsilon \geq 0, n\in\mathbb{N},\\
   &\hspace{1cm} g \in \mathcal{C}^3_b(\R_+) \mbox{ such that } \mathrm{supp}(g) \subset [0,\epsilon], g'(0)=\lambda, \ \phi\in{\mathcal C}^1_b(\R^{n\choose 2}),f\in{\mathcal C}_b(K^n)\mbox{ with }
   \\
   &\hspace{1cm}f(\kappa_1,...,\kappa_{l-1},\boldsymbol{\cdot},\kappa_{l+1},...,\kappa_n)\in{\mathcal D}(A),\,\forall l=1,...,n,\,\forall \kappa_1,...,\kappa_{l-1},\kappa_{l+1},...,\kappa_n\in K\big\}
   \subset\widetilde{\Pi}
   \end{aligned}
\end{equation}
also separates points in $\mathbb{M}^K$.

As the above subspace of $\widetilde{\Pi}$ contains constants, it vanishes nowhere. Moreover,
any function $F\in\Pi^{2,1,A}$ can be approximated by a sequence in $\widetilde{\Pi}$. As $\Pi^{2,1,A}$
is an algebra, the closure of $\widetilde{\Pi}$ clearly contains an algebra.

\end{proof}\smallskip

\begin{proposition}[Extension of the domain] If ${\mathcal X}$ is a $\mathbb{M}^K$-valued process which satisfies the $(\Omega,\widetilde{\Pi})$-martingale problem, then it satisfies the $(\Omega,{\mathcal D}(\Omega))$-martingale problem as well.
\label{P:008}
\end{proposition}

\begin{proof} Fix $F\in{\mathcal D}(\Omega)$.
We have seen that we can find a sequence $(F_n)$ in $\widetilde{\Pi}$ such that $F_n\to F$ as $n \rightarrow \infty$ and $\Omega F_n\to \Omega F$, as $n\to\infty$, boundedly pointwise. Indeed, if this is the case we can conclude from the fact that the process
\begin{equation}
\label{e:MFn}
   M^{F_n}:=\big(F_n({\mathcal X}_t)-F_n({\mathcal X}_0)-\int_0^t \Omega F_n({\mathcal X}_s)\mathrm{d}s\big)_{t\ge 0}
\end{equation}
is a martingale with respect to the canonical filtration of ${\mathcal X}$, that  the process
\begin{equation}
\label{e:MF}
\begin{aligned}
   M^{F}&:=\lim_{n\to\infty}\big(F_n({\mathcal X}_t)-F_n({\mathcal X}_0)-\int_0^t \Omega F_n({\mathcal X}_s)\mathrm{d}s\big)_{t\ge 0}
   \\
   &=\big(F({\mathcal X}_t)-F({\mathcal X}_0)-\int_0^t \Omega F({\mathcal X}_s)\mathrm{d}s\big)_{t\ge 0}
\end{aligned}
\end{equation}
is a martingale as well.
\end{proof}

% ======================================================================

{\sc Acknowledgement. } We thank Wolfgang L\"ohr for many fruitful discussions.

\bibliography{literatur}
\bibliographystyle{alpha}

%\appendix
%\section{Reminder}
%\begin{itemize}
%\item Do we need to assume that $\beta$, $\Gamma$, $\kappa$ (in case of mark function) are measurable, and if so, where do we assume it?
%\end{itemize}

% ======================================================================

\end{document}